\documentclass[a4paper,11pt,final]{amsart}

\usepackage[utf8]{inputenc}
\usepackage[T1]{fontenc}
\usepackage[english]{babel}
\usepackage[margin=2cm]{geometry}
\usepackage{amsmath}
\usepackage{amsthm}
\usepackage{amssymb}
\usepackage{dsfont}
\usepackage{mathrsfs}
\usepackage[all]{xy}
\usepackage{wasysym}
\usepackage{mathtools}
\usepackage{bbold}
\usepackage{empheq}
\usepackage{tensor}
\usepackage{bbm}
\usepackage{enumitem}
\usepackage{rotating}
\usepackage{multirow}
\usepackage{setspace}
\usepackage{fancybox}
\usepackage{xargs}
\usepackage{indentfirst}
\usepackage{xcolor}
\definecolor{InternalLinks}{rgb}{0.33, 0.29, 0.31}
\usepackage{hyperref}
\usepackage[final]{showkeys}
\usepackage[colorinlistoftodos,prependcaption,textsize=tiny]{todonotes}
\usepackage[nomargin,inline,marginclue,draft]{fixme}
\usepackage{comment}
\newif\ifshow 
\showtrue 
\ifshow
  \includecomment{wrap} 
\else
  \excludecomment{wrap} 
\fi

\theoremstyle{plain}
\newtheorem{theo}{Theorem}[section]
\newtheorem*{theo*}{Theorem}
\newtheorem{lemm}[theo]{Lemma}
\newtheorem{prop}[theo]{Proposition}
\newtheorem{coro}[theo]{Corollary}

\theoremstyle{definition}
\newtheorem{defi}[theo]{Definition}

\newtheorem{exam}{Example}[section]

\theoremstyle{remark}
\newtheorem{rema}[theo]{Remark}
\newtheorem*{rema*}{Remark}

\newtheorem*{nota*}{Notation}

\numberwithin{equation}{section}

\newcommand{\pr}{\noindent{\it Proof}\quad}
\newcommand{\fin}{\qed\smallskip}

\newcommand{\fell}[1]{\mathfrak{#1}}

\newcommand{\C}{\mathbb{C}}
\newcommand{\R}{\mathbb{R}}

\newcommand{\id}{\mathrm{id}}

\newcommand{\G}{\mathbb{G}}
\newcommand{\HH}{\mathbb{H}}
\newcommand{\dG}{\du{\mathbb{G}}}
\newcommand{\Go}{\G^{\ops}}
\newcommand{\Gc}{\G^{\con}}
\newcommand{\Gco}{\G^{\con,\ops}}
\newcommand{\dGc}{\dG^{\con}}
\newcommand{\dGo}{\dG^{\ops}}
\newcommand{\dGco}{\dG^{\con,\ops}}

\newcommand{\He}{\mathcal{H}}

\newcommand{\aYD}{\mathfrak{YD}}

\newcommand{\p}{\mathbf{p}}
\newcommand{\Pe}{\mathds{P}}
\newcommand{\pol}{\mathcal{O}}
\newcommand{\irr}{\mathrm{Irr}}

\newcommand{\pG}{\pol(\G)}

\newcommand{\dpG}{\du{\pG}}

\newcommand{\m}{\mathfrak{m}}
\newcommand{\M}{\mathrm{M}}
\newcommand{\du}[1]{\widehat{#1}}
\newcommand{\com}{\Delta}
\newcommand{\dcom}{\du{\Delta}}
\newcommand{\dS}{\du{S}}

\newcommand{\cou}{\varepsilon}
\newcommand{\dcou}{\du{\varepsilon}}
\newcommand{\uni}{\eta}
\newcommand{\var}{\varphi}
\newcommand{\dvar}{\du{\varphi}}

\newcommand{\op}{\textrm{op}}
\newcommand{\co}{\textrm{co}}
\newcommand{\ops}{\circ}
\newcommand{\con}{\mathrm{c}}

\newcommand{\alp}{\alpha}

\newcommand{\dtheta}{\du{\theta}}

\newcommand{\lbt}[2]{{#1}\indices*{^{#2}} \odo \tensor*[_{#2}]{#1}{}}
\newcommand{\rbt}[2]{{#1}\indices*{_{#2}} \odo \tensor*[^{#2}]{#1}{}}

\newcommand{\ltak}[2]{\mathbin{{#1}\indices*{^{#2}} \overline{\times} \tensor*[_{#2}]{#1}{}}}
\newcommand{\rtak}[2]{\mathbin{{#1}\indices*{_{#2}} \overline{\times} \tensor*[^{#2}]{#1}{}}}

\newcommand{\lctak}[2]{\mathbin{{#1}\indices*{^{#2}} \times \tensor*[_{#2}]{#1}{}}}
\newcommand{\rctak}[2]{\mathbin{{#1}\indices*{_{#2}} \times \tensor*[^{#2}]{#1}{}}}

\newcommand{\B}{\mathcal{B}}

\newcommand{\odo}{\otimes}
\newcommand{\oda}{\otimes}
\newcommand{\oti}{\otimes}

\newcommand{\sm}{\mathbin{\#}}

\newcommand{\blhd}{\blacktriangleleft}
\newcommand{\brhd}{\blacktriangleright}

\newcommand{\ome}{\omega}
\newcommand{\jo}{\mathrm{j}}
\newcommand{\djo}{\mathrm{\hat{j}}}
\newcommand{\ad}{\mathrm{Ad}}
\newcommand{\as}{$*$}

\newcommand{\tr}{\text{trv}}

\newcommand{\fle}[1]{\xymatrix@C=3pc@R=3pc{#1}}

\newcommand{\iso}{\cong}

\newcommand\restr[2]{{
  \left.\kern-\nulldelimiterspace 
  #1 
  \vphantom{\big|} 
  \right|_{#2} 
}}

\newcommand{\gpd}[3]{\xymatrix{\mathcal{{#1}} \ar@<0.5ex>[r]^-{{#2}}\ar@<-0.5ex>[r]_-{{#3}} & \mathcal{{#1}}^{(0)}}}

\definecolor{ToDo}{RGB}{30,144,255}
\definecolor{Question}{RGB}{220,20,60}
\definecolor{Attention}{RGB}{255,215,0}

\newcommandx{\td}[2][1=]{\hfill\todo[inline,size=\footnotesize,linecolor=ToDo,backgroundcolor=ToDo!40,bordercolor=black,#1]{{\bf ToDo:} #2}}
\newcommandx{\qs}[2][1=]{\hfill\todo[inline,size=\normalsize,linecolor=Question,backgroundcolor=Question!40,bordercolor=black,#1]{{\bf Question:} #2}}
\newcommandx{\rmq}[2][1=]{\hfill\todo[inline,size=\normalsize,linecolor=Attention,backgroundcolor=Attention!40,bordercolor=black,#1]{{\bf Remark:} #2}}

\makeatletter
\newcommand{\addresseshere}{\enddoc@text\let\enddoc@text\relax}
\makeatother
\makeatletter
\newcommand{\pushright}[1]{\ifmeasuring@#1\else\omit\hfill$\displaystyle#1$\fi\ignorespaces}
\newcommand{\pushleft}[1]{\ifmeasuring@#1\else\omit$\displaystyle#1$\hfill\fi\ignorespaces}
\makeatother

\title[Algebraic quantum transformation groupoids of compact type]{Algebraic quantum transformation groupoids \\ of compact type}
\author{Frank Taipe}
\address{Université Paris-Saclay, CNRS, Laboratoire de Mathématiques d'Orsay, 91405 Orsay, France}
\email{frank.taipe@universite-paris-saclay.fr}


\hypersetup{
pdfauthor={Frank Taipe},
pdftitle={Algebraic quantum transformation groupoids of compact type},
pdfsubject={Quantum version of transformation groupoids},
pdfkeywords={Quantum groupoids, Transformation groupoids}
}

\begin{document}

\begin{abstract}
In this work, we introduce a class of Timmermann's measured multiplier Hopf \as-algebroids called {\em algebraic quantum transformation groupoids of compact type}. Each object in this class admits a Pontrjagin-like dual called an {\em algebraic quantum transformation groupoid of discrete type}. This compact/discrete duality in the framework of algebraic quantum transformation groupoids recover the one between compact and discrete Van Daele's algebraic quantum groups. Among the non-trivial examples of algebraic quantum transformation groupoids of compact type, we give constructions arising from Fell bundles and quantum quotient spaces.
\end{abstract}

\maketitle

\section{Introduction}

Measured quantum groupoids were introduced by F. Lesieur in \cite{LPHD03} using as main ingredients two operator algebraic structures introduced by J.-M. Vallin, namely Hopf-von Neumann bimodules \cite{V96} and pseudo-multiplicative unitaries \cite{V00}. Two important features of these general quantum objects is that, on one hand, they generalize von Neumann locally compact quantum groups, in the sense of J. Kustermans and S. Vaes \cite{KV03}, and on the other hand, they are natural objects arising as quantum symmetries of general inclusions of von Neumann algebras \cite{EV00,EN00,EN05}. This last feature generalize the fact that locally compact quantum groups arise as quantum symmetries of irreducible depth 2 subfactor inclusions \cite{EN96,E98}, and also the fact that finite quantum groupoids arises as quantum symmetries of finite index depth 2 II${}_{1}$-subfactor inclusions \cite{NV00}. A nice result in the relationship between finite quantum groupoids and II${}_{1}$-subfactor inclusions is the connection of these objects via a C*-categorical framework \cite{NV00_2}. Therefore, naturally, two questions arise
\begin{enumerate}[label=\textup{(Q\arabic*)}]
\item\label{q:1} In the spirit of locally compact quantum groups, is there a nice equivalent C*-algebraic counterpart of a measured quantum groupoid?
\item\label{q:2} In the spirit of finite quantum groupoids, is it possible to make a link between measured quantum groupoids of compact type and C*-categorical structures arising from inclusions of von Neumann algebras?
\end{enumerate}

In \cite{TPHD18}, motivated by the question~\ref{q:1}, we gave the first steps of a C*-algebraic counterpart for the measured quantum transformation groupoids introduced by M. Enock and T. Timmermann \cite{ET16}. This class of measured quantum groupoids are build using braided commutative Yetter--Drinfeld structures on actions of von Neumann locally compact quantum groups. The use of Yetter--Drinfeld structures for the mentioned construction was motivated by the pioneer work of J.-H. Lu, who also coined the term ``quantum transformation groupoid'', in \cite{L96}. A nice aspect of the category of measured quantum transformation groupoids is that it is closed under a Pontrjagin-like duality following the philosophy of G. Kac. In \cite{Ta22_3}, we introduce a general definition of an algebraic quantum transformation groupoid and we show that using a GNS-type construction, we obtain examples of locally compact quantum groupoids in the sense of M. Enock \cite{EN19} which are not necessarily locally compact quantum groupoids of separable type in the sense of B.J. Kahng and A. Van Daele \cite{KVD18, KVD19}. At the moment of this paper is wrote, an Enock's locally compact quantum groupoid seems to be a good candidate to be the most general definition for a C*-algebraic counterpart of a measured quantum groupoid.

The purpose of this work is two fold, first we want to construct explicitly the Pontrjagin dual of a measured Lu's Hopf \as-algebroid associated with a braided commutative Yetter--Drinfeld structure over a no necessarily finite Hopf \as-algebra, and second we want to give the first steps on the question~\ref{q:2}. Therefore, we aim to introduce a nice class of compact quantum groupoids, the class of algebraic quantum transformation groupoids of compact type. The essential ingredient to be considered is a unital braided commutative measured Yetter--Drinfeld \as-algebra over an algebraic quantum group of compact type \cite{Ta22_1}. Despite the fact that in the present work, we will only deal with purely algebraic objects and constructions, we must keep in mind that the main source of unital braided commutative Yetter--Drinfeld \as-algebras over algebraic quantum groups needed for the construction of algebraic quantum transformation groupoids of compact type are indeed {\em algebraic cores} of Yetter--Drinfeld C*-algebras over compact quantum groups \cite{NV09,FTW23}. That means that the \as-algebras we work with are actually operator algebras. The main difference between those algebraic quantum transformation groupoids of compact type and the general ones, is that in the general case we need to be more careful dealing with non-unital Yetter--Drinfeld \as-algebras and with general algebraic quantum groups. In a further work, we will introduce and study the representation theory of algebraic quantum transformation groupoids of compact type and its connection with inclusions of von Neumann algebras.

\smallskip

This paper is organized as follows: Section~\ref{sec:preliminaries} is focus to give the most important preliminaries about Van Daele's algebraic quantum groups, and Timmermann's measured multiplier Hopf \as-algebroids. Section~\ref{sec:bcyd}, is intend to recollect only important results about braided commutative Yetter--Drinfeld \as-algebras over an algebraic quantum group necessary for the present work. Here, we give an important result concerning unital braided commutative Yetter--Drinfeld \as-algebras over algebraic quantum groups of compact type. In section~\ref{sec:mhad_yd}, we present the construction of a measured multiplier Hopf \as-algebroid associated with a unital braided commutative measured Yetter--Drinfeld \as-algebra over an algebraic quantum group of compact type. We finish the section with the construction of its Pontrjagin-like dual measured multiplier Hopf \as-algebroid. In this section, we exhibit examples arising mainly from classic transformation groupoids, Fell bundles and quantum quotient spaces.

\subsection*{Notations and conventions}

An algebra in this work means a not necessary unital algebra over the field of complex numbers $\C$. We only work with idempotent algebras which not have a unit necessarily but its multiplication map is always assumed to be non-degenerate. When we are dealing with multiplier Hopf \as--algebras or multiplier Hopf \as-algebroids, we use sometimes the Sweedler type leg notation through the {\em Van Daele's covering technique}.

Keeping that in mind, for a given non-degenerate idempotent algebra $A$, we can consider its multiplier algebra $\M(A)$ (see below for the definition), which is a unital idempotent algebra. We denote by ${}_{A}M$ a left $A$-module $M$ and by $\rhd : A \odo M \to M$ its module map. The multiplication map on an algebra $A$ is always denoted by $\m_{A}: A \odo A \to A$. We write ${}_{A}A$ (or $A_A$) when we regard $A$ as a left (or right) module over itself with respect to the left (or right) module map $\m_{A}$, respectively. For any algebra $A$, its opposite algebra, denoted by $A^{\op}$, is the vector space $A$ endowed with the multiplication map given by $\m_{A^{\op}}: A^{\op} \odo A^{\op} \to A^{\op}$, $a^{\op} \odo a'^{\op} \mapsto (a'a)^{\op}$. We will use always the superscript ${}^{\op}$ on an algebra element to indicate that the element is being regarded as an element of the opposite algebra. Sometimes, the superscript  ${}^{\op}$ also denotes the trivial anti-isomorphisms $\id_{A}: A^{\op} \to A$ or $\id_{A}: A \to A^{\op}$. 

\subsection*{Acknowledgments}

The author would like to express his thanks to Leonid Vainerman for the very fruitful discussions.

\section{Preliminaries}\label{sec:preliminaries}

\subsection{Transformation groupoids}

Since the present work is concerned with the notion of a quantum transformation groupoid, in this paragraph we recall the ``classical'' notion of a transformation groupoid or also known in the literature as an action groupoid. Let $\cdot : G \times X \to X$ be a left action of a (possible infinite) group $G$ on a (possible infinite) set $X$. The {\em transformation groupoid associated with the action of $G$ on $X$}, denoted by $G \ltimes X$, is the Cartesian product set $G \times X$ endowed with the applications
\[
\begin{array}{rccc}
\cdot : & (G \rtimes X)^{(2)} & \to & G \ltimes X \\
& ((g,x),(g',x')) & \mapsto & (gg',x')
\end{array}
, \quad
\begin{array}{rccc}
\phantom{\cdot}^{-1} : & G \ltimes X & \to & G \ltimes X \\
& (g,x) & \mapsto & (g^{-1}, g \cdot x)
\end{array},
\]
where $(G \ltimes X)^{(2)} := \left\{((g,x),(g',x')) \in (G \ltimes X)\times(G \ltimes X) \,:\, x = g' \cdot x' \right\}$ is the set of composable pairs of the transformation groupoid. The domain map and the rang map of $G \ltimes X$ are respectively given by
\[
\begin{array}{rccc}
d : & G \ltimes X & \to & (G \ltimes X)^{(0)} \\
& (g,x) & \mapsto & x
\end{array}
\quad\text { and }\quad
\begin{array}{rccc}
r : & G \ltimes X & \to & (G \ltimes X)^{(0)} \\
& (g,x) & \mapsto & g \cdot x
\end{array},
\]
where $(G \ltimes X)^{(0)} := \{e\} \times X \simeq X$ is the unit space of the transformation groupoid.

In analogy with the quantization of the notion of a group, we need to keep in mind the two canonical \as-algebras associated with the transformation groupoid $\G \ltimes X$:
\begin{enumerate}[label=\textup{(TG\arabic*)}]
\item the \as-algebra of finite supported complex-valued functions on $G \ltimes X$, denoted by $K(G \ltimes X)$, with structure given by elementwise product of functions and involution $f^{*}(g,x) = f(g^{-1},g \cdot x)$ for all $f \in K(G \ltimes X)$, $g \in G$, $x \in X$; and
\item the groupoid \as-algebra of $G \ltimes X$, denoted by $\C[G \ltimes X]$, which is the \as-algebra generated by formal elements $\{\lambda_{(g,x)}\}_{(g,x) \in G \ltimes X}$ and relations given by $\lambda_{(g,x)}\lambda_{(g',x')} = \delta_{x,g'\cdot x'}\lambda_{(gg',x')}$ and $\lambda^{*}_{(g,x)} = \lambda_{(g,x)^{-1}} = \lambda_{(x \cdot g,g^{-1})}$, for all $g,g' \in G$ and $x,x' \in X$. 
\end{enumerate}

\subsection{Algebras and its modules}

Let $A$ be an algebra over the field of complex number $\C$. A left module ${}_{A}M$ is called {\em faithful}, if for each non-zero $a \in A$ there exists $m \in M$ such that $a \rhd m$ is non-zero; {\em non-degenerate}, if for each non-zero $m \in M$ there exits $a \in A$ such that $a \rhd m$ is non-zero; {\em unital}, if $A \rhd M := \text{span}\{a \rhd m \,:\, a \in A, m \in M\} = M$. We say that ${}_{A}M$ {\em admits local units in $A$} if for every finite subset $F \subset M$ there exists $a \in A$ with $a \rhd m = m$ for all $m \in F$. For right modules, we use similar notations and terminology by identifying right $A$-modules with left $A^{\op}$-modules.

\begin{rema}\label{rem:sits_left_right}
If ${}_{A}M$ admits local units, then it is non-degenerate and unital. The module ${}_{A}A$ (or $A{}_{A}$) is non-degenerate if and only if $A{}_{A}$ (or ${}_{A}A$) is faithful, respectively. The right module $A{}_{A}$ is non-degenerate if and only if the natural map $A \to \textrm{End}(A_A)$, $a \mapsto [L_{a} : b \mapsto ab]$ is injective. Similarly, the left module ${}_{A}A$ is non-degenerate if and only if the natural map $A \to \textrm{End}(\tensor[_{A}]A{})^{\op}$, $a \mapsto [R_{a} : b \mapsto ba]$ is injective.
\end{rema}

\begin{rema}
Given a left-module ${}_{A}M$ and $T \in \mathrm{End}({}_{A}M)^{\op}$, we write $(m)T$ for the image of an element $m$ under $T$. So, in this way we have $(m)(ST) = ((m)S)T$ for all $m \in M$ and $S,T \in \mathrm{End}({}_{A}M)^{\op}$.
\end{rema}

An algebra $A$ is called {\em involutive or \as-algebra} if there exists an anti-linear map $\phantom{\cdot}^{*} : A \to A$ called {\em the involution of $A$} such that $(a^{*})^{*} = a$ and $(ab)^{*} = b^{*}a^{*}$ for all $a,b \in A$; {\em unital} if there is an element denoted by $1_{A}$ such that $a = 1_{A}a = a1_{A}$; {\em non-degenerate} if the modules ${}_{A}A$, $A_A$ both are non-degenerate; {\em idempotent}, if one of the modules ${}_{A}A$, $A_{A}$ is unital. We say that $A$ {\em admits local units} if the modules ${}_{A}A$, $A_A$ both admit local units in $A$.

\begin{rema}\label{rem:unital}
If $A$ is a unital algebra, then automatically $A$ is non-degenerate, idempotent and admits local units.
\end{rema}

A functional $\varphi$ on an algebra $A$ is called faithful if $\varphi(ab) = 0$ and $\varphi(ba) = 0$ for all $b \in A$ implies $a =0$. In case $A$ is a \as-algebra, $\varphi$ is called {\em self-adjoint} if $\varphi(a^{*}) = \overline{\varphi(a)}$ for all $a \in A$ and {\em positive} if additionally we have $\varphi(aa^{*}) \geq 0$ for all $a \in A$. We say that $\varphi$ satisfies a {\em weak KMS property}, if there is an automorphism $\sigma$ on $A$ such that $\varphi(ab) = \varphi(b\sigma(a))$ for all $a,b \in A$. If $\varphi$ is faithful and satisfies a weak KMS property, thus we say that $\varphi$ admits a {\em modular automorphism} which is denoted by $\sigma^{\varphi}$.

Let $A$ be an \as-algebra and $\gamma$ be an algebra automorphism on $A$ such that $\gamma \circ * \circ \gamma \circ * = \id_{A}$. The anti-linear map $A^{\op} \to A^{\op}$, $a^{\op} \mapsto \gamma(a^{*})^{\op}$, defines an involution on $A^{\op}$. The algebra $A^{\op}$ endowed with the involution defined above will be denoted by $A^{\op}_{\gamma}$ and it will be called {\em the $\gamma$-opposite \as-algebra of $A$}. In the case $\gamma = \id_{A}$, the $\gamma$-opposite \as-algebra $A^{\op}_{\gamma}$ will be called {\em the canonical opposite \as-algebra} and it will be denoted simply by $A^{\op}$. If $\hat\gamma$ denotes the inverse of $\gamma$, the map $\hat\gamma^{\ops}: A^{\op}_{\gamma} \to A^{\op}_{\gamma}$, defined by $\hat\gamma^{\ops}(a^{\op}) = \hat\gamma(a)^{\op}$ for all $a \in A$, yields an algebra automorphism on $A^{\op}_{\gamma}$, with inverse $\gamma^{\ops}:a^{\op} \mapsto \gamma(a)^{\op}$, such that $\hat\gamma^{\ops} \circ * \circ \hat\gamma^{\ops} \circ * = \id_{A^{\op}_{\gamma}}$. Observe that with those notations, we have $(A^{\op}_{\gamma})^{\op}_{\hat\gamma^{\ops}} = A$ as \as-algebras.

Given a functional $\mu$ on $A$, the opposite functional $\mu^{\ops}:A^{\op}_{\gamma} \to \C$ is the linear functional defined by $\mu^{\ops}(a^{\op}) = \mu(a)$ for all $a \in A$. The functional $\mu$ is faithful if and only if the functional $\mu^{\ops}$ is faithful. The functional $\mu$ is self-adjoint and $\mu\circ\gamma = \mu$ if and only if the functional $\mu^{\ops}$ is self-adjoint. The functional $\mu$ is positive and it satisfies a weak KMS property with automorphism $\gamma^{-1}$ if and only if the functional $\mu^{\ops}$ is positive and it satisfies a weak KMS property with automorphism $(\gamma^{-1})^{\op}$.

Let $A$ be a non-degenerate \as-algebra, not necessarily unital. The {\em multiplier algebra of $A$} is the vector space
\[
\M(A) := \{ (L,R) \in \textrm{End}(A_A) \times \textrm{End}(\tensor[_{A}]A{})^{\op} \; : \; R(a)b = aL(b)\,\text{ for all }a,b\in A \}.
\]
endowed with multiplication given by $(L,R)\cdot(L',R') := (L L',R' R)$, for all $(L,R),(L',R')\in \M(A)$ and involution given by $(L,R)^{*} := (R^{\dagger},L^{\dagger})$ for all $(L,R) \in \M(A)$. Here, $R^{\dagger}(a) := R(a^{*})^{*}$ and $L^{\dagger}(a) := L(a^{*})^{*}$ for all $a\in A$.

\begin{rema}\label{rem:sits}
It follows from Remark~\ref{rem:sits_left_right} that $A$ sits in $\M(A)$ as a sub-\as-algebra with the identification $a \mapsto (L_{a},R_{a})$ for all $a \in A$. In particular, if $A$ admits a unity, then $\M(A) = A$.
\end{rema}

Given an element $T=(L,R) \in \M(A)$ and $a \in A$, we use the notations $[T]a := L(a)$ and $a[T] := R(a)$. Using the identification mentioned in Remark\ref{rem:sits}, we have $Ta = [T]a$ and $aT = a[T]$. A homomorphism of algebras $f: A \to M(B)$ is called {\em non-degenerate} if $f(A)B = B$ and $B f(A) = B$. In that case, the homomorphism $f: A \to \M(B)$ can be extended to $\M(A)$ and this extension, denoted also by $f$, is defined by $[f(T)]f(a)b := f([T]a)b$ and $bf(a)[f(T)] := bf(a[T])$ for every $T \in \M(A)$, $a \in A$ and $b \in B$.

\subsection*{The left and right extended Takeuchi product}

Let $A$ be a non-degenerate algebra, $B$ and $C$ be two non-degenerate algebras related with $A$ via two non-degenerate faithful embeddings $\iota_{B}:B \to \M(A)$ and $\iota_{C}:C \to \M(A)$ whose images commute with each other, and $t_{B}:B \to C$ and $t_{C}:C \to B$ be two anti-isomorphisms.

We write $A_{B}$ and ${}^{B}A$ if we regard $A$ as a right or left $B$-module via
\[
a \lhd_{B} x := a\iota_{B}(x), \qquad \text{and} \qquad x {}^{B}\rhd a := a\iota_{C}(t_{B}(x))
\]
for all $a \in A$ and $x \in B$, respectively. We also write $A^{C}$ and ${}_{C}A$ if we regard $A$ as a right or left $C$-modules via
\[
a \lhd^{C} y := \iota_{B}(t_{C}(y))a, \qquad \text{and} \qquad y {}_{C}\rhd a := \iota_{C}(y)a
\]
for all $a \in A$ and $y \in C$, respectively. With those notations, in an obvious way the balanced tensor product $\rbt{A}{B}$ become a right $\M(A) \odo \M(A)$-module and the balanced tensor product $\lbt{A}{C}$ become a left $\M(A) \odo \M(A)$-module.

{\em The extended right Takeuchi product of $A$ over $B$}, denoted by $A \rtak{}{B} A$, is the vector subspace of elements $T \in \mathrm{End}(\rbt{A}{B})^{\op}$ such that for all $a,b \in A$, there exist elements $(a \odo 1)T, (1 \odo b)T \in \rbt{A}{B}$ verifying
\[
(a \odo b)T = (1 \odo b)((a \odo 1)T) = (a \odo 1)((1 \odo b)T).
\]
 If the left module ${}_{A}{A}$ is idempotent and non-degenerate, the modules $A_{B}$, ${}^{B}A$ are faithful and idempotent and the space $\rbt{A}{B}$ is non-degenerate as a left module over $A \odo 1$ and over $1 \odo A$, then the extended right Takeuchi product $A \rtak{}{B} A$ becomes a on-degenerate algebra. Moreover, any element $f \in A \rtak{}{B} A$ commutes with the left $A \odo A$-module action, i.e. we have $(ab \odo a'b')f = (a \odo a')((b \odo b')f)$ for all $a,a',b,b' \in A$. In the case $A$ is unital, then the map from $A \rtak{}{B} A$ to $\rbt{A}{B}$ given by $T \mapsto (1_{A} \odo 1_{A})T$ identifies $A \rtak{}{B} A$ with the classic right Takeuchi product $A \rctak{}{B} A$ which is the algebra defined by
\[
A \rctak{}{B} A = \{a \odo b \in \rbt{A}{B} : \iota_{B}(y)a \odo b = a \odo \iota_{C}(t_{B}(y))b \text{ for all } y \in B \}.
\]

Similarly, {\em the extended left Takeuchi product of $A$ over $C$}, denoted by $A \ltak{}{C} A$, is the vector subspace of elements $T \in \mathrm{End}(\lbt{A}{C})$ such that for all $a,b \in A$, there exist elements $T(a \odo 1), T(1 \odo b) \in \lbt{A}{C}$ verifying
\[
T(a \odo b) = (T(a \odo 1))(1 \odo b) = (T(1 \odo b))(a \odo 1).
\]
If the right module $A_{A}$ is idempotent and non-degenerate, the modules $A^{C}$, ${}_{C}A$ are faithful and idempotent and the space $\lbt{A}{C}$ is non-degenerate as a right module over $A \odo 1$ and over $1 \odo A$, then the extended left Takeuchi product $A \ltak{}{C} A$ becomes a non-degenerate algebra. Moreover, any element $f \in A \ltak{}{C} A$ commutes with the right $A \odo A$-module action, i.e. we have $f(ab \odo a'b') = (f(a \odo a'))(b \odo b')$ for all $a,a',b,b' \in A$. In the case $A$ is unital, then the map from $A \ltak{}{C} A$ to $\lbt{A}{C}$ given by $T \mapsto T(1_{A} \odo 1_{A})$ identifies $A \ltak{}{C} A$ with the classic left Takeuchi product $A \lctak{}{C} A$ which is the algebra defined by
\[
A \lctak{}{C} A = \{a \odo b \in \lbt{A}{C} : a\iota_{B}(t_{C}(x)) \odo b = a \odo b\iota_{C}(x) \text{ for all } x \in C \}.
\]

In the involutive case, if we additionally assume $t_{C}\circ * \circ t_{B} \circ * = \id_{B}$ or equivalently $t_{B}\circ * \circ t_{C} \circ * = \id_{C}$, then the anti-linear bijective map
\[
(-)^{*} \odo (-)^{*}: \rbt{A}{B} \to \lbt{A}{C},\; a \odo b \mapsto a^{*} \odo b^{*}
\]
is well defined and we can consider the anti-linear anti-isomorphism
\[
(* \overline{\times} *): \mathrm{End}(\rbt{A}{B}) \to \mathrm{End}(\lbt{A}{C}),\; T \mapsto T^{(* \overline{\times} *)} := ((-)^{*} \odo (-)^{*})\circ T\circ((-)^{*} \odo (-)^{*})^{-1}
\]
which satisfies
\[
(\rtak{A}{B})^{(* \overline{\times} *)} = \ltak{A}{C}.
\]

\subsection{Algebraic quantum groups and their actions}

We collect here some basics about algebraic quantum groups, for further reading we send the interested reader to the original papers \cite{VD94, VD98} or a more recently series of papers \cite{VD23_1,VD23_2,VD23_3}.

This work is concerned with algebraic quantum groups in the following sense,

\begin{defi}\label{def:aqg}
An {\em algebraic quantum group} is a triplet $\G = (\pG,\com_{\G},\varphi_{\G})$ which consists of a non-degenerate \as-algebra $\pG$, a \as-homomorphism $\com_{\G}:\pG \to \M(\pG \odo \pG)$, called {\em the multiplication of $\G$}, satisfying the conditions
\begin{enumerate}
\item $\com(h)(1 \odo g)$ and $(h \odo 1)\com(g)$ belong to $\pG \odo \pG$ for all $h,g \in \pG$;
\item $T_{\lambda} : h \odo g \mapsto \com(h)(1 \odo g), T_{\rho} : h \odo g \mapsto (h \odo 1)\com(g)$ are bijective maps;
\item the {\em associative property of $\G$} holds
\[
(h \odo 1 \odo 1)(\com \odo \id)(\com(g)(1 \odo h')) = (\id \odo \com)((h \odo 1)\com(g))(1 \odo 1 \odo h')
\]
for all $h,h',g \in \pG$;
\end{enumerate}
and a non-zero positive faithful functional $\varphi_{\G}$ on $\pG$, called the {\em left integral of $\G$}, which satisfies the left invariant condition
\[
(\id \odo \var)(\com(h)(g \odo 1)) = \var(h)g = (\id \odo \var)((g \odo 1)\com(h))
\]
for every $h,b \in \pG$.
\end{defi}

Given an algebraic quantum group $\G$, there exists a unique \as-homomorphism $\cou_{\G}: \pG \to \C$, named {\it the counit of $\G$}, and a unique anti-isomorphism $S_{\G}:\pG \to \pG$ such that $S_{\G}\circ * \circ S_{\G} \circ * = \id_{\pG}$, named {\it the antipode of $\pG$}, satisfying the conditions
\begin{enumerate}
\item $(\cou_{\G} \odo \id)(\com(h)(1 \odo g)) = hg = (\id \odo \cou_{\G})((h \odo 1)\com(g))$;
\item $\m(S_{\G} \odo \id)(\com(h)(1 \odo g)) = \cou_{\G}(h)g = \m(\id \odo S_{\G})((g \odo 1)\com(h))$;
\end{enumerate}
for all $h,g \in \pG$. The pair $(\pG,\com_{\G})$ is called the {\em underlying multiplier Hopf \as-algebra of $\G$}. If we have that the antipode $S_{\G}$ is *-preserving, or equivalently if $S^{2}_{\G} = \id_{\pG}$, then we said that {\em the algebraic quantum group $\G$ is of Kac type}.

The functional $\psi_{\G} := \varphi_{\G}\circ S_{\G}$, called the {\em right integral of $\G$}, is a non-zero positive faithful functional on $\pG$ which satisfies the right invariant condition
\[
(\psi \odo \id)(\com(h)(g \odo 1)) = \psi(h)g = (\psi \odo \id)((g \odo 1)\com(h))
\]
for every $h,g \in \pG$.

\begin{rema}\label{rem:left_right_invariance}
By \cite[Proposition~3.11]{VD98}, the left an right integral of $\G$ also satisfies the conditions
\[
S_{\G}((\id \odo \varphi_{\G})(\com_{\G}(h)(1_{\pG} \odo g))) = (\id \odo \varphi_{\G})((1_{\pG} \odo h)\com_{\G}(g))
\]
and
\[
S_{\G}((\psi_{\G} \odo \id)((h \odo 1_{\pG})\com_{\G}(g))) = (\psi_{\G} \odo \id)(\com_{\G}(h)(g \odo 1_{\pG}))
\]
for all $h,g \in \pG$, respectively.
\end{rema}

There is a unique automorphism $\sigma_{\G}$ on $\pG$, called {\em the modular automorphism of the left integral $\var_{\G}$}, such that $\varphi_{\G}(hg) = \varphi_{\G}(g\sigma_{\G}(h))$ for all $h,g \in \pG$. Additionally, the map $\sigma_{\G}$ satisfies $\varphi_{\G}\circ\sigma_{\G} = \varphi_{\G}$ and $\com_{\G}\circ\sigma_{\G} = (S^{2}_{\G} \odo \sigma_{\G})\circ\com_{\G}$. Similarly, there is a unique automorphism $\sigma'_{\G}$ on $\pG$, called {\em the modular automorphism of the right integral $\psi_{\G}$}, such that $\psi_{\G}(hg) = \psi_{\G}(g\sigma'_{\G}(h))$ for all $h,g \in \pG$. Additionally, the map $\sigma'_{\G}$ satisfies $\psi_{\G}\circ\sigma'_{\G} = \psi_{\G}$, $\com_{\G}\circ\sigma'_{\G} = (\sigma'_{\G} \odo S^{-2}_{\G})\circ\com_{\G}$ and $\sigma'_{\G} = S^{-1}_{\G}\circ\sigma^{-1}_{\G}\circ S_{\G}$.

There is an invertible element $\delta \in \M(\pG)$ such that $\psi_{\G} = \delta\cdot\varphi_{\G} = \varphi_{\G}\cdot\delta$. We also have the relations $\com_{\G}(\delta) = \delta \odo \delta$, $\cou_{\G}(\delta) = 1$, $S_{\G}(\delta) = \delta^{-1}$, $\sigma_{\G}(\delta) = \sigma'_{\G}(\delta) = \delta$ and $\sigma'_{\G}(h) = \delta\sigma_{\G}(h)\delta^{-1}$ for all $h \in \pG$.

For an element $h \in \pG$, consider the notations $h\cdot\varphi_{\G}$ and $\varphi_{\G}\cdot h$ for the functionals on $\pG$ given by $g \mapsto \varphi_{\G}(gh)$ and $g \mapsto \varphi_{\G}(hg)$, respectively. With the notations above, we have $\varphi_{\G} \cdot h = \sigma_{\G}(h)\cdot\varphi_{\G}$ for all $h \in \G$, then {\em the dual space of $\pG$} is defined as the vector subspace
\[
\dpG := \{h\cdot\varphi_{\G} \in \pG^{\vee} : h \in \pG \} = \{\varphi_{\G}\cdot h \in \pG^{\vee} : h \in \pG \}.
\]

It was shown in \cite{VD94}, that there are a structure of non-degenerate \as-algebra on $\dpG$, a \as-homomorphism $\dcom_{\G}:\dpG \to \M(\dpG \odo \dpG)$ and a non-zero positive faithful functional $\dvar_{\G}$ on $\dpG$ such that the triplet $(\dpG,\dcom_{\G},\dvar_{\G})$ yields an algebraic quantum group. In order to give the aforementioned structure in a simplified way, we prefer to make use of the notion of pairing. The {\em canonical pairing of an algebraic quantum group $\G$} is the bilinear map defined by
\[
\begin{array}{lccc}
\p: & \pG \times \dpG & \to & \C \\
& (h, \ome) & \mapsto & \ome(h)
\end{array}.
\]
With this canonical pairing, the structure of the algebraic quantum group $(\dpG,\dcom_{\G},\dvar_{\G})$ is given via the conditions
\begin{enumerate}[label=\textup{(P\arabic*)}]
\item\label{cond:P1} $\p(h,\ome\cdot\ome') = \p^{2}(\com_{\G}(h),\ome \odo \ome')$, $\p(h,\ome^{*}) = \p(S_{\G}(h^{*}),\ome)$,
\item\label{cond:P2} $\p^{2}((h \odo g),\dcom_{\G}(\ome)) = \p(hg,\ome)$,
\item\label{cond:P3} $\p(h,\dS_{\G}(\ome)) = \p(S_{\G}(h),\ome)$, $\p(h,\varphi_{\G}) = \dcou_{\G}(h\cdot\varphi_{\G})$
\end{enumerate}
for all $h,g \in \pG$, $\ome,\ome' \in \dpG$ and $h \in \pG$. Moreover, the left integral is given by $\dvar_{\G} := \du{\psi}_{\G}\circ\dS^{-1}_{\G}$, where $\du{\psi}_{\G}$ is the functional given by $h\cdot\varphi_{\G} \mapsto \cou_{\G}(h)$. In condition~\ref{cond:P2}, $\p^{2}$ denotes the pairing
\[
(h \odo g,\ome \odo \ome') \in (\pG \odo \pG) \times (\dpG \odo \dpG) \mapsto \p(h,\ome)\p(g,\ome').
\]
Associated with the canonical pairing of an algebraic quantum group $\G$, there is a unique unitary element $U \in \M(\pG \odo \dpG)$, called {\em the algebraic multiplicative unitary $\G$}, satisfying the conditions
\begin{enumerate}
\item $(\com_{\G} \odo \id)(U) = U_{13}U_{23}$,
\item $(\id \odo \Sigma\circ\dcom_{\G})(U) = U_{13}U_{12}$.
\end{enumerate}

Given an algebraic quantum group $\G$, the {\em dual algebraic quantum group of $\G$} is given by the triplet
\[
\dG = (\dpG,\dcom^{\co}_{\G}:=\Sigma\circ\dcom_{\G},\du{\psi}_{\G}),
\]
the {\em opposite algebraic quantum group of $\G$} is given by the triplet
\[
\Go = (\pG,\com^{\co}_{\G}:=\Sigma\circ\com_{\G},\psi_{\G}),
\]
and the {\em conjugate algebraic quantum group of $\G$} is given by the triplet
\[
\Gc = (\pG^{\op},\com^{\op}_{\G}:=({\;}^{\op} \odo {\;}^{\op})\circ\com_{\G}\circ{\;}^{\op}),\varphi^{\ops}_{\G}:=\varphi_{\G}\circ{\;}^{\op}),
\]
where the involution on the opposite algebra $\pG^{\op}$ is defined by $(h^{\op})^{*} = S^{-2}_{\G}(h^{*})^{\op}$ for all $h \in \pG$ (see Remark~\ref{rem:opposite_involution}). With those notations, it hold the following equalities $\G = (\Go)^{\ops} = (\Gc)^{\con}$, $\du{(\Go)} = (\dG)^{\con}$ and $\du{(\Gc)} = (\dG)^{\ops}$. Moreover, if we consider the algebraic quantum group $\G^{\ops,\con} := (\Go)^{\con} = (\Gc)^{\ops}$, the linear map $S_{\G}: \G^{\ops,\con} \overset{}{\to} \G$ yields an isomorphism of algebraic quantum groups.

An algebraic quantum group $\G$ is called of compact type if $\pG$ is a unital \as-algebra doing the underlying multiplier Hopf \as-algebra of $\G$,  $(\pG,\com_{\G})$, a unital Hopf \as-algebra. In this case, the left integral $\varphi_{\G}$ is also a right integral, i.e. $\psi_{\G} = \varphi_{\G}\circ S_{\G} = \varphi_{\G}$, and the modular automorphism $\sigma_{\G}$ satisfies $S^{-1}_{\G} \circ \sigma^{-1}_{\G} = \sigma_{\G} \circ S^{-1}_{\G}$.

\begin{exam}\label{example_group}
Let $G$ be a non necessarily finite group, $K(G)$ be the \as-algebra of complex valued finitely supported functions on $G$ and $\C[G]$ be the group \as-algebra of $G$, i.e. the unital \as-algebra generated by formal elements $\{\lambda_{g}\}_{g \in G}$ with relations given by $\lambda_{g}\lambda_{h} = \lambda_{gh}$ and $(\lambda_{g})^* = \lambda_{g^{-1}}$ for all $g,h \in G$.

In this case, $\M(K(G))$ is the \as-algebra of all complex valued functions on $G$, $K(G) \odo K(G)$ can be naturally identified with the \as-algebra of complex valued finitely supported functions on $G \times G$, and $\M(K(G) \odo K(G))$ can be identified with the \as-algebra of all complex-valued functions on $G \times G$. Take the \as-homomorphism $\com_{G}: K(G) \to \M(K(G) \odo K(G))$ defined by $p \mapsto [\com_{G}(p):(g,h) \mapsto p(gh)]$, for all $p \in K(G)$ and the non-zero positive functional $\var_{G}: K(G) \to \C$, $p \mapsto \sum_{g \in G}p(g)$, then $\G = (K(G),\com_{G},\varphi_{G})$ yields an algebraic quantum group. Its underlying multiplier Hopf \as-algebra $(K(G),\com_{G})$ is commutative in the sense that the \as-algebra $K(G)$ is commutative, and its counit and antipode are given by
\[
\begin{array}{rccc}
\cou_{G}: & K(G) & \to & \C \\
& p & \mapsto & p(e)
\end{array}
,\qquad
\begin{array}{rccc}
S_{G}: & K(G) & \to & K(G) \\
& p & \mapsto & [S(p): g \mapsto p(g^{-1})]
\end{array},
\]
respectively. On the other hand, since for each $x \in \C[G]$, there exists a unique complex valued finitely supported function $p_{x} : G \to \C$, such that $x = \sum_{g \in G}p_{x}(g)\lambda_{g}$, we can consider the \as-homomorphism $\dcom_{G}: \C[G] \to \C[G] \odo \C[G]$, $\lambda_{g} \mapsto \lambda_{g} \odo \lambda_{g}$ and the non-zero positive functional $\dvar_{G}: \C[G] \to \C$, $\lambda_{g} \mapsto \delta_{g,e}$. Thus, the triplet $(\C[G],\dcom_{G},\dvar_{G})$ yields an algebraic quantum group. Its underlying multiplier Hopf \as-algebra is cocommutative in the sense that $\Sigma_{\C[G]}\circ\dcom_{G} = \com_{G}$, and its counit and antipode are given by
\[
\begin{array}{rccc}
\dcou_{G}: & \C[G] & \to & \C \\
& \lambda_{g} & \mapsto & 1
\end{array}
, \qquad
\begin{array}{rccc}
\dS_{G}: & \C[G] & \to & \C[G] \\
& \lambda_{g} & \mapsto & \lambda_{g^{-1}}
\end{array}
,
\]
respectively. The bilinear map defined by $\p: K(G) \times \C[G] \to \C$, $(p , \lambda_{g}) \mapsto p(g)$ allows to identify $\C[G]$ with $\du{K(G)}$ as \as-algebras, indeed for each $\lambda_{g} \in \C[G]$, we have $\lambda_{g}(p) := \p(p,\lambda_{g}) = \sum_{h \in G}p(h)\delta_{g}(h) = (\delta_{g}\cdot\var)(p)$, for all $p \in K(G)$. Here $\delta_{g}$ denotes the Dirac function $h \mapsto \delta_{g,h}$. With the last identification $\p$ becomes the canonical pairing of $\G$ and we have $\dG = (\C[G],\dcom_{G},\dvar_{G})$.
\end{exam}

\subsection*{Actions of algebraic quantum groups}

A right action of an algebraic quantum group $\G$ on a non-degenerate \as-algebra $N$ is a non-degenerate \as-homomorphism $\theta: N \to \M(\pG \odo N)$ satisfying the conditions
\begin{enumerate}[label=\textup{(A\arabic*)}]
\item\label{eq:c1_coaction} $\pG \odo N = \mathrm{span}\{(h \odo 1_{\M(X)})\theta(n) : h \in \pG, m \in N \}$,
\item\label{eq:c2_coaction} $(\com_{\G} \odo \id_{\M(X)})(\theta(m)) = (\id_{\M(A)} \odo \theta)(\theta(m))$ for all $m \in N$,
\item\label{eq:c3_coaction} $(\cou_{\G} \odo \id_{\M(X)})(\theta(m))= m$ for all $m \in N$.
\end{enumerate}
Similarly, left actions of algebraic quantum groups can be defined. We keep in mind the relation of right and left actions: a map $\theta: N \to \M(N \odo \pG)$ is a left action of $\G$ on $N$ if and only if the map $\Sigma\circ\theta: N \to \M(\pG \odo N)$ is a right action of $\Go$ on $N$.

\begin{rema}
The condition~\ref{eq:c3_coaction} in the definition above is equivalent to the ask the injectivity of $\theta$. Moreover, using the language of multiplier Hopf\as-algebras, the conditions~\ref{eq:c1_coaction}, \ref{eq:c2_coaction}, \ref{eq:c3_coaction} says that the \as-homomorphism $\theta$ is a left coaction of the underlying multiplier Hopf \as-algebra $(\pG,\com_{\G})$ on $N$. 
\end{rema}

The following example explain our terminology left and right when we deal with actions of algebraic quantum groups.

\begin{exam}
Let $\cdot: G \times X \to X$ be a left action of a group $G$ on a set $X$. The linear map $\alp : G \to \text{Aut}(K(X))$, $g \mapsto [\alp_{g}(f) : x \mapsto f(g \cdot x)]$, defines a right action of the group $G$ on the \as-algebra of finitely supported complex-valued functions on $X$, $K(X)$. Denote by $\alp: K(X) \odo \C[G] \to K(X)$, the linear extension of $\alp$ to the group \as-algebra $\C[G]$ and consider the multiplier Hopf \as-algebra dual coaction associated with $\alp$, i.e. the \as-homomorphism defined by
\[
\begin{array}{lccccc}
\theta: & K(X) & \to & \M(K(G) \odo K(X)) & \iso & C(G \times X) \\
& f & \mapsto & \displaystyle\sum_{g \in G, x \in X}f(g \cdot x) \delta_{g} \odo \delta_{x} & \mapsto & [(g,x) \mapsto f(g \cdot x)] 
\end{array}
\;,
\]
which is a left coaction of the multiplier Hopf \as-algebra $(K(G),\com_{G})$ on $K(X)$. Observe that $\alp$ can also be obtain as the pullback along the action map $\cdot:G \times X \to X$. Hence, $\theta$ defines a right action of the algebraic quantum group $\G=(K(G),\com_{G},\var_{G})$ on the \as-algebra $K(X)$.
\end{exam}

\begin{rema}
Sometimes when we have to work with a right action $\theta: N \to \M(\pG \odo N)$ of $\G$ on $N$, given $m \in N$ we use the Sweedler type leg notation $\theta(m) =: m_{[-1]} \odo m_{[0]}$ in the following sense: for any $h \in \pG$, we write $hm_{[-1]} \oda m_{[0]} := (h \oda 1)\theta(m)$ and $m_{[-1]}h \oda m_{[0]} := \theta(m)(h \oda 1)$. In case we work with a left action $\theta: N \to \M(N \odo \pG)$, given $m \in N$ we use the Sweedler type leg notation $\theta(m) =: m_{[0]} \odo m_{[1]}$ in a similar way.
\end{rema}

The following technical lemma making use of the Sweedler type leg notation will be use later in computations dealing with actions of algebraic quantum groups.

\begin{lemm}
Let $\theta: N \to \M(\pG \odo N)$ be a right action of $\G$ on a non-degenerate \as-algebra $N$. It holds
\small
\begin{equation}\label{eq:easy_comp2}
m_{[-1]_{(1)}} \odo m_{[-1]_{(2)}} \odo m_{[-1]_{(3)}} \odo m_{[0]_{[-1]}} \odo m_{[0]_{[0]}} = m_{[-1]_{(1)}} \odo m_{[-1]_{(2)}} \odo m_{[0]_{[-1]_{(1)}}} \odo m_{[0]_{[-1]_{(2)}}} \odo m_{[0]_{[0]}}
\end{equation}
\normalsize
for every $m \in N$.
\end{lemm}

\smallskip

\subsection*{Smash product \as-algebra associated with an action}

Let $\G$ be an algebraic quantum group and $\theta$ be a right action of $\G$ on a non-degenerate \as-algebra $N$. The linear map
\[
\begin{array}{lccc}
\lhd_{\theta}: & N \odo \dpG & \to & N \\
& m \odo h & \mapsto & \ome \lhd_{\theta} \ome := (\p(\cdot,\ome) \odo \id)(\theta(m))
\end{array}
\]
defines a right unital non-degenerate module structure of $\dpG$ on $N$ such that
\begin{enumerate}[label=\textup{(DA\arabic*)}]
\item $(mn) \lhd_{\theta} \ome = (m \lhd_{\theta} \ome_{(1)})(n \lhd_{\theta} \ome_{(2)})$
\item $(m \lhd_{\theta} \ome)^{*} = m^{*} \lhd_{\theta} \dS_{\G}(\ome)^{*}$
\end{enumerate}
for all $m,n \in N$ and $\ome \in \dpG$. The last two conditions says that the module map $\lhd_{\theta}$ yields a right action of the multiplier Hopf \as-algebra $(\dpG,\dcom_{\G})$ on $N$. Hereinafter, we will say that $\lhd_{\theta}$ is the right multiplier Hopf \as-algebra action of $(\dpG,\dcom_{\G})$ on $N$ arising from $\theta$. It was shown in \cite{VDZ99}, that any right action of the underlying multiplier Hopf \as-algebra of $\G$ arise in this way.

Similarly, given a left action of an algebraic quantum group $\G$ on $N$, we can define a left action of the multiplier Hopf \as-algebra $(\pG,\com_{\G})$ on $N$.

\begin{rema}\label{rem:opposite_involution}
It is known that we can link left module algebras and right module algebras over multiplier Hopf algebras with integrals using the opposite algebra construction (indeed this is equivalent to use the flip map to make a link between left and right actions of algebraic quantum groups). In the involutive case, we need to pay attention to the compatibility between the module map and the antipode as it is explain in the following lines. The map $\rhd: \pG \odo N \to N$ yields a left module \as-algebra over the multiplier Hopf \as-algebra $(\pG,\com_{\G})$ if and only if the map $\dot\lhd: N \odo \pG^{\op} \to N$, $m \odo h^{\op} \mapsto h \rhd m$ yields a right module map of the opposite multiplier Hopf \as-algebra $(\pG^{\op},\com^{\op}_{\G})$ on $N$. Indeed, it is enough observe that
\begin{align*}
S_{\G}(h)^{*} \rhd m^{*} = m^{*} \dot\lhd (S^{-2}_{\G}(S^{-1}_{G}(h)^{*}))^{\op} = m^{*} \dot\lhd (S^{-1}_{\G}(h)^{\op})^{*} = m^{*} \dot\lhd S_{\Gc}(h^{\op})^{*},
\end{align*}
for each $m \in N$ and $h \in \G$.
\end{rema}

\begin{exam}[The canonical actions by convolution]
Consider the comultiplication of an algebraic quantum group $\com_{\G}$ as a right action of $\G$ on $\pG$. The right multiplier Hopf \as-algebras action of $(\dpG,\dcom)$ on $\pG$ arising from $\com_{\G}$ is given by $h \blhd \ome := (\p(\cdot,\ome) \odo \id)\com_{\G}(h)$ for all $h \in \pG$ and $\ome \in \dpG$. This action will be called {\em the canonical right action by convolution of $(\dpG,\com_{\G})$ on $\G$}. Similarly, if we consider $\com_{\G}$ as a left action of $\G$ on $\pG$, {\em the canonical left action by convolution of $(\dpG,\com_{\G})$ on $\G$} is given by the left module action $\ome \brhd h := (\id \odo \p(\cdot,\ome))\com_{\G}(h)$ for all $h \in \pG$ and $\ome \in \dpG$.
\end{exam}

\begin{exam}[The adjoint actions]
The {\em left and right adjoint action} of $(\pG,\com_{\G})$ on $\pG$ are given by $h \brhd_{\mathrm{ad}} g = h_{(1)}g S_{\G}(h_{(2)})$ and $g \blhd_{\mathrm{ad}} h = S_{\G}(h_{(1)})gh_{(2)}$, respectively.
\end{exam}

Consider the right multiplier Hopf \as-algebra action $\lhd_{\theta}$ of $(\dpG,\dcom_{\G})$ on $N$ arising from $\theta$. The vector space $\dpG \sm_{\theta} N := \dpG \odo N$, endowed with the operations given by
\[
(\ome \sm m)(\ome' \sm n) = \ome\ome'_{(1)} \sm (m \lhd_{\theta} \ome'_{(2)})n, \qquad (\ome \sm m)^{*} = \ome^{*}_{(1)} \sm (m^{*} \lhd_{\theta} \ome^{*}_{(2)})
\]
for all $m,n \in N$ and $\ome,\ome' \in \dpG$, defines a non-degenerate \as-algebra called {\em the smash product \as-algebra associated with $\theta$}. Similarly, the smash product \as-algebras associated with left actions of algebraic quantum groups can be defined.

Given two elements $\ome \in \dpG$ and $m \in N$, consider the multipliers $\iota_{\dpG}(\ome) := \ome \sm 1_{\M(N)}$ and $\iota_{N}(m) := 1_{\M(\pG)} \sm m$ of the smash product \as-algebra $\dpG \sm_{\theta} N$ defined by
\[
[\iota_{\dpG}(\ome)](\ome' \sm_{\theta} n) = \ome\ome' \sm_{\theta} n, \quad (\ome' \sm_{\theta} n)[\iota_{\dpG}(\ome)] = gh_{(1)} \sm_{\theta} (n \lhd_{\theta} h_{(2)})
\]
and
\[
[\iota_{N}(m)](\ome' \sm_{\theta} n) = \ome'_{(1)} \sm_{\theta} (m \lhd_{\theta} \ome'_{(2)})n, \quad (\ome' \sm_{\theta} n)[\iota_{N}(m)] = \ome' \sm_{\theta} nm
\]
for all $\ome' \in \dpG$, $n\in N$, respectively. With those notations, $\iota_{\dpG}: \dpG \to \M(\dpG \sm_{\theta} N)$ and $\iota_{N}: N \to \M(\dpG \sm_{\theta} N)$ yield two non-degenerate faithful \as-homomorphisms, called {\em the canonical faithful embeddings of $\dpG$ and $N$ in the smash product \as-algebra $\dpG \sm_{\theta} N$} respectively, such that
\[
\iota_{N}(m)\iota_{\dpG}(\ome) = \iota_{\dpG}(\ome_{(1)})\iota_{N}(m \lhd_{\theta} \ome_{(2)})
\]
for all $\ome \in \dpG$ and $m \in N$. With the canonical faithful embeddings, the smash product \as-algebra $\dpG \sm_{\theta} N$ is given by
\[
\dpG \sm_{\theta} N = \mathrm{span}\{\iota_{\dpG}(\ome)\iota_{N}(m) : \ome \in \dpG, m \in N\} = \mathrm{span}\{\iota_{N}(m)\iota_{\dpG}(\ome) : \ome \in \dpG, m \in N\}.
\]

\begin{exam}
The Heisenberg \as-algebra of an algebraic quantum group $\G$ is define as the smash product \as-algebra $\He(\G) := \pG \sm_{\blhd} \dpG$ associated with the canonical right action by convolution $\blhd$ of $(\dpG,\dcom_{\G})$ on $\pG$. Equivalently, we can show that $\He(\G) = \pG \sm_{\brhd} \dpG$, where $\brhd$ is the canonical left action by convolution of $(\pG,\com_{\G})$ on $\dpG$. 
\end{exam}

\subsection*{The opposite and the conjugate action}

Given a nice right action of an algebraic quantum group, we can produce two other right actions: the {\em opposite} and the {\em conjugate} action.

\begin{prop}[\cite{Ta22_1}]
Let $\theta: N \to \M(\pG \odo N)$ be a right action of $\G$ on a non-degenerate \as-algebra $N$ and $\gamma$ be an algebra automorphism on $N$, with inverse denoted by $\hat\gamma$, such that $\gamma \circ * \circ \gamma \circ * = \id_{N}$. If $\theta \circ \gamma = (S^{-2}_{\G} \odo \gamma)\circ\theta$ or equivalently $\theta \circ \hat\gamma = (S^{2}_{\G} \odo \hat\gamma)\circ\theta$, then
\begin{enumerate}
\item[\textup{(\textrm{op})}]\label{def:op} the map $\theta^{\ops}: N^{\op}_{\gamma} \to \M(\pG^{\co} \odo N^{\op}_{\gamma})$, $m^{\op} \mapsto (S_{\G} \odo {}^{\op})\theta(m)$ yields a right action of $\Go$ on $ N^{\op}_{\gamma}$ called the {\em opposite action of $\theta$ with respect to $\gamma$}.
\item[\textup{(\textrm{co})}]\label{def:co} the map $\theta^{\con}: N^{\op}_{\gamma} \to \M(\pG^{\op} \odo N^{\op}_{\gamma})$, $m^{\op} \mapsto ({}^{\op} \odo {}^{\op})\theta(m)$ yields a right action of $\Gc$ on $ N^{\op}_{\gamma}$ called the {\em conjugate action of $\theta$ with respect to $\gamma$}.
\end{enumerate}
Moreover, considering the \as-isomorphism $S_{\G}: \pG^{\op} \to \pG^{\co}$, $h^{\op} \mapsto S_{\G}(h)$, we have
\[
(S_{\G} \odo \id)\circ\theta^{\con} = \theta^{\ops}.
\]
\end{prop}

\begin{exam}
Consider the right action $\theta = \com_{\G}$ of $\G$ on $\pG$ and the algebra automorphism $\gamma = S^{-2}_{\G}$ on $\pG$. Follows from the properties of the antipode of $\G$ that $\gamma \circ * \circ \gamma \circ * = \id_{\pG}$ and $\theta\circ\gamma = (S^{-2}_{\G} \odo \gamma)\circ\theta$. In this case, the $\gamma$-opposite \as-algebra of $\pG$ is the underlying \as-algebra of the opposite multiplier Hopf \as-algebra, i.e. $\pG^{\op}_{\gamma} = \pG^{\op}$, and the opposite action is $\theta^{\con} = \com^{\op}_{\G}$.
\end{exam}

From now on, we will mainly work with right actions of algebraic quantum groups unless we mention otherwise.

\subsection{Measured multiplier Hopf \as-algebroids}\label{subsec:MMHAd}

We recall the definition of a measured multiplier Hopf \as-algebroid, our main reference is \cite{T16, T17}. A measured multiplier Hopf \as-algebroid is a collection
\[
\mathcal{A} = (A,B,C,t_{B},t_{C},\com_{B},\com_{C},\mu_{B},\mu_{C},{}_{B}\psi_{B},{}_{C}\phi_{C}),
\]
which consists of
\begin{enumerate}[label=\textup{(M\arabic*)}]
\item an idempotent \as-algebra $A$, called {\em the total algebra of $\mathcal{A}$}, and two \as-algebras $B$, $C$, called {\em the right and the left base algebra of $\mathcal{A}$} respectively, related with $A$ via two non-degenerate faithful embeddings $\iota_{B}: B \to \M(A)$, $\iota_{C}:C \to \M(A)$ whose images commute with each other and making the canonical modules ${}_{B}A$, $A_{B}$, ${}_{C}A$, $A_{C}$ idempotents;

\item two anti-isomorphisms $t_{B}: B \to C$ and $t_{C}: C \to B$ such that
\begin{enumerate}[label=\textup{(\roman*)}]
\item $t_{B}\circ *\circ t_{C}\circ * = \id_{C}$ or equivalently $t_{C}\circ *\circ t_{B}\circ * = \id_{B}$,
\item $\lbt{A}{C}$ is non-degenerate as a right module over $1 \odo A$ and over $A \odo 1$,
\item $\rbt{A}{B}$ is non-degenerate as a left module over $1 \odo A$ and over $A \odo 1$.
\end{enumerate}
Here, $A^{C}$ is the right module arising from the right module ${}_{B}A$ and the anti-isomorphism $t_{C}$, and ${}^{B}A$ is the left module arising from the right module $A_{C}$ and the anti-isomorphism $t_{B}$;

\item two morphisms $\com_{B}: A \to \rtak{A}{B}$ and $\com_{C}: A \to \ltak{A}{C}$, from $A$ to the right and left extended Takeuchi product respectively, called {\em the right and left comultiplication of $\mathcal{A}$} respectively, such that
\begin{enumerate}[label=\textup{(\roman*)}]
\item they are {\em bilinear} with respect to $B$ and $C$ in the sense
\begin{align*}
\com_{B}(\iota_{B}(x)\iota_{C}(y)a\iota_{B}(x')\iota_{C}(y')) & = (\iota_{C}(y) \odo \iota_{B}(x))\com_{B}(a)(\iota_{C}(y') \odo \iota_{B}(x')), \\
\com_{C}(\iota_{B}(x)\iota_{C}(y)a\iota_{B}(x')\iota_{C}(y')) & = (\iota_{C}(y) \odo \iota_{B}(x))\com_{C}(a)(\iota_{C}(y') \odo \iota_{B}(x')),
\end{align*}
for all $a \in A$, $x,x' \in B$ and $y,y' \in C$;

\item they are {\em coassociative} in the sense that
\begin{align*}
(a \odo 1 \odo 1)((\com_{B} \odo \id)((1 \odo c)\com_{B}(b))) & = (1 \odo 1 \odo c)((\id \odo \com_{B})((a \odo 1)\com_{B})b))), \\
(\com_{C} \odo \id)(\com_{C}(b)(1 \odo c))(a \odo 1 \odo 1) & = (\id \odo \com_{C})(\com_{C}(b)(a \odo 1))(1 \odo 1 \odo c), \\
(\com_{C} \odo \id)((1 \odo c)\com_{B}(b))(a \odo 1 \odo 1) & = (1 \odo 1 \odo c)((\id \odo \com_{B})(\com_{C}(b)(a \odo 1))), \\
(a \odo 1 \odo 1)((\com_{B} \odo \id)(\com_{C}(b)(1 \odo c))) & = ((\id \odo \com_{C})((a \odo 1)\com_{B}(b)))(1 \odo 1 \odo c),
\end{align*}
for all $a,b,c \in A$;

\item they are {\em involutive} in the sense
\begin{align*}
\com_{B}(a)^{(* \overline{\times} *)} & = \com_{C}(a^{*}),
\end{align*}
for all $a \in A$, i.e. we have $\com_{B}((b \odo c)a)^{(-)^{*} \odo (-)^{*}} = \com_{C}(a^{*})(b^{*} \odo c^{*})$ for all $a,b,c \in A$.
\end{enumerate}

\smallskip

Moreover, there exist an anti-isomorphism $S: A \to A$, called {\em the antipode of $\mathcal{A}$}, and bimodule maps $\cou_{B} \in \textrm{Hom}({}^{B}A_{B},{}_{B}B_{B})$ and ${}_{C}\cou \in \textrm{Hom}({}_{C}A^{C},{}_{C}C_{C})$, called {\em the right and left counit of $\mathcal{A}$} respectively, satisfying the conditions
\begin{enumerate}[label=\textup{(\roman*)}]
\item for all $x,x' \in B$, $y,y' \in C$, $a \in A$, it holds
\[
S(\iota_{C}(t_{B}(x))\iota_{B}(t_{C}(y))a\iota_{C}(t_{B}(x'))\iota_{B}(t_{C}(y'))) = \iota_{B}(x')\iota_{C}(y')S(a)\iota_{B}(x)\iota_{C}(y).
\]
In this case, we have $S_{B} := S|_{\iota_{B}(B)} = {}_{\iota}t^{-1}_{C}$, $S_{C} := S|_{\iota_{C}(C)} = {}_{\iota}t^{-1}_{B}$, where $S$ is regarded as an anti-isomorphism on $\M(A)$;
\item $S\circ * \circ S \circ * = \id$;
\item for all $a,b \in A$, it hold
\[
(\cou_{B} \odot \id)((1 \odo b)\com_{B}(a)) = ba = (\id \odot \cou_{B})((b \odo 1)\com_{B}(a)),
\]
\[
({}_{C}\cou \odot \id)(\com_{C}(a)(1 \odo b)) = ab = (\id \odot {}_{C}\cou)(\com_{C}(a)(b \odo 1)).
\]
Here, we are using the slice maps
\[
(\cou_{B} \odot \id)(a \odo b) := b\iota_{C}(t_{B}(\cou_{B}(a))), \quad (\id \odot \cou_{B})(a \odo b) := a\iota_{B}(\cou_{B}(b)),
\]
\[
({}_{C}\cou \odot \id)(a' \odo b') := \iota_{C}({}_{C}\cou(a'))b', \quad (\id \odot {}_{C}\cou)(a' \odo b') := \iota_{B}(t_{C}({}_{C}\cou(b')))a'
\]
for all elementary tensor $a \otimes b \in \rbt{A}{B}$ and $a' \odo b' \in \lbt{A}{C}$;

\item for all $a,b \in A$, it hold
\[
\m(S \odo \id)(\com_{C}(a)(1 \odo b)) = \iota_{B}(\cou_{B}(a))b,
\]
and
\[
\m(\id \odo S)((a \odo 1)\com_{B}(b)) = a\iota_{C}({}_{C}\cou(b)).
\]
\end{enumerate}

\item two bimodule maps ${}_{B}\psi_{B} \in \textrm{Hom}({}_{B}A_{B},{}_{B}B_{B})$ and ${}_{C}\phi_{C} \in \textrm{Hom}({}_{C}A_{C},{}_{C}C_{C})$, called {\em the right and the left partial integral of $\mathcal{A}$} respectively, and two non-zero positive faithful functionals $\mu_{B}$ on $B$ and $\mu_{C}$ on $C$, called {\em the right and the left base weight of $\mathcal{A}$} respectively, satisfying the conditions
\begin{enumerate}[label=\textup{(\roman*)}]
\item for all $a, b \in A$, it hold
\[
({}_{B}\psi_{B} \odot \id)((1 \odo b)\com_{B}(a)) = b\iota_{B}({}_{B}\psi_{B}(a)),
\]
and
\[
(\id \odot {}_{C}\phi_{C})(\com_{C}(b)(a \odo 1)) = \iota_{C}({}_{C}\phi_{C}(b))a.
\]
Here, we are using the slice maps
\[
({}_{B}\psi_{B} \odot \id): \rbt{A}{B} \to A, \; a \odo b \mapsto b\iota_{C}(t_{B}({}_{B}\psi_{B}(a))) = bS^{-1}_{C}(\iota_{B}({}_{B}\psi_{B}(a))),
\]
\[
(\id \odot {}_{C}\phi_{C}): \lbt{A}{C} \to A, \; a \odo b \mapsto \iota_{B}(t_{C}({}_{C}\phi_{C}(b)))a = S^{-1}_{B}(\iota_{C}({}_{C}\phi_{C}(b)))a;
\]

\item $\mu_{B}\circ t_{C} = \mu_{C}$, $\mu_{C}\circ t_{B} = \mu_{B}$, and $\mu_{B}\circ \cou_{B} = \mu_{C} \circ {}_{C}\cou$,

\item $\psi := \mu_{B}\circ {}_{B}\psi_{B}$ yields a positive faithful functional called {\em the right total integral of $\mathcal{A}$} and $\phi := \mu_{C} \circ {}_{C}\phi_{C}$ yields a positive faithful functional called {\em the left total integral of $\mathcal{A}$}.
\end{enumerate}
\end{enumerate}

\begin{defi}\label{def:MHAd_Kac}
The measured multiplier Hopf \as-algebroid $\mathcal{A}$ is said to be {\em of Kac type}, if its antipode is *-preserving. We say that $\mathcal{A}$ is {\em unital}, if its total algebra and its base algebras are unital \as-algebras. We say that $\mathcal{A}$ is {\em unimodular}, if its left total integral $\phi$ is equal to its right total integral $\psi$.
\end{defi}

Given a measured multiplier Hopf \as-algebroid $\mathcal{A}$, there are its modular automorphisms and its modular elements. By \cite[Theorem~2.2.3]{T17}, the base weights $\mu_{B}$, $\mu_{C}$ and the total integrals $\phi$, $\psi$ admit modular automorphisms denoted by $\sigma_{B}$, $\sigma_{C}$, $\sigma^{\phi}$ and $\sigma^{\psi}$, respectively. Moreover, they satisfy the condition $\sigma_{C} = \sigma^{\phi}|_{C} = S_{B}\circ S_{C}$ and $\sigma_{B} = \sigma^{\psi}|_{B} = S^{-1}_{B}\circ S^{-1}_{C}$. Recall that if $a \in A$, we use the notation $a \cdot \phi$, $\phi \cdot a$ for the linear maps $a' \mapsto \phi(a'a)$ and $a' \mapsto \phi(aa')$ from $A$ to $\C$, respectively. With those notations, we have $\phi\cdot a = \sigma^{\phi}(a)\cdot\phi$ and $\psi\cdot a = \sigma^{\psi}(a)\cdot\psi$ for all $a \in A$. On the other hand, by \cite[Theorem~6.4]{T16}, there exist invertible multiplier $\delta, \delta^{+}, \delta^{-} \in \M(A)$, called the {\em modular elements of $\mathcal{A}$} such that $S(\delta^{+})=(\delta^{-})^{-1}$, $(\delta^{-})^{*} = \delta^{+}$, $\psi = \delta\cdot\phi$, $\phi\circ S = \delta^{+}\cdot\phi$ and $\phi\cdot S^{-1} = \phi\cdot\delta^{-}$

\smallskip

\begin{defi}
Let $\mathcal{A'}=(A',B',C',t_{B'},t_{C'},\com_{B'},\com_{C'},\mu_{B'},\mu_{C'},{}_{B'}\psi_{B'},{}_{C'}\phi_{C'})$ be another measured multiplier Hopf \as-algebroid. Denote by ${}_{\iota}t_{B'}$ and ${}_{\iota}t_{C'}$ the lifting anti-isomorphisms arising from $t_{B'}$ and $t_{C'}$ respectively, through the faithful embeddings $\iota_{B'}$ and $\iota_{C'}$. A {\em morphism from $\mathcal{A}$ to $\mathcal{A}'$} is a non-degenerate \as-homomorphism $\pi: A \to \M(A')$ satisfying the conditions
\begin{enumerate}[label=(\arabic*)]
\item $\pi(\iota_{B}(B))\iota_{B'}(B') = \iota_{B'}(B')$ and $\pi(\iota_{C}(C))\iota_{C'}(C') = \iota_{C'}(C')$;
\item $\pi\circ\iota_{C}\circ t_{B} = {}_{\iota}t_{B'}\circ\pi\circ\iota_{B}$ and $\pi\circ\iota_{B}\circ t_{C} = {}_{\iota}t_{C'}\circ\pi\circ\iota_{C}$;
\item $\com_{B'}\circ\pi = (\rtak{\pi}{B})\circ\com_{B}$ and $\com_{C'}\circ\pi = (\ltak{\pi}{C})\circ\com_{C}$.
\end{enumerate}
\end{defi}

\subsection*{The Pontrjagin dual of a measured multiplier Hopf \as-algebroid}

Let $\mathcal{A}$ be the measured multiplier Hopf \as-algebroid
\[
(A,B,C,t_{B},t_{C},\com_{B},\com_{C},\mu_{B},\mu_{C},{}_{B}\psi_{B},{}_{C}\phi_{C}).
\]
Denote by $\phi=\mu_{C}\circ{}_{C}\psi_{C}$ its left total integral. Following \cite{T17}, the Pontrjagin dual of $\mathcal{A}$, denoted by $\du{\mathcal{A}}$, is given by the measured multiplier Hopf \as-algebroid
\[
(\hat{A},C,B,t^{-1}_{B},t^{-1}_{C},\dcom_{C},\dcom_{B},\mu_{C},\mu_{B},{}_{B}\du{\psi}_{B},{}_{C}\du{\phi}_{C})
\]
where $\du{A} = A \cdot \phi = \phi \cdot A$ and the faithful embeddings $\hat{\iota}_{B}: B \to \M(\du{A})$, $\hat{\iota}_{C}: C \to \M(\du{A})$ are defined by
\[
[\hat{\iota}_{B}(x)]\upsilon = \upsilon\cdot t_{B}(x), \qquad \upsilon[\hat{\iota}_{B}(x)] = \upsilon\cdot x, \;
\]
and
\[
[\hat{\iota}_{C}(y)]\upsilon = y\cdot\upsilon, \qquad \upsilon[\hat{\iota}_{C}(y)] = t_{C}(y)\cdot\upsilon
\]
for all $x \in B$, $y \in C$ and $\upsilon \in \du{A}$. By construction, we have $\hat{\iota}_{B}(x)\hat{\iota}_{C}(y) = \hat{\iota}_{C}(y)\hat{\iota}_{B}(x)$ inside $\M(\du{A})$, for all $x \in B$, $y \in C$.

Similar to the case of algebraic quantum groups, we prefer the use of the notion of pairing in order to give the structure of the Pontrjagin dual $\du{\mathcal{A}}$. By \cite[Theorem 5.10]{TVDW22}, the bilinear map
\[
\begin{array}{lccc}
\Pe: & A \times \du{A} & \to & \C \\
& (a,\upsilon) & \mapsto & \upsilon(a)
\end{array}
\]
yields a pairing of measured multiplier Hopf \as-algebroids called the {\em canonical pairing between $\mathcal{A}$ and its Pontrjagin dual $\du{\mathcal{A}}$}. The comultiplication maps $\dcom_{C}: \du{A} \to \du{A} \rtak{}{C} \du{A}$ and $\dcom_{B}: \du{A} \to \du{A} \ltak{}{B} \du{A}$ are given by the conditions
\[
\Pe^{2}(a \odo a',\dcom_{B}(\upsilon)(1 \odo \upsilon')) = \Pe^{2}((a \odo 1)\com_{B}(a'),\upsilon \odo \upsilon')
\]
and
\[
\Pe^{2}(a \odo a',(\upsilon \odo 1)\dcom_{C}(\upsilon')) = \Pe^{2}(\com_{C}(a)(1 \odo a'),\upsilon \odo \upsilon')
\]
for all $a,a' \in A$ and $\upsilon, \upsilon' \in \du{A}$, respectively. The involution and the antipode $\dS$ on $\du{A}$ are given by
\[
\Pe(a,\upsilon^{*}) = \overline{\Pe(S(a)^{*},\upsilon)} \quad \text{ and } \quad \Pe(S(a),\upsilon) = \Pe(a,\dS(\upsilon))
\]
for all $a \in A$ and $\upsilon \in \du{A}$, respectively. Observe that we also have $\Pe(a,\dS(\upsilon)^{*}) = \overline{\Pe(a^{*},\upsilon)}$ for all $a \in A$ and $\upsilon \in \du{A}$.

\subsection*{The canonical actions by convolution}

Similar to the case of algebraic quantum groups, there are canonical module maps associated with measured multiplier Hopf \as-algebroids. Following \cite[Proposition 23]{T16} and \cite[Lemma~3.2.2]{T17}, there are module maps $\brhd : \du{A} \odo A \to A$ and $\blhd: A \odo \du{A} \to A$ defined by
\[
(a \cdot \phi) \brhd b = (\id \odot {}_{C}\phi_{C})(\com_{C}(b)(1 \odo a)), \quad b \blhd (\phi\cdot a) = ({}_{B}\psi_{B} \odot \id)((a \odo 1)\com_{B}(b))
\]
for all $a,b \in A$, respectively. Those module maps also satisfy
\[
\upsilon \brhd S(a) = S(a \blhd \dS(\upsilon)), \qquad S(a) \blhd \upsilon = S(\dS(\upsilon) \brhd a),
\]
\[
(\upsilon \brhd a)^{*} = \du{S}(\upsilon)^{*} \brhd a^{*}, \qquad (a \blhd \upsilon)^{*} = a^{*} \blhd \du{S}(\upsilon)^{*}
\]
for all $\upsilon\in \du{A}$ and $a \in A$. Moreover, given two element $\upsilon, \upsilon' \in \du{A}$, if we have $\upsilon \brhd a = \upsilon' \brhd a$ for all $a \in A$, thus necessarily $\upsilon = \upsilon'$. The last property is known as the {\em faithfulness of the left action $\brhd$}, and it will be useful in case we need to show that two elements belonging to the underlying \as-algebra of the Pontrjagin dual of a measured multiplier Hopf \as-algebroid are equal. A similar property holds for the right action $\blhd$. 

\section{Yetter--Drinfeld \as-algebras over algebraic quantum groups}\label{sec:bcyd}

In this section, we give an important result on unital braided commutative Yetter--Drinfeld \as-algebras over algebraic quantum groups of compact type. First, we recall the basics on Yetter--Drinfeld structures over algebraic quantum groups, we refer the interested reader to \cite{Ta22_1} for more details and proofs.

Let $\G = (\pG,\com,\varphi)$ be an algebraic quantum group with algebraic multiplicative unitary $U$. A {\em (right) Yetter--Drinfeld $\G$-\as-algebra} is a triplet $(N,\theta,\dtheta)$, where $N$ is a \as-algebra, $\theta : N \to \M(\pG \odo N)$ is a right action of $\G$ on $N$ and $\dtheta: N \to \M(\dpG \odo N)$ is a right action of $\dG$ on $N$, satisfying the {\em Yetter--Drinfeld condition}
\begin{equation}\label{eq:r-YD-quantum}\tag{YD}
(\id_{\dpG} \odo \theta)\circ\dtheta = (\Sigma \odo \id_{N})\circ(\ad(U) \odo \id_{N})\circ(\id_{\pG} \odo \dtheta)\circ\theta.
\end{equation}
In other words, if the diagram
\[
\xymatrix@C=6pc@R=2pc{N\ar[r]^-{\dtheta}\ar[d]_-{\theta} & \M(\dpG \odo N)\ar[r]^-{\id_{\dpG} \odo \theta} & \M(\dpG \odo \pG \odo N) \\ \M(\pG \odo N)\ar[r]_-{\id_{\pG} \odo \dtheta} & \M(\pG \odo \dpG \odo N)\ar[r]_-{\ad(U) \odo \id_{N}} & \M(\pG \odo \dpG \odo N)\ar[u]_-{\Sigma \odo \id_{N}}}
\]
is commutative. A (right) Yetter--Drinfeld $\G$-\as-algebra $(N,\theta,\dtheta)$ is called {\em braided commutative}, if for each $m,n \in N$, we have
\begin{equation*}\label{eq:r-bc-quantum}\tag{BC}
\theta^{\con}(m^{\op})\du{\theta^{\ops}}(n^{\op}) = \dtheta^{\ops}(n^{\op})\theta^{\con}(m^{\op})
\end{equation*}
inside $\M(\He(\G) \odo N^{\op})$. A morphism of (right) Yetter--Drinfeld $\G$-\as-algebras $f: (N,\theta,\dtheta) \to (N',\theta',\dtheta')$ is a non-degenerate \as-homomorphism $f: N \to \M(N')$ such that $(\id \odo f)\circ\theta = \theta'\circ f$ and $(\id \odo f)\circ\dtheta = \dtheta'\circ f$. The category of (right) Yetter--Drinfeld $\G$-\as-algebras will be denoted by $\aYD_{\G}$ and the subcategory of braided commutative (right) Yetter--Drinfeld $\G$-\as-algebras will be denoted by $\aYD^{\textrm{bc}}_{\G}$.

\begin{defi}
A tuple $(N,\theta,\dtheta,\mu)$ is called a {\em measured (right) Yetter--Drinfeld $\G$-\as-algebra}, if $(N,\theta,\dtheta)$ is a (right) Yetter--Drinfeld $\G$-\as-algebra and $\mu$ is a non-zero positive faithful functional on $N$, called a {\em Yetter--Drinfeld integral for $(N,\theta,\dtheta)$}, which satisfies the $\theta$-invariant and $\dtheta$-invariant condition i.e. the conditions
\[
(\id_{\pG} \odo \mu)((h \odo 1_{N})\theta(m)) = \mu(m)h \quad \text{ and } \quad (\id_{\dG} \odo \mu)((\ome \odo 1_{N})\dtheta(m)) = \mu(m)\ome
\]
for all $h \in \pG$, $\ome \in \dpG$ and $m\in N$, respectively. If in addition $\mu$ satisfies a weak KMS property, then there is a unique algebra automorphism $\sigma^{\mu}: N \to N$, called the {\em modular automorphism of $\mu$}, verifying $\mu(mn) = \mu(n\sigma^{\mu}(m))$ for all $m,n \in N$. In the case $\sigma^{\mu} = \id_{N}$, the Yetter--Drinfeld integral $\mu$ is called {\em tracial}.
\end{defi}

\begin{exam}\label{ex:yd_transformation_groupoid}
Let $\cdot : G \times X \to X$ be a right action of a group $G$ on a set $X$. Consider the left action of the algebraic quantum group $\Gc = (K(G)^{\op},\com^{\op},\var^{\op})$ on $K(X)$ defined by
\[
\begin{array}{rccc}
\theta: & K(X) & \to & \M(K(G)^{\op} \oti K(X)) \\
& f & \mapsto & \displaystyle\sum_{g \in G, x \in X} f(g\cdot x)\delta_{g} \oti \delta_{x}
\end{array}
\]
with dual right action of the Hopf \as-algebra $(\C[G],\com'^{\co})$ on the \as-algebra $K(X)$ given by
\[
\begin{array}{rccc}
\lhd_{\theta}: & K(X) \odo \C[G]  & \to & K(X) \\
& f \odo \lambda_{g} & \mapsto & [x \mapsto f(g \cdot x)]  
\end{array},
\]
and the trivial right action of the algebraic quantum group $\dG^{\ops}=(\C[G],\com',\var')$ on $K(X)$
\[
\begin{array}{rccc}
\dtheta: & K(X) & \to & \M(\C[G] \otimes K(X)) \\
& f & \mapsto & \lambda_{e} \odo f
\end{array}
\]
with dual left action of the multiplier Hopf \as-algebra $(K(G),\com)$ on the \as-algebra $K(X)$ given by
\[
\begin{array}{rccc}
\lhd_{\dtheta}: & K(X) \oti K(G) & \to & K(X) \\
& f \oti p & \mapsto & p(e)f
\end{array}.
\]
The triplet $(K(X),\theta,\dtheta)$ yields a braided commutative (right) Yetter--Drinfeld $\Gc$-\as-algebra.

\begin{prop}\label{prop:YD_integral_commutative}
Let $\nu: X \to \R^{+}_{0}$ be a non-zero function such that $\nu(g\cdot x) = \nu(x)$ for all $x \in X$ and $g \in G$. Then, the non-zero positive faithful functional $\mu_{\nu}: K(X) \to \C$, given by $\mu_{\nu}(f) = \sum_{x \in X}\nu(x)f(x)$, yields a Yetter--Drinfeld integral for $(K(X),\theta,\dtheta)$. Moreover, any Yetter--Drinfeld integral for $(K(X),\theta,\dtheta)$ arises in this way.
\end{prop}
\pr
First, observe that given $x \in X$, we have $\theta(\delta_{x}) = \sum_{g \in G} \delta_{g} \odo \delta_{g^{-1}\cdot x}$ for all Dirac function $\delta_{x} \in K(X)$ with $x \in X$. Then, we have
\[
(\id \odo \mu_{\nu})(\theta(\delta_{x})) = \sum_{g \in G, y \in X}\delta_{g} \odo \nu(y)\delta_{g^{-1}\cdot x}(y) = \sum_{g \in G}\nu(g\cdot x)\delta_{g} = \nu(x)1_{G} = \mu_{\nu}(\delta_{x})1_{G}
\]
and $(\id \odo \mu_{\nu})(\dtheta(\delta_{x})) = \mu_{\nu}(\delta_{x})\lambda_{e}$ for all $x \in X$. Since $f = \sum_{x \in X} f(x)\delta_{x}$ for any $f \in K(X)$, thus $\mu_{\nu}$ is a Yetter--Drinfeld integral for $(K(X),\theta,\dtheta)$.

Now, let $\mu$ be a Yetter--Drinfeld integral for $(K(X),\theta,\dtheta)$. Since $\mu \in K(X)^{\vee} \iso \C[X]$, then there is a map $\nu: X \to \C$ such that $\mu = \sum_{x \in X}\nu(x)\mathrm{ev}_{x}$. Because $\mu$ is a non-zero positive faithful functional and $\mu(\delta_{x}) = \mu(\delta_{x}\delta^{*}_{x}) = \nu(x)$ for all $x \in X$, then $\nu$ is a non-zero function and $\mathrm{Im}(\nu) \subseteq \R^{+}_{0}$. On the other hand, because $(\id \odo \mu)(\theta(f)) = \theta(f)1_{G}$ for all $f \in K(S)$, thus for a fixed $x \in X$, we have
\[
\sum_{g \in G}\mu(\delta_{g^{-1}\cdot x})\delta_{g} = (\id \odo \mu)(\theta(\delta_{x})) = \mu(\delta_{x})1_{G} = \sum_{g \in G} \mu(\delta_{x})\delta_{g}.
\]
Since the vectors $\{\delta_{g}\}_{g \in G}$ are linear independent in $K(G)$, we have $\nu(x) = \mu(\delta_{x}) = \mu(\delta_{g^{-1}\cdot x}) = \nu(g^{-1}\cdot x)$ for all $g \in G$.
\fin
\end{exam}

\medskip

In this work, we will mainly work with braided commutative (right) Yetter--Drinfeld $\Gc$-\as-algebras, where $\G$ is an algebraic quantum group of compact type. The next proposition will be useful to obtain examples of braided commutative (right) Yetter--Drinfeld $\Gc$-\as-algebras.

\begin{prop}\label{prop:bc_yd_equivalence}
Let $\G=(\pG,\com_{\G},\varphi_{\G})$ be an algebraic quantum group and $N$ be an \as-algebra. Consider a right action $\lhd: N \odo \pG \to N$ and a right coaction $\delta: N \to \M(N \odo \pG)$. Denote by $\hat\theta_{\lhd}: N \to \M(\dpG \odo N)$ the left coaction of the multiplier Hopf \as-algebra $(\dpG,\dcom_{\G})$ on $N$ associated with the action $\lhd$ and by $\theta_{\delta}:=(S^{-1}_{\G} \odo \id)\circ\Sigma\circ\delta:N \to \M(\pG^{\op} \odo N)$ the left coaction of $(\pG^{\op},\com^{\op})$ on $N$. If $\p: \pG \times \dpG \to \C$ denotes the canonical pairing associated with the algebraic quantum group $\G$, i.e. $\p(h,\varphi_{h'}) = \varphi(h'h)$ for all $h,h' \in \pG$, then the following statements are equivalent
\begin{enumerate}[label=\textup{(E\arabic*)}]
\item $(N,\lhd,\delta)$ is a braided commutative (right-right) Yetter--Drinfeld $(\pG,\com_{\G})$-\as-algebra;
\item $(N,\lhd,\theta_{\delta})$ is a braided commutative (right-left) Yetter--Drinfeld $(\pG,\com_{\G})$-\as-algebra;
\item $(N,\delta,\hat\theta_{\lhd})$ is a braided commutative (right-left) Yetter--Drinfeld \as-algebra over $\p$;
\item $(N,\theta_{\delta},\hat\theta_{\lhd})$ is a braided commutative (right) Yetter--Drinfeld $\Gc$-\as-algebra.
\end{enumerate}
\end{prop}
\pr
Follows from \cite[Proposition 3.5, B.4 and B.5]{Ta22_1}.
\fin

\subsection{The canonical automorphisms associated with a unital Yetter--Drinfeld \as-algebra over an algebraic quantum group of compact type}

From now on, $\G=(\pG,\com_{\G},\varphi_{\G})$ will denoted an algebraic quantum group of compact type, i.e. $(\pG,\com_{\G})$ is a unital Hopf \as-algebra endowed with a left and right integral $\varphi_{\G}$, and $(N,\theta,\dtheta)$ will denoted a unital braided commutative (right) Yetter--Drinfeld $\Gc$-\as-algebra. Consider the right action $\lhd_{\dtheta}: N \odo \pG \to N$, defined by $m \lhd_{\dtheta} h = (\p(h,\cdot) \odo \id)\dtheta(m)$, of $(\pG,\com_{\G})$ on $N$ induced by $\dtheta$ and the left coaction $\theta: N \to \pG^{\op} \odo N$ of $(\pG^{\op},\com^{\op}_{\G})$ on $N$. Keep in mind also the Sweedler type leg notation $m^{\op}_{[-1]} \odo m_{[0]}:= \theta(m)$ for all $m \in N$. By Proposition~\ref{prop:bc_yd_equivalence}, the triplet $(N,\lhd_{\dtheta},\theta)$ yields a unital braided commutative (right-left) Yetter--Drinfeld $(\pG,\com_{\G})$-\as-algebra, i.e. we have
\begin{equation}\label{eq:yd_condition}
\theta(m \lhd_{\dtheta} h) =  (S^{-1}_{\G}(h_{(3)})m_{[-1]}h_{(1)})^{\op} \;\odo\; (m_{[0]} \lhd_{\dtheta} h_{(2)})
\end{equation}
and
\begin{equation}\label{eq:bc_condition}
mn = (n \lhd_{\dtheta} m_{[-1]})m_{[0]} = n_{[0]}(m \lhd_{\dtheta} S_{\G}(n_{[-1]}))
\end{equation}
for any $h \in \pG$, $m,n \in N$.

\begin{theo}\label{theo:gamma_hatgamma}
Let $(N,\theta,\dtheta)$ be a unital braided commutative Yetter--Drinfeld $\Gc$-\as-algebra. Consider the linear maps
\[
\begin{array}{lccc}
\gamma_{\theta}: & N & \to & N \\
& m & \mapsto & m_{[0]} \lhd_{\dtheta} S^{-1}_{\G}(m_{[-1]})
\end{array},
\quad
\begin{array}{lccc}
\hat\gamma_{\theta}: & N & \to & N \\
& m & \mapsto & m_{[0]} \lhd_{\dtheta} S^{2}_{\G}(m_{[-1]})
\end{array}.
\]
It holds
\begin{enumerate}[label=\textup{(CA\arabic*)}]
\item\label{itpr:igualdad_general} $* \circ \gamma_{\theta} = \hat\gamma_{\theta} \circ *$, $\theta\circ\gamma_{\theta} = (S^{2}_{\Gc} \odo \gamma_{\theta})\circ\theta$ and $\theta\circ\hat\gamma_{\theta} = (S^{-2}_{\Gc} \odo \hat\gamma_{\theta})\circ\theta$.
\end{enumerate}
Moreover, we have
\begin{enumerate}[label=\textup{(CA\arabic*)},resume]
\item\label{itpr:action_on_lambda} $\gamma_{\theta}(m \lhd_{\dtheta} h) = \gamma_{\theta}(m) \lhd_{\dtheta} S^{-2}_{\G}(h)$ and $\hat\gamma_{\theta}(m \lhd_{\dtheta} h) = \hat\gamma_{\theta}(m) \lhd_{\dtheta} S^{2}_{\G}(h)$, for all $m \in N$, $h \in \pG$;
\item\label{itpr:lamd_ofside} $m_{[0]}n \lhd_{\dtheta} S^{-1}_{\G}(m_{[-1]}) = n\gamma_{\theta}(m)$ and $nm_{[0]} \lhd_{\dtheta} S^{2}_{\G}(m_{[-1]}) = \hat\gamma_{\theta}(m)n$ for all $m, n \in N$;
\item $\gamma_{\theta}(\hat{\gamma_{\theta}}(m)) = m$ and $\hat\gamma_{\theta}(\gamma_{\theta}(m)) = m,$ for all $m \in N$;
\item $\gamma_{\theta}(mn) = \gamma_{\theta}(m)\gamma_{\theta}(n)$, $\hat\gamma_{\theta}(mn) = \hat\gamma_{\theta}(m)\hat\gamma_{\theta}(n)$ for all $m, n \in N$;
\item $\gamma_{\theta}(\gamma_{\theta}(m^{*})^{*}) = m$ and $\hat\gamma_{\theta}(\hat\gamma_{\theta}(m^{*})^{*}) = m$ for all $m \in N$;
\item $\dtheta\circ\gamma_{\theta} = (S^{2}_{\dG^{\ops}} \odo \gamma_{\theta})\circ\dtheta$ and $\dtheta\circ\hat\gamma_{\theta} = (S^{-2}_{\dG^{\ops}} \odo \hat\gamma_{\theta})\circ\dtheta$;
\item\label{idpr:action_theta} $\gamma_{\theta}(m \lhd_{\theta} \ome) = \gamma_{\theta}(m) \lhd_{\theta} \dS^{2}_{\G}(\ome)$ and $\hat\gamma_{\theta}(m \lhd_{\theta} \ome) = \hat\gamma_{\theta}(m) \lhd_{\theta} \dS^{-2}_{\G}(\ome)$, for all $m \in N$ and $\ome \in \dpG$;
\item if $\mu$ is a $\dtheta$-invariant non-zero functional on $N$, it holds $\mu = \mu\gamma_{\theta} = \mu\hat\gamma_{\theta}$ and $\mu(mn)=\mu(n\gamma_{\theta}(m)) = \mu(\hat\gamma_{\theta}(n)m)$ for all $m, n \in N$, i.e. $\mu$ satisfies a weak KMS property with modular automorphism $\gamma_{\theta}$;
\end{enumerate}
\end{theo}
\pr
\begin{enumerate}[label=\textup{(A\arabic*)}]
\item\label{it:igualdad_general} Fix $m \in N$. Because,
\[
(m^{*})^{\op}_{[-1]} \odo (m^{*})_{[0]} = \theta(m^{*})=\theta(m)^{*} = (m^{\op}_{[-1]})^{*} \odo (m_{[0]})^{*} = S^{-2}_{\G}((m_{[-1]})^{*})^{\op} \odo (m_{[0]})^{*},
\]
thus
\begin{align*}
\gamma_{\theta}(m)^{*} & = (m_{[0]} \lhd_{\dtheta} S^{-1}_{\G}(m_{[-1]}))^{*} = (m_{[0]})^{*} \lhd_{\dtheta} (m_{[-1]})^{*} \\
& = (m^{*})_{[0]} \lhd_{\dtheta} S^{2}_{\G}((m^{*})_{[-1]}) = \hat\gamma_{\theta}(m^{*}).
\end{align*}
On the other hand, we have
\small
\begin{align*}
\theta(\gamma_{\theta}(m)) & = \theta(m_{[0]} \lhd_{\dtheta} S^{-1}_{\G}(m_{[-1]})) \\
& = (S^{-2}_{\G}(m_{[-1]_{(1)}})m_{[0]_{[-1]}}S^{-1}_{\G}(m_{[-1]_{(3)}}))^{\op} \odo (m_{[0]_{[0]}} \lhd_{\dtheta} S^{-1}_{\G}(m_{[-1]_{(2)}})) & \text {by } \eqref{eq:yd_condition} \\
& = (S^{-2}_{\G}(m_{[-1]_{(1)}})m_{[0]_{[-1]_{(2)}}}S^{-1}_{\G}(m_{[0]_{[-1]_{(1)}}}))^{\op} \odo (m_{[0]_{[0]}} \lhd_{\dtheta} S^{-1}_{\G}(m_{[-1]_{(2)}})) & \text{by } \eqref{eq:easy_comp2}\\
& = S^{-2}_{\G}(m_{[-1]_{(1)}})^{\op} \odo (\cou_{\G}(m_{[0]_{[-1]}})m_{[0]_{[0]}} \lhd_{\dtheta} S^{-1}_{\G}(m_{[-1]_{(2)}})) \\
& = S^{2}_{\Gc}(m^{\op}_{[-1]}) \odo \gamma_{\theta}(m_{[0]})
\end{align*}
\normalsize
and
\small
\begin{align*}
\theta(\hat\gamma_{\theta}(m)) & = \theta(m_{[0]} \lhd_{\dtheta} S^{2}_{\G}(m_{[-1]})) \\
& = (S_{\G}(m_{[-1]_{(3)}})m_{[0]_{[-1]}}S^{2}_{\G}(m_{[-1]_{(1)}}))^{\op} \odo (m_{[0]_{[0]}} \lhd_{\dtheta} S^{2}_{\G}(m_{[-1]_{(2)}})) & \text {by } \eqref{eq:yd_condition} \\
& = (S_{\G}(m_{[0]_{[-1]_{(1)}}})m_{[0]_{[-1]_{(2)}}}S^{2}_{\G}(m_{[-1]_{(1)}}))^{\op} \odo (m_{[0]_{[0]}} \lhd_{\dtheta} S^{2}_{\G}(m_{[-1]_{(2)}})) & \text{by } \eqref{eq:easy_comp2}\\
& = S^{2}_{\G}(m_{[-1]_{(1)}})^{\op} \odo (m_{[0]} \lhd_{\dtheta} S^{2}_{\G}(m_{[-1]_{(2)}})^{\op}) \\
& = S^{-2}_{\Gc}(m^{\op}_{[-1]}) \odo \hat\gamma_{\theta}(m_{[0]}).
\end{align*}
\normalsize
\end{enumerate}

Moreover,

\begin{enumerate}[label=\textup{(CA\arabic*)},resume]
\item\label{it:action_on_lambda} Fix $m \in N$ and $h \in \pG$. By Yetter--Drinfeld condition, we have
\[
\theta(m \lhd_{\dtheta} h) = (S^{-1}_{\G}(h_{(3)})m_{[-1]}h_{(1)})^{\op} \odo (m_{[0]} \lhd_{\dtheta} h_{(2)}),
\]
thus
\begin{align*}
\gamma_{\theta}(m \lhd_{\dtheta} h) & = (m_{[0]} \lhd_{\dtheta} h_{(2)}) \lhd_{\dtheta} S^{-1}_{\G}(h_{(1)})S^{-1}_{\G}(m_{[-1]})S^{-2}_{\G}(h_{(3)}) \\
& = m_{[0]} \lhd_{\dtheta} S^{-1}_{\G}(m_{[-1]})S^{-2}_{\G}(h) \\
& = \gamma_{\theta}(m) \lhd_{\dtheta} S^{-2}_{\G}(h)
\intertext{and}
\hat\gamma_{\theta}(m \lhd_{\dtheta} h) & = (m_{[0]} \lhd_{\dtheta} h_{(2)}) \lhd_{\dtheta} S_{\G}(h_{(3)})S^{2}_{\G}(m_{[-1]})S^{2}_{\G}(h_{(1)}) \\
& = m_{[0]} \lhd_{\dtheta} S^{2}_{\G}(m_{[-1]})S^{2}_{\G}(h) \\
& = \hat\gamma_{\theta}(m) \lhd_{\dtheta} S^{2}_{\G}(h).
\end{align*}

\item\label{it:lamd_ofside} Fix $m,n \in N$. We have
\begin{align*}
m_{[0]}n \lhd_{\dtheta} S^{-1}_{\G}(m_{[-1]}) & = (n \lhd_{\dtheta} m_{[0]_{[-1]}})m_{[0]_{[0]}} \lhd_{\dtheta} S^{-1}_{\G}(m_{[-1]}) & \text{by } \eqref{eq:bc_condition} \\
& = (n \lhd_{\dtheta} m_{[0]_{[-1]}}S^{-1}_{\G}(m_{[-1]_{(2)}}))(m_{[0]_{[0]}} \lhd_{\dtheta} S^{-1}_{\G}(m_{[-1]_{(1)}})) \\
& = n(m_{[0]} \lhd_{\dtheta} S^{-1}_{\G}(m_{[-1]})) = n\gamma_{\theta}(m),
\intertext{and}
nm_{[0]} \lhd_{\dtheta} S^{2}_{\G}(m_{[-1]}) & = m_{[0]_{[0]}}(n \lhd_{\dtheta} S_{\G}(m_{[0]_{[-1]}})) \lhd_{\dtheta} S^{2}_{\G}(m_{[-1]})  & \text{by } \eqref{eq:bc_condition} \\
& = (m_{[0]_{[0]}} \lhd_{\dtheta} S^{2}_{\G}(m_{[-1]_{(1)}}))(n \lhd_{\dtheta} S_{\G}(m_{[0]_{[-1]}})S^{2}_{\G}(m_{[-1]_{(2)}})) \\
& = (m_{[0]} \lhd_{\dtheta} S^{2}_{\G}(m_{[-1]}))n = \hat\gamma_{\theta}(m)n.
\end{align*}

\item Fix $m \in N$. By item~\ref{it:igualdad_general}, it follows
\[
\gamma_{\theta}(m)^{\op}_{[-1]} \odo \gamma_{\theta}(m)_{[0]} = \theta(\gamma_{\theta}(m)) = (S^{2}_{\Gc} \odo \gamma_{\theta})(\theta(m)) = S^{-2}_{\G}(m_{[-1]})^{\op} \odo \gamma_{\theta}(m_{[0]}),
\]
\[
\hat\gamma_{\theta}(m)^{\op}_{[-1]} \odo \hat\gamma_{\theta}(m)_{[0]} = \theta(\hat\gamma_{\theta}(m)) = (S^{-2}_{\Gc} \odo \hat\gamma_{\theta})(\theta(m)) = S^{2}_{\G}(m_{[-1]})^{\op} \odo \hat\gamma_{\theta}(m_{[0]}).
\]
Hence,
\begin{align*}
\gamma_{\theta}(\hat\gamma_{\theta}(m)) & = \hat\gamma_{\theta}(m)_{[0]} \lhd_{\dtheta} S^{-1}_{\G}(\hat\gamma_{\theta}(m)_{[-1]}) \\
& = (m_{[0]} \lhd_{\dtheta} S^{2}_{\G}(m_{[-1]_{(2)}})) \lhd_{\dtheta} S_{\G}(m_{[-1]_{(1)}}) \\
& = \cou_{\G}(m_{[-1]})m_{[0]}  = m,
\\\\
\hat\gamma_{\theta}(\gamma_{\theta}(m)) & = \gamma_{\theta}(m)_{[0]} \lhd_{\dtheta} S^{2}_{\G}(\gamma_{\theta}(m)_{[-1]}) \\
& = (m_{[0]} \lhd_{\dtheta} S^{-1}_{\G}(m_{[-1]_{(2)}})) \lhd_{\dtheta} m_{[-1]_{(1)}} \\
& = \cou_{H}(m_{[-1]})m_{[0]} = m.
\end{align*}

\item Fix $m,n \in N$. By item~\ref{it:lamd_ofside}, it follows
\begin{align*}
\gamma_{\theta}(mn) & = (m_{[0]}n_{[0]}) \lhd_{\dtheta} S^{-1}_{\G}(n_{[-1]}m_{[-1]}) \\
& = (m_{[0]}n_{[0]} \lhd_{\dtheta} S^{-1}_{\G}(m_{[-1]})) \lhd_{\dtheta} S^{-1}_{\G}(n_{[-1]}) \\
& = n_{[0]}\gamma_{\theta}(m) \lhd_{\dtheta} S^{-1}_{\G}(n_{[-1]}) \\
& = \gamma_{\theta}(m)\gamma_{\theta}(n).
\end{align*}
Hence, we also have
\[
\hat\gamma_{\theta}(m)\hat\gamma_{\theta}(n)  = \hat\gamma_{\theta}(\gamma_{\theta}(\hat\gamma_{\theta}(m)\hat\gamma_{\theta}(n))) = \hat\gamma_{\theta}(\gamma_{\theta}(\hat\gamma_{\theta}(m))\gamma_{\theta}(\hat\gamma_{\theta}(n))) = \hat\gamma_{\theta}(mn).
\]

\item\label{it:inversos} Fix $m \in N$. By item~\ref{it:igualdad_general}, we have
\[
\gamma_{\theta}(\gamma_{\theta}(m^{*})^{*})^{*} = \hat\gamma_{\theta}(\gamma_{\theta}(m^{*})) = m^{*}
\]
and
\[
\hat\gamma_{\theta}(\hat\gamma_{\theta}(m^{*})^{*}) = \gamma_{\theta}(\hat\gamma_{\theta}(m^{*}))^{*} = m.
\]
The claim follows directly from this.

\item By item~\ref{it:action_on_lambda}, we have
\begin{align*}
(\p(h,\cdot) \odo \id)\dtheta(\gamma_{\theta}(m)) & = \gamma_{\theta}(m) \lhd_{\dtheta} h = \gamma_{\theta}(m \lhd_{\dtheta} S^{2}_{\G}(h)) \\
& = (\p(S^{2}_{\G}(h),\cdot) \odo \gamma_{\theta})\dtheta(m) = (\p(h,\cdot) \odo \id)(\dS^{2}_{\G} \odo \gamma_{\theta})\dtheta(m)
\end{align*}
for all $h \in \pG$ and $m \in N$. The last implies $\dtheta\circ\gamma_{\theta} = (S^{2}_{\du{\G}^{\ops}} \odo \gamma_{\theta})\circ\dtheta$, and applying the inverse map $\hat\gamma_{\theta}$, we also have $\dtheta\circ\hat\gamma_{\theta} = (S^{-2}_{\du{\G}^{\ops}} \odo \hat\gamma_{\theta})\circ\dtheta$.

\item Consider the action of the multiplier Hopf \as-algebra $(\dpG,\dcom^{\co}_{\G})$ associated with $\theta$, i.e. the linear map $\lhd_{\theta}: N \odo \dpG \to N$ given by $m \lhd_{\theta} \ome = \p(m_{[-1]},\ome)m_{[0]}$ for all $m \in N$ and $\ome \in \dpG$. Hence, by direct computation
\begin{align*}
\gamma_{\theta}(m \lhd_{\theta} \ome) & = \p(m_{[-1]},\ome)\gamma_{\theta}(m_{[0]}) = \p(S^{-2}_{\G}(m_{[-1]}),\dS^{2}_{\G}(\ome))\gamma_{\theta}(m_{[0]}) \\
& = \p(\gamma_{\theta}(m)_{[-1]},\dS^{2}_{\G}(\ome))\gamma_{\theta}(m)_{[0]} = \gamma_{\theta}(m) \lhd_{\theta} \dS^{2}_{\G}(\ome)
\end{align*}
for $m \in N$ and $\ome \in \dpG$. Similarly, we can obtain the other equality.

\item Let $\mu$ be a $\dtheta$-invariant non-zero functional on $N$. We have $\mu(m \lhd_{\dtheta} h) = \mu(m)\cou_{\G}(h)$ for all $m \in N$, $h \in \pG$. Fix $m,n \in N$, thus $\mu(\gamma_{\theta}(m)) = \mu(m_{[0]} \lhd_{\dtheta} S^{-1}_{\G}(m_{[-1]})) = \mu(m)$ and by item~\ref{it:inversos} we also have $\mu(\hat\gamma_{\theta}(m)) = \mu(m)$. On the other hand, by item~\ref{it:lamd_ofside}, $\mu(n\gamma_{\theta}(m)) = \mu(m_{[0]}n \lhd_{\dtheta} S^{-1}_{\G}(m_{[-1]})) = \mu(mn)$ and similarly $\mu(\hat\gamma_{\theta}(n)m) = \mu(mn)$.
\end{enumerate}
\fin

\begin{defi}
Given a unital braided commutative Yetter--Drinfeld $\Gc$-\as-algebra $(N,\theta,\dtheta)$, the maps $\gamma_{\theta}$ and $\hat\gamma_{\theta}$ on $N$, defined in Theorem~\ref{theo:gamma_hatgamma}, are called {\em the canonical automorphisms of $(N,\theta,\dtheta)$}.
\end{defi}

\begin{exam}\label{ex:canonical_commutative_trivial}
Let $N$ be a unital commutative \as-algebra and $\theta: N \to \pG^{\op} \otimes N$ be any action of an algebraic quantum group of compact type $\Gc$ such that $(N,\theta,\tr_{\dG^{\ops}})$ is a unital braided commutative Yetter--Drinfeld $\Gc$-\as-algebra. The canonical automorphisms $\gamma_{\theta}$ and $\hat\gamma_{\theta}$ on $N$ are given by the trivial maps, indeed we have $\gamma_{\theta}(m) = m_{[0]} \lhd_{\tr_{\dG^{\ops}}}S^{-1}_{\G}(m_{[-1]})  = m$ and $\hat\gamma_{\theta}(m) = m_{[0]} \lhd_{\tr_{\dG^{\ops}}} S^{2}_{\G}(m_{[-1]}) = m$ for all $m \in N$, respectively.
\end{exam}

\begin{exam}\label{ex:canonical_commutative_trivial_dual}
Let $N$ be a unital commutative \as-algebra and $\dtheta: N \to \M(\dpG \otimes N)$ be any action of an algebraic quantum group of discrete type $\dG^{\ops}$ such that $(N,\tr_{\Gc},\dtheta)$ is a unital braided commutative Yetter--Drinfeld $\Gc$-\as-algebra. The canonical automorphisms $\gamma_{\tr_{\Gc}}$ and $\hat\gamma_{\tr_{\Gc}}$ on $N$ are given by the trivial maps, indeed we have $\gamma_{\tr_{\Gc}}(m) = m \lhd_{\dtheta} S^{-1}_{\G}(1) = m$ and $\hat\gamma_{\tr_{\G}}(m) = m \lhd_{\dtheta} S^{2}_{\G}(1) = m$ for all $m \in N$, respectively.
\end{exam}

\subsection{Duality for Yetter--Drinfeld \as-algebras}

Let $(N,\theta,\dtheta)$ be a unital braided commutative Yetter--Drinfeld $\Gc$-\as-algebra. Using the canonical automorphism $\hat\gamma_{\theta}$ on $N$, consider the $\hat\gamma_{\theta}$-opposite \as-algebra $N^{\op}_{\hat\gamma_{\theta}}$, i.e. the vector space $N$ with non-degenerate \as-algebra structure given by $m^{\op}n^{\op} := (nm)^{\op}$ and $(m^{\op})^{*} := \hat\gamma_{\theta}(m^{*})^{\op}$ for all $m,n \in N$. By Theorem~\ref{theo:gamma_hatgamma}, we have $\theta\circ\hat\gamma_{\theta} = (S^{-2}_{\Gc} \odo \hat\gamma_{\theta})\circ\theta$ and $\hat\theta\circ\hat\gamma_{\theta} = (S^{-2}_{\dG^{\ops}} \odo \hat\gamma_{\theta})\circ\hat\theta$, thus we can consider the conjugate actions $\theta^{\con}: N^{\op}_{\hat\gamma_{\theta}} \to \pG \odo N^{\op}_{\hat\gamma_{\theta}}$, $m^{\op} \mapsto ({}^{\op} \odo {}^{\op})\theta(m)$ and $\dtheta^{\con}: N^{\op}_{\hat\gamma_{\theta}} \to \M(\dpG^{\op} \odo N^{\op}_{\hat\gamma_{\theta}})$, $m^{\op} \mapsto ({}^{\op} \odo {}^{\op})\dtheta(m)$.

\begin{prop}\label{prop:dual_bc_yd}
Let $\G$ be an algebraic quantum group of compact type, $\theta: N \to \pG^{\op} \odo N$ be an action of $\Gc$ on $N$ and $\dtheta: N \to \M(\dpG \odo N)$ be an action of $\dGo$ on $N$. The following statements are equivalent:
\begin{enumerate}[label=\textup{(\roman*)}]
\item $(N,\theta,\dtheta)$ is a unital braided commutative Yetter--Drinfeld $\Gc$-\as-algebra;
\item $(N^{\op}_{\hat\gamma_{\theta}},\dtheta^{\con},\theta^{\con})$ is a unital braided commutative Yetter--Drinfeld $\dGco$-\as-algebra.
\end{enumerate}
Moreover, in that case, given a non-zero positive faithful functional $\mu$ on $N$, the following statements are equivalent:
\begin{enumerate}[label=\textup{(\roman*)}]
\item $\mu$ is a Yetter--Drinfeld integral for $(N,\theta,\dtheta)$;
\item $\mu^{\ops}$ is a Yetter--Drinfeld integral for $(N^{\op}_{\hat\gamma_{\theta}},\dtheta^{\con},\theta^{\con})$.
\end{enumerate}
\end{prop}
\pr
The first part follows from \cite[Proposition~3.8]{Ta22_1}. On the other hand, if $\mu$ is a Yetter--Drinfeld integral for $(N,\theta,\dtheta)$, we have
\begin{align*}
(\id \odo \mu^{\ops}){\theta^{\con}}(m^{\op}) & = ({}^{\op} \odo \mu)\dtheta(m) = \mu(m)1^{\op}_{\dpG} = \mu^{\ops}(m^{\op})1_{\pol(\dGco)}
\end{align*}
for all $m \in N$. Similarly, $(\id \odo \mu^{\ops})\theta^{\con}(m^{\op}) = \mu^{\ops}(m^{\op})1_{\pG}$ for all $m \in N$. This shows that $\mu^{\ops}$ is a Yetter--Drinfeld integral for $(N^{\op}_{\hat\gamma_{\theta}},\dtheta^{\con},\theta^{\con})$. The converse statement follows easily by replacing $(N,\theta,\dtheta)$ and $\mu$ by $(N^{\op}_{\hat\gamma_{\theta}},\dtheta^{\con},\theta^{\con})$ and $\mu^{\ops}$, respectively.
\fin

\begin{defi}
Let $(N,\theta,\dtheta,\mu)$ be a unital braided commutative measured Yetter--Drinfeld $\Gc$-\as-algebra. The unital braided commutative measured Yetter--Drinfeld $\dGco$-\as-algebra $(N^{\op}_{\hat\gamma_{\theta}},\dtheta^{\con},\theta^{\con},\mu^{\ops})$ is called {\em the dual conjugate of $(N,\theta,\dtheta,\mu)$}.
\end{defi}

\subsection{Examples}
Here, we collect the main examples of unital braided commutative measured Yetter--Drinfeld \as-algebras over algebraic quantum groups of compact type.

\begin{exam}[Finite transformation groupoids]\label{ex:bc_yd_finite_transformation_groupoid}
Let $\cdot : G \times X \to X$ be a left action of a finite group $G$ on a finite set $X$. Consider the braided commutative Yetter--Drinfeld $\Gc$-\as-algebra $(K(X),\theta,\dtheta)$ given by Example~\ref{ex:yd_transformation_groupoid}. By Example~\ref{ex:canonical_commutative_trivial}, the canonical automorphisms on $K(X)$ are $\gamma_{\theta} = \du\gamma_{\theta} = \id_{K(X)}$. Then by Proposition~\ref{prop:YD_integral_commutative}, given a non-zero $G$-invariant function $\nu:X \to \R^{+}_{0}$, the tuple $(K(X),\theta,\dtheta,\mu_{\nu})$ yields a unital braided commutative measured Yetter--Drinfeld $\Gc$-\as-algebra. On the other hand, by Proposition~\ref{prop:dual_bc_yd}, $(K(X)^{\op}_{\hat\gamma_{\theta}}=K(X),\dtheta^{\con}=\tr_{\dGco},\theta^{\con}=\theta,\mu^{\ops}_{\nu}=\mu_{\nu})$ yields a unital braided commutative measured Yetter--Drinfeld $\dGco$-\as-algebra.
\end{exam}

\begin{exam}[Compact transformation groupoids]\label{ex:bc_yd_compact_transformation_groupoid}
Let $G$ be a compact group. Denote by $\varphi$ the functional on $C(G)$ arising from the left Haar measure $h$ on $G$, i.e. $\varphi(f) = \int_{G}f \mathrm{d}h$ for all $f \in C(G)$, and by $\irr(G)$ the set of equivalence classes of irreducible continuous representations of $G$. Given an irreducible continuous representation $\pi: G \to \B(H_{\pi})$ on the vector space $H_{\pi}$ of dimension $\mathrm{d}_{\pi}$, fix an orthonormal basis $(e^{\pi}_{i})_{1 \leq i \leq \mathrm{d}_{\pi}}$ of $H_{\pi}$. For any pair $1 \leq i,j \leq \mathrm{d}_{\pi}$, consider the continuous functions $u^{\pi}_{ij}: G \to \C$, defined by $u^{\pi}_{ij}(g) = \langle \pi(g)e_{j},e_{i} \rangle$ for all $g \in G$ ,and called the {\em matrix coefficients functions of $\pi$ with respect to the basis $(e_{i})$}. With these notations, the vector subspace
\[
\mathscr{P}(G) := \mathrm{span}\{u^{\pi}_{ij} \in C(G) : \pi \in \irr(G), 1 \leq i,j \leq \mathrm{d}_{\pi}\}
\]
yields a unital dense sub-\as-algebra of $C(G)$, called the {\em \as-algebra of polynomial functions on $G$}. Then, considering the restrictions from $C(G)$ to $\mathscr{P}(G)$, of $\com$ and $\var$, the triplet $\G = (\mathscr{P}(G),\com,\varphi)$ yields an algebraic quantum group of compact type.

Now, let $X$ be a compact Hausdorff space equipped with a left action of $G$. There is a continuous action of the compact quantum group $(C(G),\com)$ on the C*-algebra of continuous functions on $X$, $C(X)$, defined by $\theta: C(X) \to C(G) \otimes C(X) \iso C(G \times X)$, $\theta(f): G \times X \to \C$, $(g,x) \mapsto f(g\cdot x)$ for every $f \in C(X)$. It can be shown that the vector subspace of $C(X)$
\[
K_{a}(X) := \mathrm{span}\left\{x \mapsto \int_{G} u^{\pi}_{ij}\theta(f)(-,x)\mathrm{d}h : \pi \in \irr(G), 1 \leq i,j \leq \mathrm{d}_{\pi}, f \in C(X) \right\}
\]
yields a unital sub-\as-algebra of $C(X)$, called the {\em Podle{\'s} \as-algebra of $C(X)$ through the action $\theta$}. Moreover, the restriction of $\theta$ to $K_{a}(X)$ yields an action $\theta: K_{a}(X) \to \mathscr{P}(G) \odo K_{a}(X)$ of the algebraic quantum group $\G$ on $K_{a}(X)$.

Consider the above action $\theta: K_{a}(X) \to \mathscr{P}(G)^{\op} \odo K_{a}(X)$ of $\Gc$ on $C(X)$ and the trivial action $\dtheta: K_{a}(X) \to \M(\du{\mathscr{P}(G)} \odo K_{a}(X))$ of $\dGo$ on $K_{a}(X)$, thus the triplet $(K_{a}(X),\theta,\dtheta)$ yields a unital braided commutative Yetter--Drinfeld $\Gc$-\as-algebra. Let $\nu: X \to \C$ be a left $G$-invariant measure on $X$, i.e. a positive measure on $X$ such that $\nu(g\cdot S) = \nu(S)$ for all $g \in G$ and all measurable set $S$ of $X$. Hence, the restriction to $K_{a}(X)$ of the non-zero positive functional $\mu_{\nu}: C(X) \to \C, f \mapsto \int_{X}f d\nu$ yields a Yetter--Drinfeld integral for $(K_{a}(X),\theta,\dtheta)$.
\end{exam}

\begin{exam}[Fell bundles]\label{ex:bc_yd_fell_bundle}
Let $\fell{A}$ be a separable Fell bundle over a discrete group $G$, $\G = (\C[G],\dcom,\dvar)$ be the algebraic quantum group of compact type associated with $G$ and $\dG = (K(G),\com^{\co},\var\circ S)$ be the dual algebraic quantum group of $\G$. By an adapted version of \cite[Proposition~3.1]{KMQW10}, there is an action of $\Gc = (\C[G]^{\op},\dcom^{\op},\dvar^{\ops})$ on the unital \as-algebra of continuous compactly supported cross sections of $\fell{A}$, $\theta_{\fell{A}}: \Gamma_{c}(G,\fell{A}) \to \C[G]^{\op} \odo \Gamma_{c}(G,\fell{A})$, defined by $\theta_{\fell{A}}(\iota(a)) = \lambda_{g^{-1}} \odo \iota(a)$ for each $g \in G$ and $a \in \fell{A}_{g}$. Recall that the \as-algebra structure in $\Gamma_{c}(G,\fell{A}) \subset C_{c}(G,\fell{A})$ is given by
\[
(f \cdot f')(g) = \sum_{h\in G}f(h)f'(h^{-1}g), \quad \text{ and } \quad f^{*}(g) = f(g^{-1})^{*},
\]
for all $f,f \in \Gamma_{c}(G,\fell{A})$ and $g \in G$; and $\iota: \fell{A} \to \Gamma_{c}(g,\fell{A})$ denotes the canonical inclusion defined by $\iota(a)(h) = \delta_{g,h}a$ for all $a \in \fell{A}_{g}$ and $g,h \in G$.

Given an element $g\in G$, denote by $\fell{A}^{\ad(g^{-1})}$ the separable Fell bundle over $G$ defined by
\[
(\fell{A}^{\ad(g^{-1})})_{g'} := \fell{A}_{g^{-1}g'g}
\]
for all $g' \in G$. Suppose now that there exists a family of Fell bundle isomorphisms, in the sense of \cite{AE01}, $\rho = \{\rho_{g}:\fell{A} \to \fell{A}^{\ad(g^{-1})}\}_{g\in G}$ such that
\begin{enumerate}
\item $\rho_{e} = \id_{\fell{A}}$;
\item $\rho_{g}\circ\rho_{g'} = \rho_{g'g}$ for all $g,g' \in G$;
\item for any $g,h \in G$, we have $ab = \rho_{g^{-1}}(b)a$ for all $a \in \fell{A}_{g}$ and $b \in \fell{A}_{h}$.
\end{enumerate}
Hence, the linear map $\dtheta_{\fell{A},\rho}: \Gamma_{c}(G,\fell{A}) \to \M(K(G) \odo \Gamma_{c}(G,\fell{A}))$, defined by
\[
\dtheta_{\fell{A},\rho}(\iota(a)) = \sum_{g \in G} \delta_{g} \odo \iota(\rho_{g}(a))
\]
for all $a \in \fell{A}$, gives an action of $\dGo = (K(G),\com,\varphi)$ on the \as-algebra $\Gamma_{c}(G,\fell{A})$ such that the triplet $(\Gamma_{c}(G,\fell{A}),\theta_{\fell{A}},\dtheta_{\fell{A},\rho})$ yields a braided commutative Yetter--Drinfeld $\Gc$-\as-algebra. In this example, the canonical automorphisms on $\Gamma_{c}(G,\fell{A})$ are given by $\gamma_{\theta_{\fell{A}}}(\iota(a)) = \iota(\rho_{g}(a))$ and $\hat\gamma_{\theta_{\fell{A}}}(\iota(a)) = \iota(\rho^{-1}_{g}(a))$ for all $a \in \fell{A}_{g}$ and $g \in G$.

Let $\nu: \fell{A}_{e} \to \C$ be a faithful state on the unital C*-algebra $\fell{A}_{e}$. Then, the functional $\mu_{\nu}$ on $\Gamma_{c}(G,\fell{A})$, defined by $f \mapsto \nu(f(e))$, yields a Yetter--Drinfeld integral for $(\Gamma_{c}(G,\fell{A}),\theta_{\fell{A}},\dtheta_{\fell{A},\rho})$. Reciprocally, any Yetter--Drinfeld integral for $(\Gamma_{c}(G,\fell{A}),\theta_{\fell{A}},\dtheta_{\fell{A},\rho})$ arises in this way.
\end{exam}

\begin{exam}[Trivial examples]\label{ex:bc_yd_trivial}
Let $\G$ be an algebraic quantum group of compact type and $N$ be a unital \as-algebra. Consider the trivial actions $\tr_{\Gc}: N \to \pG^{\op} \odo N$, $m \mapsto 1^{\op}_{\pG} \odo m$ and $\tr_{\dGo}: N \to \M(\dpG \odo N)$, $m \mapsto 1_{\M(\dpG)} \odo m$. The triplet $(N,\tr_{\Gc},\tr_{\dGo})$ yields a unital Yetter--Drinfeld $\Gc$-\as-algebra, which it is braided commutative if an only if the \as-algebra $N$ is commutative. Moreover, the canonical automorphims on $N$ are given by the trivial maps and, any non-zero positive faithful functional $\mu$ on $N$ yields a Yetter--Drinfeld integral for $(N,\tr_{\Gc},\tr_{\dGo})$.
\end{exam}

\begin{exam}[Canonical examples]\label{ex:bc_yd_canonical_aqg}
Let $\G$ be an algebraic quantum group of compact type. Consider $N = \pG$, $\theta = (S^{-1}_{\G} \odo \id)\circ\Sigma\circ\com_{\G}$ and $\dtheta = \ad_{\Sigma(U^{*})}$. The triplet $(N,\theta,\dtheta)$ yields a unital braided commutative Yetter--Drinfeld $\Gc$-\as-algebra and its canonical automorphisms are given by $\gamma_{\theta} = S^{-2}_{\G}$ and $\hat\gamma_{\theta} = S^{2}_{\G}$ (see Corollary~\ref{coro:canonical_example} in Example \ref{ex:bc_yd_quotient_type}).
\end{exam}

\begin{exam}[Quotient type coideals]\label{ex:bc_yd_quotient_type}
Let $\G = (\pG,\com_{\G},\varphi_{\G})$ be an algebraic quantum group of compact type and $\HH = (\pol(\HH),\com_{\HH},\varphi_{\HH})$ be an algebraic quantum subgroup of $\G$ defined by the surjective \as-homomorphism $\pi : \pG \to \pol(\HH)$. The linear maps ${}^{\pi}\lambda : \pG \to \pol(\HH) \odo \pG$, $h \mapsto (\pi \odo \id)\com_{\G}(h)$ and $\lambda^{\pi} : \pG \to \pG \odo \pol(\HH)$, $h \mapsto (\id \odo \pi)(\com_{\G}(h))$, define a left and right coaction of the unital Hopf \as-algebra $(\pol(\HH),\com_{\HH})$ on $\pG$, respectively, such that $(\id \odo \com_{\G}){}^{\pi}\lambda = ({}^{\pi}\lambda \odo \id)\com_{\G}$ and $(\com_{\G} \odo \id)\lambda^{\pi} = (\id \odo \lambda^{\pi})\com_{\G}$. Consider the following two sub-\as-algebras of $\pG$
\[
\pol(\G/\HH) := \{h \in \pG : \lambda^{\pi}(h) = h \odo 1_{\pol(\HH)} \}, \quad \pol(\HH\backslash\G) := \{h \in \pG : {}^{\pi}\lambda(h) = 1_{\pol(\HH)} \odo h \},
\]
called {\em the quotient spaces related to the ``inclusion'' of algebraic quantum groups $\HH < \G$}, respectively.

\begin{prop}\label{prop:YD_quotient}
Consider the following algebraic quantum group actions $\com_{\G}: \pG \to \pG \odo \pG$ and $\ad_{\Sigma(U^{*})}: \pG \to \M(\dpG \odo \pG)$ related to $\G$. For each $m \in \pol(\HH\backslash\G)$, we have $\com_{\G}(m) \in \pol(\HH\backslash\G) \odo \pG$ and $\ad_{\Sigma(U^{*})}(m) \in \M(\dpG \odo \pol(\HH\backslash\G))$. Moreover, if $\delta: \pol(\HH\backslash\G) \to \pol(\HH\backslash\G) \odo \pG$ denotes the restriction of $\com_{\G}$ to $\pol(\HH\backslash\G)$ and $\dtheta: \pol(\HH\backslash\G) \to \M(\du{\pol(\G}) \odo \pol(\HH\backslash\G))$ denotes the restriction of $\ad_{\Sigma(U^{*})}$ to $\pol(\HH\backslash\G)$, then $(\pol(\HH\backslash\G),\theta_{\delta},\dtheta)$ yields a unital braided commutative Yetter--Drinfeld $\Gc$-\as-algebra with canonical automorphisms given by $\gamma_{\theta_{\delta}} = S^{-2}_{\G}|_{\pol(\HH\backslash\G)}$ and $\hat\gamma_{\theta_{\delta}} = S^{2}_{\G}|_{\pol(\HH\backslash\G)}$. Here, we are using the notation $\theta_{\delta} :=(S^{-1}_{\G}\odo\id)\circ\Sigma\circ\delta$.
\end{prop}
\pr
Fix $m \in \pol(\HH\backslash\G)$, we have
\[
({}^{\pi}\lambda \odo \id)(\com_{\G}(m)) = (\pi \odo \id \odo \id)(\com_{\G} \odo \id)\com_{\G}(m) = (\id \odo \com_{\G}){}^{\pi}\lambda(m) = 1_{\pol(\HH)} \odo \com_{\G}(m)
\]
then $\com_{\G}(m) \in \pol(\HH\backslash\G) \odo \pG$. On the other hand, because $(\com_{\G} \odo \id)(U) = U_{13}U_{23}$, thus
\begin{align*}
(\id \odo {}^{\pi}\lambda)(\ad_{\Sigma(U^{*})}(m)) & = (\id \odo \pi \odo \odo \id)(\id \odo \com_{\G})(\Sigma(U^{*})(1 \odo m)\Sigma(U)) \\
& = (\id \odo \pi \odo \id)(\Sigma(U^{*})_{13}\Sigma(U^{*})_{12}(1 \odo \com_{\G}(m))\Sigma(U)_{12}\Sigma(U)_{13}) \\
& = \Sigma(U^{*})_{13}W^{*}_{12}(1 \odo {}^{\pi}\lambda(m))W_{12}\Sigma(U)_{13} \\
& = \Sigma(U^{*})_{13}(1 \odo 1 \odo m)\Sigma(U)_{13} \\
& = \ad_{\Sigma(U^{*})}(m)_{13},
\end{align*}
where $W := (\id \odo \pi)(\Sigma(U))$, so $\ad_{\Sigma(U^{*})}(m) \in \M(\dpG \odo \pol(\HH\backslash\G))$. Equivalently, note that the dual action $\lhd_{\ad_{\Sigma(U^{*})}}: \pG \odo \pG \to \pG$ is given by $h' \lhd_{\ad_{\Sigma(U^{*})}} h = S_{\G}(h_{(1)})h'h_{(2)}$ for all $h,h' \in \pG$, then it is straightforward to show that $m \lhd_{\ad_{\Sigma(U^{*})}} h \in \pol(\HH\backslash\G)$ for every $h \in \pG$, then necessarily $\ad_{\Sigma(U^{*})}(m) \in \M(\dpG \odo \pol(\HH\backslash\G))$.

By simple computation, we have
\begin{align*}
\ad(\Sigma(U^{*})_{13})(\dtheta \odo \id)(\delta(m)) & = \Sigma(U^{*})_{13}\Sigma(U^{*})_{12}(1 \odo \com_{\G}(m))\Sigma(U)_{12}\Sigma(U)_{13} \\
& = (\id \odo \delta)(\dtheta(m)),
\end{align*}
for all $m \in \pol(\HH\backslash\G)$, which is equivalent to have
\[
\delta(m \lhd_{\dtheta} h) = \delta(S_{\G}(h_{(1)})mh_{(2)}) = (m_{(1)} \lhd_{\dtheta} h_{(2)}) \odo S_{\G}(h_{(1)})m_{(2)}h_{(3)},
\]
for all $m \in \pol(\HH\backslash\G)$ and $h \in \pG$. Additionally $n_{(1)}(m \lhd_{\dtheta} n_{(2)}) = n_{(1)}S_{\G}(n_{(2)})mn_{(3)} = mn$, then $(\pol(\HH\backslash\G),\lhd_{\dtheta},\delta)$ is a unital braided commutative (right-right) Yetter--Drinfeld \as-algebra over $(\pG,\com_{\G})$ or equivalently $(\pol(\HH\backslash\G),\delta,\dtheta)$ yields a unital braided commutative (right-left) Yetter--Drinfeld \as-algebra over the canonical pairing $\p: \pG \times \dpG \to \C$, $\p(h,\ome) = \ome(h)$. By Proposition~\ref{prop:bc_yd_equivalence}, $(\pol(\HH\backslash\G),\theta_{\delta},\dtheta)$ is a unital braided commutative Yetter--Drinfeld $\Gc$-\as-algebra.

Using the Sweedler type leg notation, we have $\theta_{\delta}(m) = S^{-1}_{\G}(m_{(2)})^{\op} \odo m_{(1)}$ for all $m \in \pol(\HH\backslash\G)$ and on the other hand, the action $\lhd_{\dtheta}: \pol(\HH\backslash\G) \otimes \pG \to \pol(\HH\backslash\G)$ is given by $m \lhd_{\dtheta} h = S_{\G}(h_{(1)})mh_{(2)}$ for all $m \in \pol(\HH\backslash\G)$, $h \in \pG$. Thus by Theorem~\ref{theo:gamma_hatgamma}, the canonical automorphisms $\gamma_{\theta_{\delta}}$ and $\hat\gamma_{\theta_{\delta}}$ on $\pol(\HH\backslash\G)$ of $(\pol(\HH\backslash\G),\theta_{\delta},\dtheta)$ are given by $\gamma_{\theta_{\delta}}(m) = m_{(1)} \lhd_{\dtheta} S^{-2}_{\G}(m_{(2)}) = S^{-1}_{\G}(m_{(2)})m_{(1)}S^{-2}_{\G}(m_{(3)}) = S^{-2}_{\G}(m)$ and $\hat\gamma_{\theta_{\delta}}(m) = m_{(1)} \lhd_{\dtheta} S_{\G}(m_{(2)}) = S^{2}_{\G}(m_{(3)})m_{(1)}S_{\G}(m_{(2)}) = S^{2}_{\G}(m)$ for all $m \in \pol(\HH\backslash\G)$, respectively.
\fin

\begin{rema}
The braided commutative Yetter--Drinfeld $\Gc$-\as-algebra $(\pol(\HH\backslash\G),\theta_{\delta},\dtheta)$ defined above is called a {\em quotient type left coideal of $\G$ by $\HH$}.
\end{rema}

\begin{coro}\label{coro:canonical_example}
Let $\G = (\pG,\com_{\G},\varphi_{\G})$ be an algebraic quantum group of compact type. The triplet $(\pG,\theta=(S^{-1}_{\G}\odo\id)\circ\Sigma\circ\com_{\G},\dtheta=\ad_{\Sigma(U^{*})})$ yields a unital braided commutative Yetter--Drinfeld $\Gc$-\as-algebra. Its canonical automorphisms are given by $\gamma_{\theta} = S^{-2}_{\G}$ and $\hat\gamma_{\theta} = S^{2}_{\G}$.
\end{coro}
\pr
Considering the canonical surjective \as-homomorphism $\cou_{\G} : \pG \to \C$, the trivial algebraic quantum group $\HH_{\bullet}=(\C,\id_{\C},\id_{\C})$ is an algebraic quantum subgroup of $\G$. Because ${}^{\cou_{\G}}\lambda = \id_{\pG}$, we obtain $\pol(\HH_{\bullet}\backslash \G) = \pG$ as quotient space. The statement follows from Proposition~\ref{prop:YD_quotient}.
\fin
\end{exam}

\section{Measured multiplier Hopf \as-algebroids arising from Yetter--Drinfeld \as-algebras}\label{sec:mhad_yd}

Let $\G=(\pG,\com_{\G},\varphi_{\G})$ be an algebraic quantum group of compact type and $(N,\theta,\dtheta)$ be a unital braided commutative Yetter--Drinfeld $\Gc$-\as-algebra with canonical automorphism $\gamma_{\theta}$ and $\hat\gamma_{\theta}$. By Proposition~\ref{prop:bc_yd_equivalence}, $(N,\lhd_{\dtheta},\theta)$ yields a unital braided commutative (right-left) Yetter--Drinfeld $(\pG,\com_{\G})$-\as-algebra. Using the right Hopf \as-algebra action $\lhd_{\dtheta}$, consider the unital smash product \as-algebra $A:= \pG \sm_{\dtheta} N$, with algebraic structure given by
\[
(h \sm m)(g \sm n) = hg_{(1)} \sm (m \lhd_{\dtheta} g_{(2)})n,
\]
\[
(h \sm m)^{*} = h^{*}_{(1)} \sm (m^{*} \lhd_{\dtheta} h^{*}_{(2)}) = (h_{(1)})^{*} \sm (m \lhd_{\dtheta} S^{-1}_{\G}(h_{(2)}))^{*}
\]
for every $h,g \in \pG$ and $m, n \in N$. Observe that, we also have
\[
(1_{\pG} \sm m)(h \sm 1_{N}) = h_{(1)} \sm (m \lhd_{\dtheta} h_{(2)}), \;\;\; (1_{\pG} \sm (m \lhd_{\dtheta} S^{-1}_{\G}(h_{(2)})))(h_{(1)} \sm 1_{N}) = h \sm m
\]
for all $h \in \pG$, $m \in N$. The vector subspaces
\[
\text{ and } \quad B := \{1_{\pG} \sm m \in A : m \in N\} \quad C := \{m_{[-1]} \sm m_{[0]} \in A : m \in N\}
\]
yield two sub-\as-algebras of $A$. Indeed, we have
\[
(1_{\pG} \sm m)(1_{\pG} \sm n) = 1_{\pG} \sm mn, \quad (1_{\pG} \sm m)^* = 1_{\pG} \sm m^{*};
\]
and
\begin{align*}
(m_{[-1]} \sm m_{[0]})(n_{[-1]} \sm n_{[0]}) & = m_{[-1]}n_{[-1]_{(1)}} \sm (m_{[0]} \lhd_{\dtheta} n_{[-1]_{(2)}})n_{[0]} \\
& = m_{[-1]}n_{[-1]} \sm n_{[0]}m_{[0]}, & \text{by } \eqref{eq:bc_condition} \\
& = (nm)_{[-1]} \sm (nm)_{[0]}
\\\\
(m_{[-1]} \sm m_{[0]})^* & = (m_{[-1]_{(1)}})^{*} \sm (m_{[0]} \lhd_{\dtheta} S^{-1}_{\G}(m_{[-1]_{(2)}}))^{*} \\
& = S^{2}_{\G}((m^{*})_{[-1]}) \sm \hat\gamma_{\theta}((m^{*})_{[0]}) \\
& = \hat\gamma_{\theta}(m^{*})_{[-1]} \# \hat\gamma_{\theta}(m^{*})_{[0]} & \text{by } \ref{itpr:igualdad_general}
\end{align*}
for all $m, n \in N$. Moreover, $B$ and $C$ commute elementwise inside $A$, indeed we have
\begin{align*}
(1_{\pG} \sm m)(n_{[-1]} \sm n_{[0]}) & = n_{[-1]_{(1)}} \sm (m \lhd_{\dtheta} n_{[-1]_{(2)}})n_{[0]} \\
& = n_{[-1]_{(1)}} \sm (m \lhd_{\dtheta} n_{[0]_{[-1]}})n_{[0]_{[0]}} \\
& = n_{[-1]} \sm n_{[0]}m = (n_{[-1]} \sm n_{[0]})(1_{\pG} \sm m) & \text{by } \eqref{eq:bc_condition}
\end{align*}
for every $m, n \in N$. The conditions above can be summarize in the next proposition.

\begin{prop}\label{prop:multiplier_algebra}
Let $\G=(\pG,\com_{\G},\varphi_{\G})$ be an algebraic quantum group of compact type and $(N,\theta,\dtheta)$ be a unital braided commutative Yetter--Drinfeld $\Gc$-\as-algebra with canonical automorphism $\gamma_{\theta}$ and $\hat\gamma_{\theta}$. Consider the unital \as-algebra $A = \pG \sm_{\dtheta} N$. Thus, the linear maps
\begin{align*}\label{maps:alpha_beta}
\begin{array}{lccc}
\alp: & N & \to & A \\
& m & \mapsto & 1_{\pG} \sm m
\end{array},
\qquad
\begin{array}{lccc}
\beta: & N & \to & A \\
& m & \mapsto & m_{[-1]} \sm m_{[0]}
\end{array}
\end{align*}
are injective and satisfy the following conditions
\begin{enumerate}[label=\textup{(C\arabic*)}]
\item\label{eq:C1} $\alp(m)\alp(n) = \alp(mn)$ and $\alp(m)^{*} = \alp(m^{*})$,
\item\label{eq:C2} $\beta(m)\beta(n) = \beta(nm)$ and $\beta(m)^{*} = \beta(\hat\gamma_{\theta}(m^{*})) = \beta(\gamma_{\theta}(m)^{*})$,
\item\label{eq:C3} $\alp(m)\beta(n) = \beta(n)\alp(m)$
\end{enumerate}
for every $m, n \in N$. Moreover, by Theorem~\ref{theo:gamma_hatgamma}, it is straightforward to show the following equalities

\begin{enumerate}[label=\textup{(I\arabic*)}]
\item\label{eq:I1} $\alp(m)(h \sm 1) = (h_{(1)} \sm 1)\alp(m \lhd_{\dtheta} h_{(2)})$,
\item\label{eq:I2} $(h \sm 1)\alp(m) = \alp(m \lhd_{\dtheta} S^{-1}_{\G}(h_{(2)}))(h_{(1)} \sm 1)$,
\item\label{eq:I3} $\alp(m) = (S_{\G}(m_{[-1]}) \sm 1)\beta(m_{[0]}) = \beta(\hat\gamma_{\theta}(m_{[0]}))(S_{\G}(m_{[-1]}) \sm 1)$,
\item\label{eq:I4} $\alp(\gamma_{\theta}(m)) = (S^{-1}_{\G}(m_{[-1]}) \sm 1)\beta(\gamma_{\theta}(m_{[0]})) = \beta(m_{[0]})(S^{-1}_{\G}(m_{[-1]}) \sm 1)$,
\item\label{eq:I5} $\alp(\hat\gamma_{\theta}(m)) = (S^{3}_{\G}(m_{[-1]}) \sm 1)\beta(\hat\gamma(m_{[0]})) = \beta(\hat\gamma^{2}_{\theta}(m_{[0]}))(S^{3}_{\G}(m_{[-1]}) \sm 1)$,
\item\label{eq:I6} $\beta(m) = (m_{[-1]} \sm 1)\alp(m_{[0]}) = \alp(\gamma_{\theta}(m_{[0]}))(m_{[-1]} \sm 1)$,
\item\label{eq:I7} $\beta(\gamma_{\theta}(m)) = (S^{-2}_{\G}(m_{[-1]}) \sm 1)\alp(\gamma_{\theta}(m_{[0]})) = \alp(\gamma^{2}_{\theta}(m_{[0]}))(S^{-2}_{\G}(m_{[-1]}) \sm 1)$
\item\label{eq:I8} $\beta(\hat\gamma_{\theta}(m)) = (S^{2}_{\G}(m_{[-1]}) \sm 1)\alp(\hat\gamma_{\theta}(m_{[0]})) = \alp(m_{[0]})(S^{2}_{\G}(m_{[-1]}) \sm 1)$
\end{enumerate}
for all $m \in N$, $h \in H$. In particular, besides the usual representation of the smash product \as-algebra $A = \{(h \sm 1_{N})\alp(m) : h \in \pG, m \in N \}$, we also have
\[
A = \{ (h \sm 1_{N})\beta(m) : h \in \pG, m \in N \} = \{ \beta(m)(h \sm 1_{N}) : h \in \pG, m \in N \}.
\]
\end{prop}

\begin{rema}
Sometimes, the linear maps defined in the proposition above will be denoted by $\alp_{\dtheta}$ and $\beta_{\theta}$, respectively, in order to emphasize its relation with the actions coming from the Yetter--Drinfeld structure $(N,\theta,\dtheta)$.
\end{rema}

The next lemma will be important for the next theorem.

\begin{lemm}\label{lem:technical_beta}
It holds
\begin{enumerate}[label=\textup{(\alph*)}]
\item \begin{equation}\label{eq:YD_condition_smash_1}
\beta(m)(h \sm 1) = (h_{(2)} \sm 1)\beta(m \lhd_{\dtheta} h_{(1)}),
\end{equation}
\begin{equation}\label{eq:YD_condition_smash_2}
(h \sm 1)\beta(m) = \beta(m \lhd_{\dtheta} S_{\G}(h_{(1)}))(h_{(2)} \sm 1);
\end{equation}

\item
\begin{align}\label{eq:beta_antipode}
\beta(m \lhd_{\dtheta} h_{(2)})(S^{-1}_{\G}(h_{(1)}) \sm 1) = (S^{-1}_{\G}(h) \sm 1)\beta(m)
\end{align}
\begin{align}\label{eq:beta_hatgamma_antipode}
\beta(\hat\gamma_{\theta}(m \lhd_{\dtheta} h_{(2)}))(S_{\G}(h_{(1)}) \sm 1) = (S_{\G}(h) \sm 1)\beta(\hat\gamma_{\theta}(m))
\end{align}
\end{enumerate}
for all $m \in N$ and $h \in \pG$.
\end{lemm}
\pr
Fix $m \in N$ and $h \in H$.
\begin{enumerate}[label=\textup{(\alph*)}]
\item The Yetter--Drinfeld condition
\[
\theta(m \lhd_{\dtheta} h) = (S^{-1}_{\G}(h_{(1)})m_{[-1]}h_{(1)})^{\op} \odo (m_{[0]} \lhd_{\dtheta} h_{(2)})
\]
is equivalent to have the condition
\begin{equation}\label{eq:YD_condition_classic}
m_{[-1]}h_{(1)} \odo (m_{[0]} \lhd_{\dtheta} h_{(2)}) = h_{(2)}(m \lhd_{\dtheta} h_{(1)})_{[-1]} \odo (m \lhd_{\dtheta} h_{(1)})_{[0]}.
\end{equation}
Thus
\begin{align*}
\beta(m)(h \sm 1) & = m_{[-1]}h_{(1)} \sm (m_{[0]} \lhd_{\dtheta} h_{(2)}) \\
& = h_{(2)}(m \lhd_{\dtheta} h_{(1)})_{[-1]} \sm (m \lhd_{\dtheta} h_{(1)})_{[0]} \\
& = (h_{(2)} \sm 1)\beta(m \lhd_{\dtheta} h_{(1)})
\end{align*}
and
\begin{align*}
\beta(m \lhd_{\dtheta} S_{\G}(h_{(1)}))(h_{(2)} \sm 1) & = (h_{(3)} \sm 1)\beta((m \lhd_{\dtheta} S_{\G}(h_{(1)})) \lhd_{\dtheta} h_{(2)}) \\
& = (h \sm 1)\beta(m).
\end{align*}

\item Follows directly from the item above.
\end{enumerate}
\fin

\begin{rema}
Observe that if $(N,\theta,\dtheta,\mu)$ is a unital braided commutative measured Yetter--Drinfeld $\Gc$-\as-algebra, i.e. since the functional $\mu$ verify $(\id \odo \mu)\theta(m) = \mu(m)1^{\op}_{\pG}$ and $(\id \odo \mu)\dtheta(m) = \mu(m)1_{\dpG}$ for all $m \in N$, then we have $\mu(m \lhd_{\theta} \ome) = \dcou_{\G}(\ome)\mu(m)$ and $\mu(m \lhd_{\dtheta} h) = \cou_{\G}(h)\mu(m)$ for every $m \in N$, $\ome \in \dpG$ and $h \in \pG$.
\end{rema}

The next theorem is the involutive non necessarily finite case of \cite[Theorem~4.1]{BM02}, \cite[Examples~2.6, 3.4, 5.5]{T16} and \cite[Section~5.4]{T17}.

\begin{theo}\label{theo:MHAd_YD}
Let $\G=(\pG,\com_{\G},\varphi_{\G})$ be an algebraic quantum group of compact type and $(N,\theta,\dtheta,\mu)$ be a unital braided commutative measured Yetter--Drinfeld $\Gc$-\as-algebra. Consider the notations of Proposition~\ref{prop:multiplier_algebra} and take the following linear maps
\[
\begin{array}{lccc}
t_{B}: & B & \to & C \\
& \alp(m) & \mapsto & \beta(m)
\end{array},
\qquad
\begin{array}{lccc}
t_{C}: & C & \to & B \\
& \beta(m) & \mapsto & \alp(\gamma_{\theta}(m))
\end{array},
\]
\[
\begin{array}{lccc}
\com_{B}: & A & \to & A \rtak{}{B} A \\
& h \sm m & \mapsto & (h_{(1)} \sm 1_{N}) \rtak{}{B} (h_{(2)} \sm m)
\end{array},
\,
\begin{array}{lccc}
\com_{C}: & A & \to & A \ltak{}{C} A \\
& h \sm m & \mapsto & (h_{(1)} \sm 1_{N}) \ltak{}{C} (h_{(2)} \sm m)
\end{array},
\]
\[
\begin{array}{lccc}
S: & A & \to & A \\
& h \sm m & \mapsto & \beta(\hat\gamma_{\theta}(m))(S_{\G}(h) \sm 1_{N})
\end{array},
\]
\[
\begin{array}{lccc}
\cou_{B}: & A & \to & B \\
& \alp(m)(h \sm 1_{N}) & \mapsto & \alp(m \lhd_{\dtheta} h)
\end{array},
\qquad
\begin{array}{lccc}
{}_{C}\cou: & A & \to & C \\
& (h \sm 1_{N})\beta(m) & \mapsto & \beta(m \lhd_{\dtheta} S_{\G}(h))
\end{array},
\]
\[
\begin{array}{lccc}
\mu_{B}: & B & \to & \C \\
& \alp(m) & \mapsto & \mu(m)
\end{array},
\qquad
\begin{array}{lccc}
\mu_{C}: & C & \to & \C \\
& \beta(m) & \mapsto & \mu(m)
\end{array},
\]
and
\[
\begin{array}{lccc}
{}_{B}\psi_{B}: & A & \to & B \\
& \alp(m)(h \sm 1_{N}) & \mapsto & \varphi_{\G}(h)\alp(m)
\end{array},
\qquad
\begin{array}{lccc}
{}_{C}\phi_{C}: & A & \to & C \\
& (h \sm 1_{N})\beta(m) & \mapsto & \varphi_{\G}(h)\beta(m)
\end{array}.
\]
Hence, the collection
\[
\mathcal{A}(N,\theta,\dtheta,\mu) := (A,B,C,t_{B},t_{C},\com_{B},\com_{C},\mu_{B},\mu_{C},{}_{B}\psi_{B},{}_{C}\phi_{C})
\]
yields a unital measured multiplier Hopf \as-algebroid, called the {\em measured multiplier Hopf \as-algebroid associated with the braided commutative Yetter--Drinfeld $\Gc$-\as-algebra $(N,\theta,\dtheta)$ and the Yetter--Drinfeld integral $\mu$}.
\end{theo}

\pr
For sake of completeness, we check each condition in the definition of a measured multiplier Hopf \as-algebroid. We have
\begin{enumerate}[label=\textup{(M\arabic*)}]
\item The condition follows mainly from Proposition~\ref{prop:multiplier_algebra}. Because $A$ is a unital \as-algebra, thus the canonical modules ${}_{B}A$, $A_{B}$, ${}_{C}A$ and $A_{C}$ are evidently idempotent.

\item Fix $m \in N$. Because $t_{B}(\alp(m)^{*}) = \beta(m^{*})$ and $t_{C}(\beta(m)^{*}) = t_{C}(\beta(\hat\gamma_{\theta}(m^{*}))) = \alp(m)^{*}$, thus $t_{C}(t_{B}(\alp(m)^{*})^{*}) = t_{C}(\beta(\hat\gamma_{\theta}(m))) = \alp(m)$. Equivalently, we have $t_{B}(t_{C}(\beta(m)^{*})^{*}) = t_{B}(\alp(m)) = \beta(m)$.

\item \textbf{The left and right comultiplications:} For each $h,g \in \pG$ and $m,n \in N$, we have
\begin{align*}
\com_{B}((h \sm m)(g \sm n)) & = \com_{B}(hg_{(1)} \sm (m \lhd_{\dtheta} g_{(2)})n) \\
& = (h_{(1)}g_{(1)} \sm 1_{N}) \rtak{}{B} (h_{(2)}g_{(2)} \sm (m \lhd_{\dtheta} g_{(3)})n) \\
& = (h_{(1)} \sm 1_{N})(g_{(1)} \sm 1_{N}) \rtak{}{B} (h_{(2)} \sm m)(g_{(2)} \sm n) \\
& = \com_{B}(h \sm m)\com_{B}(g \sm n),
\end{align*}
and
\begin{align*}
\com_{C}((h \sm m)(g \sm n)) & = \com_{C}(hg_{(1)} \sm (m \lhd_{\dtheta} g_{(2)})n) \\
& = (h_{(1)}g_{(1)} \sm 1_{N}) \ltak{}{C} (h_{(2)}g_{(2)} \sm (m \lhd_{\dtheta} g_{(3)})n) \\
& = (h_{(1)} \sm 1_{N})(g_{(1)} \sm 1_{N}) \ltak{}{C} (h_{(2)} \sm m)(g_{(2)} \sm n) \\
& = \com_{C}(h \sm m)\com_{C}(g \sm n).
\end{align*}
In particular,
\begin{align*}
\com_{B}(\alp(m)(h \sm 1_{N}))
& = (h_{(1)} \sm 1_{N}) \rtak{}{B} \alp(m)(h_{(2)} \sm 1_{N}),
\end{align*}
\begin{align*}
\com_{C}(\alp(m)(h \sm 1_{N}))
& = (h_{(1)} \sm 1_{N}) \ltak{}{C} \alp(m)(h_{(2)} \sm 1_{N}),
\end{align*}
and
\begin{align*}
\com_{B}(\beta(m)) & = \com_{B}(\alp(\gamma_{\theta}(m_{[0]}))(m_{[-1]} \sm 1_{N})) \\
& = (m_{[-1]} \sm 1_{N}) \rtak{}{B} \beta(m_{[0]}) \\
& = (m_{[-1]} \sm 1_{N}) \rtak{}{B} t_{B}(\alp(m_{[0]})) \\
& = (m_{[-1]} \sm 1_{N})\alp(m_{[0]}) \rtak{}{B} (1_{\pG} \sm 1_{N}) \\
& = \beta(m) \rtak{}{B} (1_{\pG} \sm 1_{N}),
\end{align*}
\begin{align*}
\com_{C}(\beta(m)) & = \com_{C}(\alp(\gamma_{\theta}(m_{[0]}))(m_{[-1]} \sm 1_{N})) \\
& = (m_{[-1]} \sm 1_{N}) \ltak{}{C} \beta(m_{[0]}) \\
& = t_{C}(\beta(m_{[0]}))(m_{[-1]} \sm 1_{N}) \ltak{}{C} (1_{\pG} \sm 1_{N}) \\
& = \alp(\gamma_{\theta}(m_{[0]}))(m_{[-1]} \sm 1_{N}) \ltak{}{C} (1_{\pG} \sm 1_{N}) \\
& = \beta(m) \ltak{}{C} (1_{\pG} \sm 1_{N}),
\end{align*}

Hence,

\begin{enumerate}[label=\textup{(\roman*)}]

\item Fix $m,m',n,n' \in N$ and $a \in A$, it follows
\begin{align*}
\com_{B}(\alp(m)\beta(n)a\alp(m')\beta(n')) & = \com_{B}(\alp(m))\com_{B}(\beta(n))\com_{B}(a)\com_{B}(\alp(m'))\com_{B}(\beta(n')) \\
& = (\beta(n) \rtak{}{B} \alp(m))\com_{B}(a)(\beta(n') \rtak{}{B} \alp(m')).
\end{align*}
Similarly,
\[
\com_{C}(\alp(m)\beta(n)a\alp(m')\beta(n')) = (\beta(n) \ltak{}{C} \alp(m))\com_{C}(a)(\beta(n') \ltak{}{C} \alp(m')).
\]

\item Because $A$ is unital and $\com_{\G}$ is coassociative, it is straightforward to show
\[
(\com_{B} \odo \id)\com_{B} = (\id \odo \com_{B})\com_{B}, \qquad (\com_{C} \odo \id)\com_{C} = (\id \odo \com_{C})\com_{C}.
\]
We can also show that
\begin{align*}
(\com_{C} \odo \id)(\com_{B}(h \sm m)) & = (\id \odo \com_{B})(\com_{C}(h \sm m)), \\
(\com_{B} \odo \id)(\com_{C}(h \sm m)) & = (\id \odo \com_{C})(\com_{B}(h \sm m)),
\end{align*}
for all $h \sm m \in A$.

\item Considering the conjugate linear map
\[
\begin{array}{lccc}
(* \overline{\times} *): & A \ltak{}{C} A & \to & A \rtak{}{B} A \\
& a \ltak{}{C} a' & \mapsto & a^{*} \rtak{}{B} a'^{*}
\end{array},
\]
we have
\begin{align*}
\com_{C}((h \sm m)^{*})^{(* \overline{\times} *)} & = \com_{C}((h_{(1)})^{*} \sm (m \lhd_{\dtheta} S^{-1}_{\G}(h_{(2)}))^{*})^{(* \overline{\times} *)} \\
& = ((h_{(1)})^{*} \sm 1_{N})^{*} \rtak{}{B} ((h_{(2)})^{*} \sm (m \lhd_{\dtheta} S^{-1}_{\G}(h_{(3)}))^{*})^{*} \\
& = ((h_{(1)})^{*} \sm 1_{N})^{*} \rtak{}{B} (h_{(2)} \sm ((m \lhd_{\dtheta} S^{-1}_{\G}(h_{(4)})) \lhd_{\dtheta} h_{(3)}))  \\
& = (h_{(1)} \sm 1_{N}) \rtak{}{B} (h_{(2)} \sm  m) \\
& = \com_{B}(m \sm h).
\end{align*}
\end{enumerate}

\noindent\textbf{The antipode:} Fix $h,g \in \pG$ and $m,n \in N$. By direct computation, we have
\begin{align*}
S((h \sm m)(g \sm n)) & = S((hg_{(1)} \sm 1_{N})\alp((m \lhd_{\dtheta} g_{(2)})n)) \\
& = \beta(\hat\gamma_{\theta}(m \lhd_{\dtheta} g_{(2)})\hat\gamma_{\theta}(n))(S_{\G}(hg_{(1)}) \sm 1_{N}) \\
& = \beta(n)\beta(\hat\gamma_{\theta}(m \lhd_{\dtheta} g_{(2)}))(S_{\G}(g_{(1)}) \sm 1_{N})(S_{\G}(h) \sm 1_{N}) \\
& = \beta(\hat\gamma_{\theta}(n))(S_{\G}(g) \sm 1_{N})\beta(\hat\gamma_{\theta}(m))(S_{\G}(h) \sm 1_{N}) & \text{by } \eqref{eq:beta_hatgamma_antipode} \\
& = S(g \sm n)S(h \sm m).
\end{align*}

\noindent\textbf{The left and right counits:} By Proposition~\ref{prop:multiplier_algebra} and Lemma~\ref{lem:technical_beta}, the following formulas
\begin{enumerate}
\item $\cou_{B}((h \sm 1_{N})\alp(m)) = \cou_{\G}(h)\alp(m)$,
\item $\cou_{B}(\alp(m)(h \sm 1_{N})) = \alp(m \lhd_{\dtheta} h)$,
\item $\cou_{B}((h \sm 1_{N})\beta(m)) = \cou_{\G}(h)\alp(m)$,
\item $\cou_{B}(\beta(m)(h \sm 1_{N})) = \alp(m \lhd_{\dtheta} h)$
\end{enumerate}
are equivalent, for all $h \in \pG$ and $m \in N$. Similarly, the formulas
\begin{enumerate}
\item ${}_{C}\cou((h \sm 1_{N})\alp(m)) = \beta(\hat\gamma_{\theta}(m) \lhd_{\dtheta} S_{\G}(h))$,
\item ${}_{C}\cou(\alp(m)(h \sm 1_{N})) = \cou_{\G}(h)\beta(\hat\gamma_{\theta}(m))$,
\item ${}_{C}\cou((h \sm 1_{N})\beta(m)) = \beta(m \lhd_{\dtheta} S_{\G}(h))$,
\item ${}_{C}\cou(\beta(m)(h \sm 1_{N})) = \cou_{\G}(h)\beta(m)$,
\end{enumerate}
are equivalent, for all $h \in \pG$ and $m \in N$. Hence, for any $h \in \pG$ and $m,n \in N$, we have
\begin{align*}
\cou_{B}(\alp(n) {}\indices*{^{B}}\rhd (h \sm 1_{N})\beta(m)) & = \cou_{B}((h \sm 1_{N})\beta(nm)) \\
& = \cou_{\G}(h)\alp(nm) \\
& = \alp(n)\cou_{B}((h \sm 1_{N})\beta(m))
\end{align*}
and
\begin{align*}
\cou_{B}((h \sm 1_{N})\alp(m) \lhd_{B} \alp(n)) & = \cou_{B}((h \sm 1_{N})\alp(mn)) \\
& = \cou_{\G}(h)\alp(mn) \\
& = \cou_{B}((h \sm 1_{N})\alp(m))\alp(n).
\end{align*}
Similarly,
\begin{align*}
{}_{C}\cou(\beta(n) {}_{C}\rhd \beta(m)(h \sm 1_{N})) & = {}_{C}\cou(\beta(mn)(h \sm 1_{N})) \\
& = \cou_{\G}(h)\beta(mn) \\
& = \beta(n){}_{C}\cou(\beta(m)(h \sm 1_{N}))
\end{align*}
and
\begin{align*}
{}_{C}\cou(\alp(m)(h \sm 1_{N}) \lhd^{C} \beta(n)) & = {}_{C}\cou(\alp(\gamma_{\theta}(n)m)(h \sm 1_{N})) \\
& = \cou_{\G}(h)\beta(n\hat\gamma_{\theta}(m)) \\
& = {}_{C}\cou(\alp(m)(h \sm 1_{N}))\beta(n)
\end{align*}
for all $h \in \pG$ and $m,n \in N$. Those equalities imply that $\cou_{B} \in \textrm{Hom}({}^{B}A_{B},{}_{B}B_{B})$ and ${}_{C}\cou \in \textrm{Hom}({}_{C}A^{C},{}_{C}C_{C})$. Moreover,

\begin{enumerate}[label=\textup{(\roman*)}]
\item For $h \in \pG$ and $m \in N$, it holds $S(\beta(m)) = \beta(\hat\gamma_{\theta}(m_{[0]}))(S_{\G}(m_{[-1]}) \sm 1_{N}) = \alp(m) = t^{-1}_{B}(\beta(m))$ and $S(\alp(m)) = \beta(\hat\gamma_{\theta}(m)) = t^{-1}_{C}(\alp(m))$. Hence
\begin{align*}
S(t_{B}(\alp(m))t_{C}(\beta(n))at_{B}(\alp(m'))t_{C}(\beta(n'))) & = \alp(m')\beta(n')S(a)\alp(m)\beta(n),
\end{align*}
for all $m,m',n,n' \in N$ and $a \in A$.

\item It follows
\begin{align*}
S(S(h \sm m)^{*})^{*}
& = S(S(h \sm 1_{N})^{*})^{*}S(S(\alp(m))^{*})^{*} \\
& = (h \sm 1_{N})S(\beta(\hat\gamma_{\theta}(m))^{*})^{*} \\
& = (h \sm 1_{N})\alp(m^{*})^{*} \\
& = h \sm m,
\end{align*}
for all $h \in \pG$ and $m \in N$.

\item Fix $h \sm m, g \sm n \in A$. It hold
\begin{align*}
(\cou_{B} \odot \id)((1 \odo (g \sm n))\com_{B}(h \sm m)) & = (\cou_{B} \odot \id)((h_{(1)} \sm 1_{N}) \odo (g \sm n)(h_{(2)} \sm m)) \\
& = (g \sm n)(h_{(2)} \sm m)t_{B}(\cou_{B}(h_{(1)} \sm 1_{N})) \\
& = (g \sm n)(h \sm m)
\\\\
(\id \odot \cou_{B})(((g \sm n) \odo 1)\com_{B}(h \sm m)) & = (\id \odot \cou_{B})((g \sm n)(h_{(1)} \sm 1_{N}) \odo (h_{(2)} \odo m)) \\
& = (g \sm n)(h_{(1)} \sm 1_{N})\cou_{B}(h_{(2)} \sm m) \\
& = (g \sm n)(h_{(1)} \sm 1_{N})\cou_{\G}(h_{(2)})\alp(m) \\
& = (g \sm n)(h \sm m)
\end{align*}
and
\begin{align*}
({}_{C}\cou \odot \id)(\com_{C}(h \sm m)(1 \odo (g \sm n))) & = ({}_{C}\cou \odot \id)((h_{(1)} \sm 1_{N}) \odo (h_{(2)} \sm m)(g \sm n)) \\
& = \cou_{C}(h_{(1)} \sm 1_{N})(h_{(2)} \sm m)(g \sm n) \\
& = (h \sm m)(g \sm n),
\\\\
(\id \odot {}_{C}\cou)(\com_{C}(h \sm m)((g \sm n) \odo 1)) & = (\id \odot {}_{C}\cou)((h_{(1)} \sm 1_{N})(g \sm n) \odo (h_{(2)} \sm m)) \\
& = t_{C}({}_{C}\cou(h_{(2)} \sm m))(h_{(1)} \sm 1_{N})(g \sm n) \\
& = \alp(m \lhd_{\dtheta} S^{-1}_{\G}(h_{(2)}))(h_{(1)} \sm 1_{N})(g \sm n) \\
& = (h \sm m)(g \sm n).
\end{align*}

\item Fix $h \sm m, g \sm n \in A$. It holds
\begin{align*}
\m_{A}(S \odo \id)(\com_{C}(h \sm m)(1 \odo (g \sm n))) & = \m_{A}(S \odo \id)((h_{(1)} \sm 1_{N}) \odo (h_{(2)} \sm m)(g \sm n)) \\
& = (S_{\G}(h_{(1)}) \sm 1_{N})(h_{(2)} \sm m)(g \sm n) \\
& = \cou_{\G}(h)\alp(m)(g \sm n) \\
& = \cou_{B}(h \sm m)(g \sm n)
\intertext{and}
\m_{A}(\id \odo S)(((h \sm m) \odo 1)\com_{B}(g \sm n)) & = \m_{A}(\id \odo S)((h \sm m)(g_{(1)} \sm 1_{N}) \odo (g_{(2)} \sm n)) \\
& = (h \sm m)(g_{(1)} \sm 1_{N})\beta(\hat\gamma_{\theta}(n))(S_{\G}(g_{(2)}) \sm 1_{N}) \\
& = (h \sm m)(g_{(1)}S_{\G}(g_{(2)}) \sm 1_{N})\beta(\hat\gamma_{\theta}(n) \lhd_{\dtheta} S_{\G}(g_{(3)})) & \text{by } \eqref{eq:YD_condition_smash_1} \\
& = (h \sm m)\beta(\hat\gamma_{\theta}(n) \lhd_{\dtheta} S_{\G}(g)) \\
& = (h \sm m){}_{C}\cou(g \sm n).
\end{align*}
\end{enumerate}

\item By Proposition~\ref{prop:multiplier_algebra}, Lemma~\ref{lem:technical_beta} and left and right invariance of $\varphi_{\G}$ (see Remark~\ref{rem:left_right_invariance}), the formulas
\begin{enumerate}
\item ${}_{B}\psi_{B}((h \sm 1_{N})\alp(m)) = \varphi_{\G}(h)\alp(m)$,
\item ${}_{B}\psi_{B}(\alp(m)(h \sm 1_{N})) = \varphi_{\G}(h)\alp(m)$,
\item ${}_{B}\psi_{B}((h \sm 1_{N})\beta(m)) = \varphi(hm_{[-1]})\alp(m_{[0]})$,
\item ${}_{B}\psi_{B}(\beta(m)(h \sm 1_{N}))= \varphi_{\G}(h_{(2)}m_{[-1]})\alp(m_{[0]} \rhd_{\dtheta} h_{(1)})$,
\end{enumerate}
are equivalent for all $m \in N$ and $h \in \pG$. Similar, the formulas
\begin{enumerate}
\item ${}_{C}\phi_{C}((h \sm 1_{N})\alp(m)) = \varphi_{\G}(hS_{\G}(m_{[-1]}))\beta(m_{[0]})$,
\item ${}_{C}\phi_{C}(\alp(m)(h \sm 1_{N})) = \varphi_{\G}(S_{\G}(m_{[-1]})h)\beta(\hat\gamma_{\theta}(m_{[0]}))$,
\item ${}_{C}\phi_{C}((h \sm 1_{N})\beta(m)) = \varphi_{\G}(h)\beta(m)$,
\item ${}_{C}\phi_{C}(\beta(m)(h \sm 1_{N})) = \varphi_{\G}(h)\beta(m)$,
\end{enumerate}
are equivalent for all $m \in N$ and $h \in \pG$. Then,
\begin{align*}
{}_{B}\psi_{B}(\alp(n) {}_{B}\rhd \alp(m)(h \sm 1_{N})) & = \alp(n){}_{B}\psi_{B}(\alp(m)(h \sm 1_{N})), \\
{}_{B}\psi_{B}((h \sm 1_{N})\alp(m) \lhd_{B} \alp(n')) & = {}_{B}\psi_{B}((h \sm 1_{N})\beta(m))\alp(n')
\end{align*}
and
\begin{align*}
{}_{C}\phi_{C}(\beta(n) {}_{C}\rhd \beta(m)(h \sm 1_{N})) & = \beta(n){}_{C}\phi_{C}(\beta(m)(h \sm 1_{N})), \\
{}_{C}\phi_{C}((h \sm 1_{N})\beta(m) \lhd_{C} \beta(n')) & = {}_{C}\phi_{C}((h \sm 1_{N})\beta(m))\beta(n')
\end{align*}
for all $m,n \in N$ and $h \in \pG$. The equalities above shows that ${}_{B}\psi_{B} \in \textrm{Hom}({}_{B}A_{B},{}_{B}B_{B})$ and ${}_{C}\phi_{C} \in \textrm{Hom}({}_{C}A_{C},{}_{C}C_{C})$. On the other hand, since the Yetter--Drinfeld integral $\mu$ is a non-zero positive faithful functional, then the maps $\mu_{B}$ and $\mu_{C}$ are too. Moreover, we have

\begin{enumerate}[label=\textup{(\roman*)}]
\item For all for all $m,n \in N$ and $h,g \in \pG$, it holds
\small
\begin{align*}
({}_{B}\psi_{B} \odot \id)((1 \odo (h \sm m))\com_{B}((g \sm 1_{N})\alp(n)))
& = t_{B}({}_{B}\psi_{B}(g_{(1)} \sm 1_{N}))(h \sm m)(g_{(2)} \sm 1_{N})\alp(n) \\
& = \varphi_{\G}(g_{(2)})(h \sm m)(g_{(1)} \sm 1_{N})\alp(n) \\
& = (h \sm m){}_{B}\psi_{B}((g \sm 1_{N})\alp(n))
\end{align*}
\normalsize
and
\small
\begin{align*}
(\id \odot {}_{C}\phi_{C})(\com_{C}((h \sm 1_{N})\beta(m))((g \sm n) \odo 1))
& = t_{C}({}_{C}\phi_{C}(h_{(2)} \sm 1_{N}))(h_{(1)} \sm 1_{N})\beta(m)(g \sm n) \\
& = \varphi_{\G}(h_{(2)})(h_{(1)} \sm 1_{N})\beta(m)(g \sm n) \\
& = {}_{C}\phi_{C}((h \sm 1_{N})\beta(m))(g \sm n).
\end{align*}
\normalsize

\item Easily $\mu_{B}(t_{C}(\beta(m))) = \mu_{B}(\alp(\gamma_{\theta}(m))) = \mu(m)$ and $\mu_{C}(t_{B}(\alp(m))) = \mu_{C}(\beta(m)) = \mu(m)$ for all $m \in N$. Moreover, $\mu_{C}({}_{C}\cou(h \sm m)) = \mu(\hat\gamma_{\theta}(m) \lhd_{\dtheta} S_{\G}(h)) = \cou_{\G}(h)\mu(m) = \mu_{B}(\cou_{B}(h \sm m))$ for all $m \in N$ and $h \in \pG$.

\item Since $\psi(h \sm m) = \varphi_{\G}(h)\mu(m)$ and $\phi(h \sm m) = \mu_{C}({}_{C}\phi_{C}(h \sm m)) = \varphi_{\G}(hS_{\G}(m_{[-1]}))\mu(m_{[0]}) = \varphi_{\G}(h)\mu(m)$ for all $m \in N$ and $h \in \pG$, then the maps $\psi$ and $\phi$ are positive faithful functional on $A$. Moreover, in this case we have $\phi = \psi$, and easily we can verify that
\begin{align*}
\psi(\beta(m)(h \sm 1_{N})) & = \phi(\beta(m)(h \sm 1_{N})) = \varphi_{\G}(h)\mu(m),
\end{align*}
and
\[
\phi((h \sm 1_{N})\beta(m)) = \psi((h \sm 1_{N})\beta(m)) = \varphi_{\G}(h)\mu(m)
\]
for any $h \in \pG$ and $m \in N$.
\end{enumerate}
\end{enumerate}
\fin

\begin{prop}
Let $\G$ be an algebraic quantum group of compact type and $(N,\theta,\dtheta,\mu)$ be an unital braided commutative measured Yetter--Drinfeld $\Gc$-\as-algebra. The following statements are equivalent
\begin{enumerate}[label=\textup{(\roman*)}]
\item $\G$ is an algebraic quantum group of Kac type and $\mu$ is a tracial Yetter--Drinfeld integral;
\item $\mathcal{A}(N,\theta,\theta,\mu)$ is a measured multiplier Hopf \as-algebroid of Kac type in sense of Definition~\ref{def:MHAd_Kac}.
\end{enumerate}
\end{prop}
\pr
Denote by $S$ the antipode of the measured multiplier Hopf \as-algebroid $\mathcal{A}(N,\theta,\dtheta,\mu)$. Observe that the square of the antipode $S$ is given by
\[
S^{2}(h \sm m) = S(\beta(\hat\gamma_{\theta}(m))(S_{\G}(h) \sm 1_{N})) = (S^{2}_{\G}(h) \sm 1_{N})\alp(\hat\gamma_{\theta}(m)) = S^{2}_{\G}(h) \sm \hat\gamma_{\theta}(m)
\]
for all $h \sm m \in A = \pG \sm_{\dtheta} N$. Suppose that $\G$ is of Kac type and $\mu$ is a tracial Yetter--Drinfeld integral. Since $S^{2}_{\G} = \id_{\G}$ and $\gamma_{\theta} = \sigma^{\mu} = \id_{N}$, we have directly $S^{2} = \id_{A}$.

Reciprocally, suppose now that $S^{2} = \id_{A}$. By linearity of $\varphi_{\G}$ and $\mu$, we have
\[
0 = (\id \odo \mu)((h \sm m)-(S^{2}_{\G}(h) \sm \gamma_{\theta}(m))) = (S^{2}_{\G}(h)-h)\mu(m)
\]
and
\[
0 = (\varphi_{\G} \odo \id)((h \sm m)-(S^{2}_{\G}(h) \sm \gamma_{\theta}(m))) = \varphi_{\G}(h)(m-\gamma_{\theta}(m))
\]
for all $h \in \pG$ and $m \in N$. Since $\varphi_{\G}$ and $\mu$ are non-zero functionals, then necessarily $S^{2}_{\G}(h) = h$ and $\gamma_{\theta}(m) = m$ for all $h \in \pG$ and $m \in N$.
\fin

The next proposition gives the modular automorphisms of the measured multiplier Hopf \as-algebroid $\mathcal{A}(N,\theta,\dtheta,\mu)$. Those modular automorphisms will come in handy later when we have to deal with the Pontrjagin dual of $\mathcal{A}(N,\theta,\dtheta,\mu)$.

\begin{prop}
The modular automorphisms associated with the measured multiplier Hopf \as-algebra $\mathcal{A}(N,\theta,\dtheta,\mu)$ are given by
\begin{enumerate}[label=\textup{(\roman*)}]
\item $\sigma_{B}(\alp(m)) = (S^{-1}_{B} \circ S^{-1}_{C})(\alp(m)) = \alp(\gamma_{\theta}(m))$ and $\sigma_{C}(\beta(m)) = (S_{B}\circ S_{C})(\beta(m)) = \beta(\hat\gamma_{\theta}(m))$ for all $m \in N$; and
\item $\sigma^{\phi}(h \sm m) = \sigma_{\G}(h) \sm \gamma_{\theta}(m) = \sigma^{\psi}(h \sm m)$ for all $h \in \pG$ and $m \in N$.
\end{enumerate}
\end{prop}
\pr
\begin{enumerate}[label=\textup{(\roman*)}]
\item By direct computation
\begin{align*}
\mu_{B}(\alp(m)\alp(n)) & = \mu_{B}(\alp(mn)) = \mu(mn) \\
& = \mu(n\gamma_{\theta}(m)) = \mu_{B}(\alp(n)\alp(\gamma_{\theta}(m))) \\
\intertext{and}
\mu_{C}(\beta(m)\beta(n)) & = \mu_{C}(\beta(nm)) = \mu(nm) \\
& = \mu(\hat\gamma_{\theta}(m)n) = \mu_{C}(\beta(n)\beta(\hat\gamma_{\theta}(m))),
\end{align*}
for all $m,n \in N$. On the other hand, because $(S^{-1}_{B}\circ S^{-1}_{C})(\alp(m)) = (t_{C} \circ t_{B})(\alp(m)) = \alp(\gamma_{\theta}(m))$ and $(S_{B}\circ S_{C})(\beta(m)) = (t^{-1}_{C} \circ t^{-1}_{B})(\beta(m)) = \beta(\hat\gamma_{\theta}(m))$ for all $m \in N$, thus $\sigma_{B}(\alp(m)) = (S^{-1}_{B}\circ S^{-1}_{C})(\alp(m)) = \alp(\gamma_{\theta}(m))$ and $\sigma_{C}(\beta(m)) = (S_{B}\circ S_{C})(\beta(m)) = \beta(\hat\gamma_{\theta}(m))$, for all $m \in N$.

\item Fix $h \in \pG$ and $m \in N$. We have
\begin{align*}
\phi((g \sm n)(\sigma_{\G}(h) \sm 1_{N})) & = \phi(g\sigma_{\G}(h)_{(1)} \sm (n \lhd_{\dtheta} \sigma_{\G}(h)_{(2)}) \\
& = \varphi_{\G}(g\sigma_{\G}(h)_{(1)})\mu(n \lhd_{\dtheta} \sigma_{\G}(h)_{(2)}) \\
& = \varphi_{\G}(g\sigma_{\G}(h)_{(1)})\cou_{\G}(\sigma_{\G}(h)_{(2)})\mu(n) \\
& = \varphi_{\G}(g\sigma_{\G}(h))\mu(n) \\
& = \phi((h \sm 1_{N})(g \sm n)) \\
\intertext{and}
\phi((g \sm n)(1_{\pG} \sm \gamma_{\theta}(m))) & = \varphi_{\G}(g)\mu(n\gamma_{\theta}(m)) \\
& = \varphi_{\G}(g)\mu(mn) \\
& = \varphi_{\G}(g_{(1)})\mu((m \lhd_{\dtheta} g_{(2)})n) \\
& = \phi(g_{(1)} \sm (m \lhd_{\dtheta} g_{(2)})n) \\
& = \phi((1_{\pG} \sm m)(g \sm n))
\end{align*}
for all $g \in \pG$ and $n \in N$. The last equalities together imply that
\begin{align*}
\phi((g \sm n)(\sigma_{\G}(h) \sm \gamma_{\theta}(m)))
& = \phi((h \sm m)(n \sm g))
\end{align*}
for all $h,g \in \pG$ and $m,n \in N$. By unicity of the automorphism $\sigma^{\phi}$, we have $\sigma^{\phi}(h \sm m) = \sigma_{\G}(h) \sm \gamma_{\theta}(m)$ for all $h \in \pG$, $m \in N$. Observe that in particular, we obtain $\sigma^{\phi}(\alp(m)) = \alp(\gamma_{\theta}(m)) = \sigma_{B}(\alp(m))$ and $\sigma^{\phi}(\beta(m)) = (S^{-1}\circ(\sigma^{\phi})^{-1}\circ S)(\beta(m)) = \beta(\hat\gamma_{\theta}(m)) = \sigma_{C}(\beta(m))$ for all $m \in N$. The last equality follows from the fact that $\phi\circ S = \phi = \phi \circ S^{-1}$.
\end{enumerate}
\fin

\begin{prop}
Let $\G$ be an algebraic quantum group of compact type and $(N,\theta,\dtheta,\mu)$ be a unital braided commutative measured Yetter--Drinfeld $\Gc$-\as-algebra, with canonical automorphisms denoted by $\gamma_{\theta}$ and $\hat\gamma_{\theta}$. Consider the measured multiplier Hopf \as-algebroid $\mathcal{A}(N,\theta,\dtheta,\mu)$ construct by Theorem~\ref{theo:MHAd_YD}. The following statements are equivalent
\begin{enumerate}[label=\textup{(\roman*)}]
\item $* \circ S_{C} \circ {}_{C}\cou = \cou_{B} \circ *$;
\item $* \circ S_{B} \circ \cou_{B} = {}_{C}\cou \circ *$;
\item $\mu$ is a tracial Yetter--Drinfeld integral, i.e. $\gamma_{\theta} = \hat\gamma_{\theta} = \id_{N}$.
\end{enumerate}
\end{prop}
\pr
Since $S^{-1}_{C} \circ * = * \circ S_{B}$, it is evident that the statements (i) and (ii) are equivalent. Then, it is enough to show that the statements (ii) and (iii) are equivalent. For that, first observe that
\begin{align*}
S_{B}(\cou_{B}(h \sm m))^{*} & = t^{-1}_{C}(\cou_{\G}(h)\alp(m))^{*} \\
& = \cou_{\G}(h^{*})\beta(\hat\gamma_{\theta}(m))^{*} \\
& = \cou_{\G}(h^{*})\beta(m^{*}) \\
\intertext{and}
{}_{C}\cou((h \sm m)^{*}) & = {}_{C}\cou(h^{*}_{(1)} \sm (m^{*} \lhd_{\dtheta} h^{*}_{(2)})) \\
& = \beta(\hat\gamma_{\theta}(m^{*} \lhd_{\dtheta} h^{*}_{(2)}) \lhd_{\dtheta} S_{\G}(h^{*}_{(1)})) \\
& = \beta(\hat\gamma_{\theta}(m^{*}) \lhd_{\dtheta} S^{2}_{\G}(h^{*}_{(2)})S_{\G}(h^{*}_{(1)})) \\
& = \cou_{\G}(h^{*})\beta(\hat\gamma_{\theta}(m^{*}))
\end{align*}
for all $h \in \pG$ and $m \in N$. If the statement (ii) is true, follows from the equalities above that $\beta(\hat\gamma_{\theta}(m)) = \beta(m)$ for all $m \in H$. Thus, from the injectivity of $\beta$ we have $\hat\gamma_{\theta} = \id_{N}$, and then also $\gamma_{\theta} = \hat\gamma^{-1}_{\theta} = \id_{N}$. Reciprocally, if the statement (iii) is true, thus the statement (i) follows directly.
\fin

\begin{rema}
Since there are examples of unital braided commutative Yetter--Drinfeld \as-algebras which admits non-trivial canonical automorphisms (see Example~\ref{ex:bc_yd_quotient_type}), thus the proposition above shows that \cite[Proposition~7.3]{TVD18} is not true in general.
\end{rema}

\begin{exam}[Finite transformation groupoids]\label{ex:MHAd_finite_classic}
Let $G$ be a finite group acting by the left on a finite space $X$ and $\nu:X \to \R^{+}_{0}$ be a non-zero $G$-invariant function. Consider the unital braided commutative measured Yetter--Drinfeld $\Gc$-\as-algebra $(K(X),\theta,\hat\theta,\mu_{\nu})$ arising from the action of $G$ on $X$ (see Example~\ref{ex:bc_yd_finite_transformation_groupoid}). The measured multiplier Hopf \as-algebroid $\mathcal{A}(K(X),\theta,\dtheta,\mu_{\nu})$ is given by
\[
\begin{array}{cccccc}
A = & K(G) \odo_{\dtheta} K(X) & \cong  & K(G \ltimes X) \\
& p \odo f & \mapsto & \displaystyle\sum_{g \in G,x \in X} p(g)f(x)\delta_{(g,x)}
\end{array},
\]
\[
\begin{array}{rccc}
\alpha_{\dtheta}: & K(X) & \to & \M(K(G)\,\odo_{\dtheta}\,K(X)) \\
& f & \mapsto & \displaystyle\sum_{g \in G,x \in X} f(d(g,x))\delta_{(g,x)}
\end{array},
\begin{array}{rccc}
\beta_{\theta}: & K(X) & \to & \M(K(G)\,\odo_{\dtheta}\,K(X)) \\
& f & \mapsto & \displaystyle\sum_{g \in G, x \in X} f(r(g,x))\delta_{(g,x)}
\end{array}
\]
\[
B := \{\alpha_{\dtheta}(f) = d^{\bullet}(f) : f \in K(X)\} \cong K((G \ltimes X)^{(0)}) \cong \{\beta_{\theta}(f) = r^{\bullet}(f) : f \in K(X)\} =: C,
\]
\[
t_{B}: d^{\bullet}(f) \mapsto r^{\bullet}(f), \quad t_{C}=t^{-1}_{B}: r^{\bullet}(f) \mapsto d^{\bullet}(f).
\]
Given $p \odo f \in K(G \ltimes X)$,
\begin{enumerate}
\item Since $A \rtak{}{B} A = A \ltak{}{C} A = K((G \ltimes X)^{(2)})$, we have
\begin{align*}
\com_{B}(p \odo f)((g,x),(g',x')) = \com_{C}(p \odo f)((g,x),(g',x')) & = (p_{(1)} \odo 1_{X})(g,x)(p_{(2)} \odo f)(g',x') \\
& = p_{(1)}(g)p_{(2)}(g')f(x') \\
& = p(gg')f(x') \\
& = (p \odo f)((g,x)(g',x'))
\end{align*}
for all $((g,x),(g',x')) \in (G \ltimes X)^{(2)}$.

\item It hold
\begin{align*}
S(p \odo f)(g,x) & = \beta_{\theta}(f)(S_{\G}(p) \odo 1_{X})(g,x) \\
& = r^{\bullet}(f)(g,x)p(g^{-1}) \\
& = (p \odo f)(g^{-1},g\cdot x) \\
& = (p \odo f)((g,x)^{-1}),
\end{align*}
for all $(g,x) \in G \ltimes X$, and
\begin{align*}
\cou_{B}(\alp(f)(p \odo 1_{X}))(e,x) & = \cou_{\G}(p)f(x) = p(e)f(x) = \cou_{\G}(p)f(e\cdot x) = \cou_{C}((p \odo 1_{X})\beta(f))(e,x) \\\\
\mu_{B}(\alp(f)) & = \mu_{\nu}(f) = \sum_{x \in X}\nu(x)f(x) = \mu_{\nu}(f) = \mu_{C}(\beta(f)) \\\\
{}_{B}\psi_{B}(\alp(f)(p \odo 1_{X}))(e,x) & = \sum_{g \in G}p(g)f(x) = \sum_{g \in G}p(g)f(e \cdot x) = {}_{C}\phi_{C}((p \odo 1_{X})\beta(f))(e,x)
\end{align*}
for all $(e,x) \in (G \ltimes X)^{(0)}$.
\end{enumerate}
\end{exam}

\begin{exam}[Fell bundles]\label{ex:MHAd_fell}
Let $\fell{A}$ be a separable Fell bundle over a discrete group $G$ and $\nu: \fell{A}_{e} \to \C$ be a faithful state on the C*-algebra $\fell{A}_{e}$. Assume that $\fell{A}$ is equipped with a family of Fell bundle isomorphisms $\rho = \{\rho_{g}:\fell{A} \to \fell{A}^{\ad(g^{-1})}\}_{g\in G}$ verifying the conditions
\begin{enumerate}
\item $\rho_{e} = \id_{\fell{A}}$;
\item $\rho_{g}\circ\rho_{g'} = \rho_{g'g}$ for all $g,g' \in G$;
\item for any $g,h \in G$, we have $ab = \rho_{g^{-1}}(b)a$ for all $a \in \fell{A}_{g}$ and $b \in \fell{A}_{h}$.
\end{enumerate}
Using the unital braided commutative measured Yetter--Drinfeld $\Gc$-\as-algebra $(\Gamma_{c}(G,\fell{A}),\theta_{\fell{A}},\hat\theta_{\fell{A},\rho},\mu_{\nu})$ (see Example~\ref{ex:bc_yd_fell_bundle}), we obtain the measured multiplier Hopf \as-algebroid $\mathcal{A}(\Gamma_{c}(G,\fell{A}),\theta_{\fell{A}},\hat\theta_{\fell{A},\rho},\mu_{\nu})$ with total algebra given by
\[
A = \C[G] \#_{\dtheta_{\fell{A},\rho}} \Gamma_{c}(G,\fell{A}).
\]
\end{exam}

\begin{exam}[Algebraic quantum groups of compact type]\label{ex:MHAd_aqg}
In general, any algebraic quantum group is a particular example of a measured multiplier Hopf \as-algebroid (see \cite{T17}). Indeed, for a given algebraic quantum group $\G=(\pG,\com_{\G},\varphi_{\G})$, the collection
\[
\mathcal{A}(\G) := (\pG,\C,\C,\id_{\C},\id_{\C},\com_{\G},\com_{\G},\id_{\C},\id_{\C},\psi_{\G},\varphi_{\G})
\]
yields a measured multiplier Hopf \as-algebroid, called {\em the canonical measured multiplier Hopf \as-algebroid associated with $\G$}. In this example, we show that algebraic quantum groups can also be obtained as examples of measured multiplier Hopf \as-algebroids associated with braided commutative Yetter--Drinfeld \as-algebras.

Suppose now that $\G$ is an algebraic quantum group of compact type. Consider the trivial unital braided commutative measured Yetter--Drinfeld $\Gc$-\as-algebra $(\C,\tr_{\Gc},\tr_{\dGo},\id_{\C})$. Thus, the measured multiplier Hopf \as-algebroid $\mathcal{A}(\C,\tr_{\Gc},\tr_{\dGo},\id_{\C})$ is given by
\[
A = \pG \sm_{\tr} \C = \pG, \quad \alp = \beta: 1_{\C} \mapsto 1_{\pG}, \\
\]
\[
\quad B = \alp(\C) = \C\cdot1_{\pG} = \beta(\C) = C,
\]
\[
t_{B} = t_{C} = \id_{\C}, \quad \com_{B} = \com_{C} = \com_{\G},
\]
\[
S = S_{\G}, \quad \cou_{B} = \cou_{C} = \cou_{\G},
\]
\[
\mu_{B} = \mu_{C} = \id_{\C} , \quad {}_{B}\psi_{B} = {}_{C}\phi_{C} = \varphi_{\G}, \quad \psi = \phi = \varphi_{\G}.
\]
In other words, we have
\[
\mathcal{A}(\C,\tr_{\Gc},\tr_{\dGo},\id_{\C}) = \mathcal{A}(\G).
\]
\end{exam}

\begin{exam}[Quotient type left coideals]\label{ex:MHAd_quotient}
Let $\G$ be an algebraic quantum group of compact type and $\HH$ be a quantum subgroup defined by the surjective \as-morphism $\pi:\pG \to \pol(\HH)$. From Example~\ref{ex:bc_yd_quotient_type}, the triplet $(\pol(\HH\backslash\G),\theta=(S^{-1}_{\G} \odo \id)\circ\Sigma\circ\com_{\G}|_{\pol(\HH\backslash\G)},\dtheta=\ad_{\Sigma(U^{*})}|_{\pol(\HH\backslash\G)})$ gives an example of a unital braided commutative Yetter--Drinfeld $\Gc$-\as-algebra with canonical automorphism given by $\gamma_{\theta} = S^{-2}_{\G}|_{\pol(\HH\backslash\G)}$, $\hat\gamma_{\theta} = S^{2}_{\G}|_{\pol(\HH\backslash\G)}$. If we consider the Yetter--Drinfeld integral $\mu := \cou_{\G}|_{\pol(\HH\backslash\G)}$, thus, we obtain the measured multiplier Hopf \as-algebroid
\[
\mathcal{A}(\pol(\HH\backslash\G),(S^{-1}_{\G} \odo \id)\circ\Sigma\circ\com_{\G}|_{\pol(\HH\backslash\G)},\ad_{\Sigma(U^{*})}|_{\pol(\HH\backslash\G)},\cou_{\G}|_{\pol(\HH\backslash\G)}),
\]
where
\[
A= \pG \sm_{\dtheta} \pol(\HH\backslash\G),
\]
\[
\alp: m \in \pol(\HH\backslash\G) \mapsto 1_{\G} \sm m \in A, \quad \beta: m \in \pol(\HH\backslash\G) \mapsto S^{-1}_{\G}(m_{(2)}) \sm m_{(1)} \in A,
\]
\[
B= 1_{\pG} \sm \pol(\HH\backslash\G), \quad C= \beta(\pol(\HH\backslash\G)),
\]
\[
t_{B}: 1 \sm m \mapsto S^{-1}_{\G}(m_{(2)}) \sm m_{(1)}, \quad t_{C}: S^{-1}_{\G}(m_{(2)}) \sm m_{(1)} \mapsto 1 \sm S^{-2}_{\G}(m),
\]
\[
S: h \sm m \mapsto S^{-1}_{\G}(m_{(2)})S_{\G}(h_{(3)}) \sm S^{2}_{\G}(h_{(2)})S^{2}(m_{(1)})S_{\G}(h_{(1)}),
\]
\[
\cou_{B}: \alp(m)(h \sm 1) \mapsto \alp(S_{\G}(h_{(1)})mh_{(2)}), \quad \cou_{C}: (h \sm 1)\beta(m) \mapsto \cou_{\G}(h)\beta(S^{2}_{\G}(h_{(2)})mS_{\G}(h_{(1)})),
\]
\[
\mu_{B}: 1 \sm m \mapsto \cou_{\G}(m), \qquad \mu_{C}: S^{-1}_{\G}(m_{(2)}) \sm m_{(1)} \mapsto \cou_{\G}(m)
\]
\[
{}_{B}\psi_{B}: \alp(m)(h \sm 1) \mapsto \varphi_{\G}(h)\alp(m), \qquad {}_{C}\phi_{C}: (h \sm 1)\beta(m) \mapsto \varphi_{\G}(h)\beta(m),
\]
\[
\psi = \phi: h \sm m \mapsto \varphi_{\G}(h)\cou_{\G}(m).
\]
\end{exam}

\begin{rema}
When we work with the smash product \as-algebra $\pG \sm_{\dtheta} N$, sometimes we will use the notation $\iota^{\dtheta}_{\G}$ and $\iota^{\dtheta}_{N}$ for the canonical faithful embeddings $h \in \pG \mapsto h \sm_{\dtheta} 1_{N} \in \pG \sm_{\dtheta} N$ and $m \in N \mapsto 1_{N} \sm_{\dtheta} n \in \pG \sm_{\dtheta} N$, respectively.
\end{rema}

\begin{prop}\label{prop:morphism_of_MHAd}
Let $\G$ be an algebraic quantum group of compact type, $(N_{i},\theta_{i},\dtheta_{i},\mu_{i})$ be a unital braided commutative measured Yetter--Drinfeld $\Gc$-\as-algebra for $i = 1,2$, and $f:N_{1} \to N_{2}$ be a morphism between the Yetter--Drinfeld $\Gc$-\as-algebras $(N_{1},\theta_{1},\dtheta_{1})$ and $(N_{2},\theta_{2},\dtheta_{2})$. Then, the linear map
\[
\begin{array}{lccc}
\pi_{f} := \id_{\pG} \sm f: & \pG \sm_{\dtheta_{1}} N_{1} & \to & \pG \sm_{\dtheta_{2}} N_{2} \\
& h \sm m_{1} & \mapsto & h \sm f(m_{1})
\end{array}
\]
yields a morphism between $\mathcal{A}(N_{1},\theta_{1},\dtheta_{1},\mu_{1})$ and $\mathcal{A}(N_{2},\theta_{2},\dtheta_{2},\mu_{2})$ that preserve elementwise the sub-\as-algebra $\pG$, i.e. such that $\pi_{f}\circ\iota^{\dtheta_{1}}_{\G} = \iota^{\dtheta_{2}}_{\G}$. Moreover, any morphism $\pi$ between measured multiplier Hopf \as-algebroids of the form $\mathcal{A}(N,\theta,\dtheta,\mu)$, with $(N,\theta,\dtheta,\mu)$ being a unital braided commutative measured Yetter--Drinfeld $\Gc$-\as-algebra, that preserve elementwise the sub-\as-algebra $\pG$ arise in this way. 
\end{prop}
\pr
Keep the notations $\alp_{i} = \alp_{\dtheta_{i}}$, $\beta_{i} = \beta_{\theta_{i}}$, $A_{i} = \pG \sm_{\dtheta_{i}} N_{i}$, $B_{i}=\alp_{i}(N)$ and $C_{i}=\beta_{i}(N)$. Since $f$ is a $\dGo$-equivariant \as-homomorphism, the linear map $\pi_{f}$ is well defined and we have $\pi_{f}\circ\alp_{\dtheta_{1}} = \alp_{\dtheta_{2}}\circ f$. By $\Gc$-equivariance of $f$, we also have $\pi_{f}\circ\beta_{\theta_{1}} = \beta_{\theta_{2}}\circ f$, thus the \as-homomorphism $\pi_{f}: A_{1} \to A_{2}$ satisfies the condition $\pi_{f}\circ t_{B_{1}} = t_{B_{2}}\circ \pi_{f}$ and $\pi_{f}\circ t_{C_{1}} = t_{C_{2}}\circ\pi_{f}$. Moreover, we claim that we also have $\com_{B_{2}}\circ\pi_{f} = (\pi_{f} \rtak{}{B_{1}} \pi_{f})\circ\com_{B_{1}}$ and $\com_{C_{2}}\circ\pi_{f} = (\pi_{f} \ltak{}{C_{1}} \pi_{f})\circ\com_{C_{1}}$. Since $\com_{B}$ and $\com_{C}$ have similar definitions, we only give the proof of one of the conditions, the other have a similar proof. For any $h \sm m_{1} \in A_{1}$, we have
\[
\com_{B_{2}}(\pi_{f}(h \sm m_{1})) = h_{(1)} \sm 1_{N} \rtak{}{B_{1}} h_{(2)} \sm f(m_{1}) = (\pi_{f} \rtak{}{B_{1}} \pi_{f})(\com_{B_{1}}(h \sm m_{1})).
\]

Now, consider a morphism $\pi: A_{1} \to A_{2}$ between $\mathcal{A}(N_{1},\theta_{1},\dtheta_{1},\mu_{1})$ and $\mathcal{A}(N_{2},\theta_{2},\dtheta_{2},\mu_{2})$ such that $\pi\circ\iota^{\dtheta_{1}}_{\G} = \iota^{\dtheta_{2}}_{\G}$. By definition of morphism in the unital case, we have $\pi(B_{1}) \subseteq B_{2}$, then for each $m \in N$, there is a unique  $m_{\pi} \in N$ such that $\pi(1_{\pG} \sm m) = 1_{\pG} \sm m_{\pi}$. Take the map $f_{\pi}: m \in N \mapsto m_{\pi} \in N$. Since $\pi$ is a homomorphism of \as-algebras, then by unicity the map $f_{\pi}$ is necessarily a homomorphism of \as-algebras. Since $\pi\circ\iota^{\dtheta_{1}}_{\G} = \iota^{\dtheta_{2}}_{\G}$, thus
\[
m_{[-1]} \sm f_{\pi}(m_{[0]}) = \pi(\beta_{1}(m)) = \pi(t_{B_{1}}(\alp_{1}(m))) = t_{B_{2}}(\pi(\alp_{1}(m))) = \beta_{2}(f_{\pi}(m)) = f_{\pi}(m)_{[-1]} \sm f_{\pi}(m)_{[0]},
\]
for all $m \in N$, i.e. $f_{\pi}$ is a $\Gc$-equivariant \as-homomorphism. Finally, because
\[
h_{(1)} \sm f_{\pi}(m \lhd_{\dtheta_{2}} h_{(2)}) = \pi(\alp(m)(h \sm 1_{N}) = \pi(\alp(m))\pi(h \sm 1_{N}) = h_{(1)} \sm (f_{\pi}(m) \lhd_{\dtheta_{2}} h_{(2)})
\]
then, by implying the counit of $\G$ on the first leg, we have $f_{\pi}(m \lhd_{\dtheta_{2}} h) = f_{\pi}(m) \lhd_{\dtheta_{2}} h$ for all $h \in \pG$ and $m \in N$. This last statement implies that $f_{\pi}$ is a $\dGo$-equivariant \as-homomorphism.
\fin

The next corollary shows that any measured multiplier Hopf \as-algebroid associated with a braided commutative measured Yetter--Drinfeld $\Gc$-\as-algebra have $\mathcal{A}(\G)$ as a multiplier Hopf sub-\as-algebroid.

\begin{coro}\label{cor:group_inside_MHAd}
Let $\G$ be an algebraic quantum group of compact type, $(N,\theta,\dtheta,\mu)$ be a unital braided commutative measured Yetter--Drinfeld $\Gc$-\as-algebra. Then, the linear map
\[
\begin{array}{lccc}
\jo_{\G}: & \pG & \to & \pG \sm_{\dtheta} N \\
& h & \mapsto & h \sm 1_{N}
\end{array}
\]
yields an injective morphism of measured multiplier Hopf \as-algebroids between $\mathcal{A}(\G)$ and $\mathcal{A}(N,\theta,\dtheta,\mu)$. Moreover, the total algebra of the measured multiplier Hopf \as-algebroid $\mathcal{A}(N,\theta,\dtheta,mu)$ is given by
\[
A = \mathrm{span}\{\jo_{\G}(h)\alp_{\dtheta}(m) : h \in \pG, m \in N\} = \mathrm{span}\{\alp_{\dtheta}(m)\jo_{\G}(h) : h \in \pG, m \in N\}.
\]
\end{coro}
\pr
The trivial \as-homomorphism $\uni_{N}: \C \to N$, $\lambda \mapsto \lambda\cdot 1_{N}$ is a morphism between the unital braided commutative Yetter--Drinfeld $\Gc$-\as-algebras $(\C,\tr_{\Gc},\tr_{\dGo})$ and $(N,\theta,\dtheta)$. By the proposition above, $\jo_{\G} = \id_{\G} \sm \uni_{N}: \pG \sm_{\tr} \C \to \pG \sm_{\dtheta} N$ yields a morphism between the measured multiplier Hopf \as-algebroids $\mathcal{A}(\G) = \mathcal{A}(\C,\tr_{\Gc},\tr_{\dGo},\id_{\C})$ and $\mathcal{A}(N,\theta,\dtheta,\mu)$. If, we have $\jo_{\G}(h) = 0$, then by applying the total left integral $\phi$, we have $\varphi_{\G}(h)\mu(n) = \phi(\jo_{\G}(h)\alp_{\dtheta}(n)) = 0$ for all $n \in N$. The injectivity of $\jo_{\G}$ follows from the faithfulness of the left Haar integral $\varphi_{\G}$ and because $\mu$ is a non-zero functional.

The last part of the corollary follows directly since $\jo_{\G}$ and $\alp_{\dtheta}$ are the canonical faithful embeddings associated with the smash product $\pG \sm N$.
\fin

\subsection{The Pontrjagin dual}

In this section, we will compute the Pontrjagin dual of the measured multiplier Hopf \as-algebroid $\mathcal{A}(N,\theta,\dtheta,\mu)$ construct in Theorem~\ref{theo:MHAd_YD}. An important difference with the previous section is that here, we will be dealing mainly with non-unital \as-algebras, and then we need to make use of multipliers and pairings as tools for computations.

Let $\G=(\pG,\com_{\G},\varphi_{\G})$ be an algebraic quantum group of compact type and $(N,\theta,\dtheta,\mu)$ be a unital braided commutative measured Yetter--Drinfeld $\Gc$-\as-algebra. From now on, denote by 
\[
\du{\mathcal{A}}(N,\theta,\dtheta,\mu) = (\hat{A},C,B,t^{-1}_{B},t^{-1}_{C},\dcom_{C},\dcom_{B},\mu_{C},\mu_{B},{}_{B}\du{\psi}_{B},{}_{C}\du{\phi}_{C})
\]
the Pontrjagin dual of the measured multiplier Hopf \as-algebroid $\mathcal{A}(N,\theta,\dtheta,\mu)$. If $\phi$ denotes the left total integral of $\mathcal{A}(N,\theta,\dtheta,\mu)$, recall that $\du{A} = A \cdot \phi = \phi \cdot A$ and the canonical faithful embeddings $\hat{\iota}_{B}: B \to \M(\du{A})$, $\hat{\iota}_{C}: C \to \M(\du{A})$ are defined by
\[
[\hat{\iota}_{B}(\alp(m))](\phi\cdot a) = \phi\cdot(a\beta(m)), \quad
(\phi\cdot a)[\hat{\iota}_{B}(\alp(m))] = \phi\cdot(a\alp(m))
\]
\[
[\hat{\iota}_{B}(\alp(m))](a\cdot\phi) = (a\beta(\hat\gamma_{\theta}(m)))\cdot\phi, \quad
(a\cdot\phi)[\hat{\iota}_{B}(\alp(m))] = (a\alp(\gamma_{\theta}(m)))\cdot\phi
\]
and
\[
[\hat{\iota}_{C}(\beta(m))](\phi\cdot a) = \phi\cdot(\beta(\gamma_{\theta}(m))a), \quad (\phi\cdot a)[\hat{\iota}_{C}(\beta(m))] = \phi\cdot(\alp(m)a),
\]
\[
[\hat{\iota}_{C}(\beta(m))](a\cdot\phi) = (\beta(m)a)\cdot\phi, \quad (a\cdot\phi)[\hat{\iota}_{C}(\beta(m))] = (\alp(\gamma_{\theta}(m))a)\cdot\phi,
\]
for all $m \in N$, $h \in \pG$ and $a \in A$, respectively. In the remainder of this section we want to compute in a explicit form the Pontrjagin dual of the measured multiplier Hopf \as-algebroid $\mathcal{A}(N,\theta,\dtheta,\mu)$. In order to do that, first we give a universal property for smash products \as-algebras in the non-unital case to be used in Proposition~\ref{prop:dual_smash}.

\begin{lemm}\label{lem:mor_smash}
Let $\G = (\pG,\com_{\G},\varphi_{\G})$ be an algebraic quantum group, $C$ be a non-degenerate \as-algebra endowed with a multiplier Hopf \as-algebra action $\lhd: C \odo \pG \to C$ and $D$ be a non-degenerate \as-algebra. Denote by $\iota_{\pG}:\pG \to \M(\pG\sm_{\lhd}C)$ and $\iota_{C}:C \to \M(\pG\sm_{\lhd}C)$ the two canonical faithful embeddings given by $h \mapsto h \sm_{\lhd} 1_{\M(C)}$ and $c \mapsto 1_{\M(\pG)} \sm_{\lhd} c$, respectively. Assume there are two injective non-degenerate \as-homomorphism $\jo_{1}:\pG \to \M(D)$ and $\jo_{2}:C \to \M(D)$ satisfying the relation, called the {\em smash relation},
\begin{align}\label{cond:morphism_smash}
\jo_{2}(c)\jo_{1}(h) = \jo_{1}(h_{(1)})\jo_{2}( c \lhd h_{(2)})
\end{align}
as element in $\M(D)$, for all $h\in \pG$ and $c \in C$. Then, there is a unique non-degenerate injective \as-homomorphism $\jo: \pG \sm_{\lhd} C \to \M(D)$ such that $\jo\circ\iota_{\pG} = \jo_{1}$ and $\jo\circ\iota_{C} = \jo_{2}$.
\end{lemm}
\pr
Take the linear map $\jo := \jo_{1} \sm_{\lhd} \jo_{2}: \pG \sm_{\lhd} C \to \M(D)$, $h \sm_{\lhd} c \mapsto \jo_{1}(h)\jo_{2}(c)$. The proof that this is the unique non-degenerate injective \as-homomorphism satisfying the statement is straightforward. 
\fin

The next proposition give the explicit structure of the \as-algebra $\du{A}$. The non involutive case is given by \cite[Proposition~5.4.1]{T17}.

\begin{prop}\label{prop:dual_smash}
The linear map
\[
\begin{array}{lccc}
\lhd: & C \odo \dpG & \to & C \\
& \beta(m) \otimes \ome  & \mapsto & \beta(m \lhd_{\theta} \ome)
\end{array}
\]
yields an action of the multiplier Hopf \as-algebra $(\dpG,\dcom_{\G})$ on $C$, making the \as- algebra $\du{A}$ isomorphic to the smash product \as-algebra $\dpG \sm_{\lhd} C$ via the non-degenerate linear map
\[
\begin{array}{lccc}
\Xi: & \du{A} & \to & \dpG \sm_{\lhd} C \\
& (\alp(m)(h \sm 1_{N}))\cdot\phi & \mapsto & (h\cdot\varphi_{\G}) \sm _{\lhd} \beta(\hat\gamma_{\theta}(m))
\end{array}.
\]
Moreover, the linear maps $\iota_{C}: C \to \M(\dpG \sm_{\lhd} C)$ and $\iota'_{B}: B \to \M(\dpG \sm_{\lhd} C)$, defined respectively by
\[
(\ome \sm_{\lhd} \beta(n))[\iota_{C}(\beta(m))] = \ome \sm_{\lhd} \beta(mn), \quad\quad
(\ome \sm_{\lhd} \beta(n))[\iota'_{B}(\alp(m))] = \ome m^{[-1]} \sm_{\lhd} \beta(nm^{[0]}) 
\]
%
for all $m \in N$ and $\ome \in \dpG$, are compatible with the canonical faithful embeddings $\hat{\iota}_{C}$ and $\hat{\iota}_{B}$, in the sense that it satisfy $\Xi \circ \hat{\iota}_{C} = \iota_{C}$ and $\Xi \circ \hat{\iota}_{B} = \iota'_{B}$. Here, we are using the Sweedler type leg notation $\dtheta(m) = m^{[-1]} \odo m^{[0]}$ for $m \in N$.
\end{prop}
\pr
To show that $\lhd$ yields an action of multiplier Hopf \as-algebra, it is only necessary check the compatibility with the involution, because the others conditions follow directly from the definition of the map $\lhd$ and from the fact that $\lhd_{\theta}$ is an action of the multiplier Hopf \as-algebra $(\dpG,\dcom^{\co}_{\G})$ on $N$. 
By condition~\ref{idpr:action_theta}, we have
\begin{align*}
(\beta(m) \lhd \ome)^{*} & = \beta(m \lhd_{\theta} \ome)^{*} = \beta(\hat\gamma_{\theta}((m \lhd_{\theta} \ome)^{*})) \\
& = \beta(\hat\gamma_{\theta}(m^{*} \lhd_{\theta} \dS^{-1}_{\G}(\ome)^{*})) = \beta(\hat\gamma_{\theta}(m^{*} \lhd_{\theta} \dS_{\G}(\ome^{*}))) \\
& = \beta(\hat\gamma_{\theta}(m^{*}) \lhd_{\theta} \dS^{-1}_{\G}(\ome^{*})) = \beta(\hat\gamma_{\theta}(m^{*})) \lhd \dS^{-1}_{\G}(\ome^{*}) \\
& = \beta(m)^{*} \lhd \dS_{\G}(\ome)^{*}
\end{align*}
for all $m \in N$ and $\ome \in \dpG$.

Now, define the linear map $\djo_{\G}: \dpG \to \du{A}$ by $h\cdot\varphi_{\G} \mapsto (h \sm 1_{N})\cdot\phi$ for all $h \in \pG$. Observe that, since $\phi$ is faithful then $\djo_{\G}$ is an injective map. Also, note that the map $\djo_{\G}$ satisfies
\begin{align*}
\Pe(g \sm m,\djo_{\G}(\ome)) = \djo_{\G}(\ome)(g \sm m) = \p(g,\ome)\mu(m),
\end{align*}
for all $\ome \in \dpG$, $g \in \pG$ and $m \in N$. In the next lines, we will proof that the map $\djo_{\G}$ defined above is a non-degenerate homomorphism of \as-algebras satisfying the smash relation~\eqref{cond:morphism_smash} with the canonical faithful embedding $\hat{\iota}_{C}$. We have
\begin{enumerate}[label=\textup{(\roman*)}]
\item \textbf{The map $\djo_{\G}$ is a homomorphism}. For that, fix $\ome = h\cdot\varphi_{\G}$, we have
\begin{align*}
\djo_{\G}(\ome) \brhd ((g \sm 1_{N})\alp(m)) & = (\id \odot {}_{C}\phi_{C})(\com_{C}(g \sm m)(1_{A} \odo (h \sm 1_{N}))) \\
& = (\id \odo {}_{C}\phi_{C})((g_{(1)} \sm m)\odo (g_{(2)}h \sm 1_{N})) \\
& = t_{C}({}_{C}\phi_{C}(g_{(2)}h \sm 1_{N}))(g_{(1)} \sm m) \\
& = \varphi_{\G}(g_{(2)}h)g_{(1)} \sm m \\
& = (\id \odo \p(\cdot,h\cdot\varphi_{\G}))\com_{\G}(g) \sm m \\
& = ((\ome \brhd g) \sm 1_{N})\alp(m)
\end{align*}
and
\begin{align*}
(\alp(m)(g \sm 1_{N})) \blhd \djo_{\G}(\ome) & = ({}_{B}\psi_{B} \odot \id)(\com_{B}(\alp(m)(g \sm 1_{N}))((h \sm 1_{N}) \odo 1_{A})) \\
& = ({}_{B}\psi_{B} \odot \id)(\alp(m)(g_{(1)}h \sm 1_{N}) \odo (g_{(2)} \sm 1_{N})) \\
& = (g_{(2)} \sm 1_{N})t_{B}({}_{B}\psi_{B}(\alp(m)(g_{(1)}h \sm 1_{N}))) \\
& = (\varphi_{\G}(g_{(1)}h)g_{(2)} \sm 1_{N})\beta(m) \\
& = ((\p(\cdot,h\cdot\varphi_{\G}) \odo \id)\com_{\G}(g) \sm 1_{N})\beta(m) \\
& = ((g \blhd \ome) \sm 1_{N})\beta(m)
\end{align*}
for all $g\in \pG$, $m\in N$. Thus, fixing $m\in N$ and $g\in \pG$, we have
\[
\djo_{\G}(\ome)\djo_{\G}(\ome') \brhd (g \sm m) = \djo_{\G}(\ome) \brhd ((\ome' \brhd g) \sm m) = (\ome\ome' \brhd g) \sm m= \djo_{\G}(\ome\ome') \brhd (g \sm m)
\]
for all $\ome,\ome' \in \dpG$. By faithfulness of the map $\brhd$, the claim follows.

\item \textbf{The map $\djo_{\G}$ is an \as-homomorphism}. From the item above, condition~\ref{eq:I1} and condition~\ref{eq:I2}, it follows
\[
\djo_{\G}(\ome) \brhd (\alp(m)(g \sm 1_{N})) = \alp(m)((\ome \brhd g) \sm 1_{N})
\]
for all $\ome \in \dpG$, $g \in \pG$ and $m \in N$. Now, fix $\ome \in \dpG$. By direct computation
\begin{align*}
\djo_{\G}(\ome^{*}) \brhd (\alp(m)(g \sm 1_{N})) & = \alp(m)((\ome^{*} \brhd g) \sm 1_{N}) \\
& = (\djo_{\G}(\dS^{-1}_{\G}(\ome)) \brhd (\alp(m)(g \sm 1_{N}))^{*})^{*}
\end{align*}
for all $\ome \in \dpG$, $g \in \pG$ and $m \in N$. Hence, in particular, we have
\begin{align*}
\djo_{\G}(\ome^{*}) \brhd S(a)^{*} = (\djo_{\G}(\dS^{-1}_{\G}(\ome)) \brhd S(a))^{*}
\end{align*}
for all $\ome \in \dpG$ and $a \in A$. On the other hand, if $\ome = h\cdot\varphi_{\G}$ for a certain $h \in \pG$, since $\dS^{-1}_{\G}(\ome) = \sigma_{\G}(S_{\G}(h))\cdot\varphi_{\G}$, it follows
\begin{align*}
\dS^{-1}(\djo_{\G}(\ome))(\alp(m)(g \sm 1_{N})) & = \djo_{\G}(\ome)(S^{-1}(\alp(m)(g \sm 1_{N}))) \\
& = \phi((S^{-1}_{\G}(g) \sm 1_{N})\beta(m)(h \sm 1_{N})) \\
& = \phi(\beta(m)(h\sigma_{\G}(S^{-1}_{\G}(g)) \sm 1_{N})) \\
& = \mu(m)\varphi_{\G}(h\sigma_{\G}(S^{-1}_{\G}(g))) \\
& = \mu(m)\varphi_{\G}(g\sigma_{\G}(S_{\G}(h))) \\
& = \phi(\alp(m)(g \sm 1_{N})(\sigma_{\G}(S_{\G}(h)) \sm 1_{N})) \\
& = \djo_{\G}(\dS^{-1}_{\G}(\ome))(\alp(m)(g \sm 1_{N}))
\end{align*}
for all $m \in N$ and $g \in \pG$, i.e. $\dS^{-1}(\djo_{\G}(\ome)) = \djo_{\G}(\dS^{-1}_{\G}(\ome))$. Then, we have
\begin{align*}
\djo_{\G}(\ome)^{*} \brhd S(a)^{*} & = (\dS^{-1}(\djo_{\G}(\ome)) \brhd S(a))^{*} = (\djo_{\G}(\dS^{-1}_{\G}(\ome)) \brhd S(a))^{*} = \djo_{\G}(\ome^{*}) \brhd S(a)^{*}
\end{align*}
for all $a \in A$. The claim follows from the faithfulness of the action $\brhd$.

\item \textbf{The \as-homomorphism $\djo_{\G}$ is non-degenerate.} Since $\djo$ is an \as-homomorphism, it is enough to prove that $\djo_{\G}(\dpG)\du{A} = \du{A}$. Given $h,g \in \pG$ and $n \in N$, take elements $e^{i},f^{i} \in \pG$ such that $h \odo g = \sum_{i}\com_{\G}(e^{i})(f^{i} \odo 1_{\G})$, then
\[
(h \sm 1_{N}) \odo (g \sm 1_{N})\alp(n) = \sum_{i}\com_{C}((e^{i} \sm 1_{N})\alp(n))((f^{i} \sm 1_{N}) \odo 1_{A})
\]
which implies
\small
\begin{align*}
(((h \sm 1_{N})\cdot\phi)(((g \sm 1_{N})\alp(n))\cdot\phi))(a)
& = \sum_{i} \phi({}_{C}\phi_{C}(a(e^{i} \sm 1)\alp(n))(f^{i} \sm 1)) \\
& = \sum_{i} \phi(a(e^{i} \sm 1)\alp(n){}_{C}\phi_{C}(f^{i} \sm 1)) \\
& = \sum_{i} \phi(a(e^{i} \sm 1)\alp(n))\varphi_{\G}(f^{i}) \\
& = \sum_{i} ((\varphi_{\G}(f^{i})(e^{i} \sm 1)\alp(n))\cdot\phi)(a)
\end{align*}
\normalsize
for all $a \in A$. In particular, we have
\small
\begin{align*}
((h \sm 1)\alp(n))\cdot\phi = ((h_{(1)} \sm 1)\cdot\phi)(((h_{(2)} \sm 1)\alp(n))\cdot\phi) = \djo_{\G}(h_{(1)}\cdot\varphi_{\G})(((h_{(2)} \sm 1)\alp(n))\cdot\phi)
\end{align*}
\normalsize
for all $h \in \pG$, and $n \in M$. Since the elements of the form $((h \sm 1)\alp(n))\cdot\phi$ span the \as-algebra $\du{A}$, then we have the result.

\item \textbf{The injective \as-homomorphism $\djo_{\G}$ and the canonical faithful embedding $\hat{\iota}_{C}$ satisfy the smash relation~\eqref{cond:morphism_smash}.} Fix $m \in N$ and $\ome = h\cdot\varphi_{\G}$ with $h \in \pG$. On one hand, we have
\begin{align*}
[\hat{\iota}_{C}(\beta(m))\djo_{\G}(\ome)](g \sm n) & = [(\beta(m)(h \sm 1_{N})\cdot\phi](g \sm n) \\
& = \phi((g \sm n)\beta(m)(h \sm 1_{N})) \\
& = \phi((\sigma^{-1}_{\G}(h)g \sm 1_{N})\alp(n)\beta(m)) \\
& = \phi((\sigma^{-1}_{\G}(h)g \sm 1_{N})\beta(m)\alp(n)) \\
& = \phi((\sigma^{-1}_{\G}(h)gm_{[-1]} \sm 1_{N})\alp(m_{[0]}n)) \\
& = \varphi_{\G}(\sigma^{-1}_{\G}(h)gm_{[-1]})\mu(m_{[0]}n) \\
& = \p(gm_{[-1]},\ome)\mu(m_{[0]}n),
\intertext{and on the other}
[\djo_{\G}(\ome)\hat\iota_{C}(\beta(m))](g \sm n) & = [(\alp(\gamma_{\theta}(m))(h \sm 1_{N}))\cdot\phi](g \sm n) \\
& = \phi((g \sm n)\alp(\gamma_{\theta}(m))(h \sm 1_{N})) \\
& = \phi((\sigma^{-1}_{\G}(h)g \sm 1_{N})\alp(n\gamma_{\theta}(m))) \\
& = \varphi_{\G}(\sigma^{-1}_{\G}(h)g)\mu(n\gamma_{\theta}(m)) \\
& = \varphi_{\G}(gh)\mu(mn) \\
& = \p(g,\ome)\mu(mn),
\end{align*}
for all $n\in N$, $g \in \pG$. Thus, for a fixed $m \in N$ and a fixed $\ome \in \dpG$, it holds
\begin{align*}
[\djo_{\G}(\ome_{(1)})\hat{\iota}_{C}(\beta(m) \lhd \ome_{(2)})](g \sm n) & = [\djo_{\G}(\ome_{(1)})\hat{\iota}_{C}(\beta(m \lhd_{\theta} \ome_{(2)}))](g \sm n) \\
& = \p(g,\ome_{(1)})\mu((m \lhd_{\theta} \ome_{(2)})n) \\
& = \p(g,\ome_{(1)})\p(m_{[-1]},\ome_{(2)})\mu(m_{[0]}n) \\
& = \p(gm_{[-1]},\ome)\mu(m_{[0]}n) \\
& = [\hat{\iota}_{C}(\beta(m))\djo_{\G}(\ome)](g \sm n)
\end{align*}
for all $n\in N$, $g \in \pG$. The claim follows from this last equality.
\end{enumerate}

Finally, by Lemma~\ref{lem:mor_smash}, that the map
\[
\begin{array}{lccc}
\djo_{\G} \sm_{\lhd} \hat\iota_{C}: & \dpG \sm_{\lhd} C & \to & \M(\du{A}) \\
& (h\cdot\varphi_{\G}) \sm_{\lhd} \beta(m) & \mapsto & \djo_{\G}(h\cdot\varphi_{\G})\hat\iota_{C}(\beta(m)) = (\alp(\gamma_{\theta}(m))(h \sm 1_{N}))\cdot\phi
\end{array}
\]
yields an injective \as-homomorphism. Moreover, since the elements of the form $(\alp(m)(1_{N} \sm h))\cdot\phi$ spam all the \as-algebra $\du{A}$, then the image of $\djo_{\G} \sm \iota_{B}$ is the \as-algebra $\du{A}$ and $\Xi = (\djo_{\G} \sm_{\lhd} \hat\iota_{C})^{-1}$ is an isomorphism of \as-algebras.

In order to proof the last statement of the proposition, first observe that by direct computation
\begin{align*}\label{eq:calculo_integral_izquierda}
\p(g,\p(S_{\G}(h_{(2)}),\ome)(h_{(1)}\cdot\varphi_{\G})) & = \p(g,h_{(1)}\cdot\varphi_{\G})\p(S_{\G}(h_{(2)}),\ome) \\
& = \p((\varphi_{\G} \odo S_{\G})((g \odo 1)\com_{\G}(h)),\ome) \\
& = \p((\varphi_{\G} \odo \id)(\com_{\G}(g)(h \odo 1)),\ome) \\
& = \p(g,(h\cdot\varphi_{\G})\ome),
\end{align*}
for all $h,g \in \pG$ and $\ome \in \dpG$. Then, by considering the extension \as-homomorphism of $\Xi$ to multipliers, we have
\begin{align*}
\Xi((\alp(n)(h \sm 1_{N}))\cdot\phi)\Xi([\hat{\iota}_{C}(\beta(m))]) & = \Xi(((\alp(n)(h \sm 1_{N}))\cdot\phi)[\hat{\iota}_{C}(\beta(m))]) \\
& = \Xi((\alp(\gamma_{\theta}(m)n)(h \sm 1_{N}))\cdot\phi) \\
& = (h\cdot\varphi_{\G}) \sm_{\lhd} \beta(m\hat\gamma_{\theta}(n)) \\
& = ((h\cdot\varphi_{\G}) \sm_{\lhd} \beta(\hat\gamma_{\theta}(n)))[\iota_{C}(\beta(m))] \\
& = \Xi((\alp(n)(h \sm 1_{N}))\cdot\phi)[\iota_{C}(\beta(m))]
\end{align*}
and
\begin{align*}
\Xi((\alp(n)(h \sm 1_{N}))\cdot\phi)\Xi([\hat{\iota}_{B}(\alp(m))]) & = \Xi(((\alp(n)(h \sm 1_{N}))\cdot\phi)[\hat{\iota}_{B}(\alp(m))]) \\
& = \Xi((\alp(n)(h \sm 1_{N})\alp(\gamma_{\dtheta}(m)))\cdot\phi) \\
& = \Xi((\alp(n(\gamma_{\theta}(m) \lhd_{\dtheta} S^{-1}_{\G}(h_{(2)})))(h_{(1)} \sm 1_{N}))\cdot\phi) \\
& = (h_{(1)}\cdot\varphi_{\G}) \sm_{\lhd} \beta(\hat\gamma_{\theta}(n)(m \lhd_{\dtheta} S_{\G}(h_{(2)}))) \\
& = (\p(S_{\G}(h_{(2)}),m^{[-1]})(h_{(1)}\cdot\varphi_{\G})) \sm_{\lhd} \beta(\hat\gamma_{\theta}(n)m^{[0]}) \\
& = (h\cdot\varphi_{\G})m^{[-1]} \sm_{\lhd} \beta(\hat\gamma_{\theta}(n)m^{[0]}) \\
& = ((h\cdot\varphi_{\G}) \sm_{\lhd} \beta(\hat\gamma_{\theta}(n)))[\iota'_{B}(\alp(m))] \\
& = \Xi((\alp(n)(h \sm 1_{N}))\cdot\phi)[\iota'_{B}(\alp(m))]
\end{align*}
for all $m,n \in N$ and $h \in \pG$. The last two equalities implies that $\Xi\circ\hat{\iota}_{C} = \iota_{C}$ and $\Xi\circ\hat{\iota}_{B} = \iota'_{B}$.
\fin

Similarly to Corollary~\ref{cor:group_inside_MHAd}, the next Corollary shows that the Pontrjagin dual of any measured multiplier Hopf \as-algebroid associated with a unital braided commutative Yetter--Drinfeld $\Gc$-\as-algebra have $\mathcal{A}(\dGo)$ as as a multiplier Hopf sub-\as-algebroid.

\begin{coro}\label{cor:dual_group_inside_dual_MHAd}
The linear map
\[
\begin{array}{lccc}
\djo_{\G}: & \dpG & \to & \du{A} \\
& h\cdot\varphi_{\G} & \mapsto & (h \sm 1_{N})\cdot\phi
\end{array}
\]
yields an injective morphism between the measured multiplier Hopf \as-algebroids $\mathcal{A}(\dGo)$ and $\du{\mathcal{A}}(N,\theta,\dtheta,\mu)$.
\end{coro}
\pr
In the proof of Proposition~\ref{prop:dual_smash}, we have shown that the linear map $\djo_{\G}$ is a non-degenerate injective \as-homomorphism. Moreover, observe that directly $\djo_{\G}(\C)B = B$, $\djo_{\G}(\C)C = C$, $\djo_{\G}\circ\id_{\C} = t_{B}\circ\djo_{\G}$ and $\djo_{\G}\circ\id_{\C} = t_{C}\circ\djo_{\G}$. Then, we only have to prove that the \as-homomorphism $\djo_{\G}$ is compatible with the left and right comultiplication. First, observe that by Proposition~\ref{prop:multiplier_algebra}, Corollary~\ref{cor:group_inside_MHAd} and Proposition~\ref{prop:dual_smash}, the \as-algebra $A$ and $\du{A}$ are given by
\begin{align*}
A & = \mathrm{span}\{\jo_{\G}(h)\alp(m) : h \in \pG, m \in N\} = \mathrm{span}\{\alp(m)\jo_{\G}(h) : h \in \pG, m \in N\} \\
& = \mathrm{span}\{\jo_{\G}(h)\beta(m) : h \in \pG, m \in N\} = \mathrm{span}\{\beta(m)\jo_{\G}(h) : h \in \pG, m \in N\},
\end{align*}
and
\[
\du{A} = \mathrm{span}\{\djo_{\G}(\ome)\hat{\iota}_{C}(\beta(m)) : \ome \in \pG, m \in N\} = \mathrm{span}\{\hat{\iota}_{C}(\beta(m))\djo_{\G}(\ome) : \ome \in \pG, m \in N\},
\]
respectively. Consider the canonical pairing $\Pe$ between $\mathcal{A}(N,\theta,\dtheta,\mu)$ and $\du{\mathcal{A}}(N,\theta,\dtheta,\mu)$, and fix $\ome,\ome'\in\dpG$. Since $\mu$ is $\dtheta$-invariant, we have
\[
\Pe(\jo_{\G}(h)\beta(m)) = \p(gm_{[-1]},\ome)\mu(m_{[0]}) = \p(g,\ome_{(1)})\mu(m \lhd_{\theta} \ome_{(2)}) = \p(g,\ome)\mu(m)
\]
for all $g \in \pG$ and $m \in N$. Then, it follows
\begin{align*}
\Pe^{2}(\jo_{\G}(g)\beta(m) \otimes \jo_{\G}(g')\alp(m') &,(\djo_{\G}(\ome') \otimes 1)\dcom_{C}(\djo_{\G}(\ome))) \\
& = \Pe^{2}(\com_{C}(\jo_{\G}(g)\beta(m))(1 \otimes \jo_{\G}(g)\alp(m')),\djo_{\G}(\ome') \otimes \djo_{\G}(\ome)) \\
& = \Pe^{2}(\jo_{\G}(g_{(1)})\beta(m) \otimes \jo_{\G}(g_{(2)}g')\alp(m')),\djo_{\G}(\ome') \otimes \djo_{\G}(\ome)) \\
& = \Pe(\jo_{\G}(g_{(1)})\beta(m),\djo_{\G}(\ome'))\Pe(\jo_{\G}(g_{(2)}g')\alp(m'),\djo_{\G}(\ome)) \\
& = \p(g_{(1)},\ome')\mu(m)\p(g_{(2)}g',\ome)\mu(m') \\
& = \p(g,\ome'\ome_{(1)})\mu(m)\p(g',\ome_{(2)})\mu(m') \\
& = \Pe(\jo_{\G}(g)\beta(m),\djo_{\G}(\ome'\ome_{(1)}))\Pe(\jo_{\G}(g')\alp(m'),\djo_{\G}(\ome_{(2)})) \\
& = \Pe^{2}(\jo_{\G}(g)\beta(m) \otimes \jo_{\G}(g')\alp(m'),(\djo_{\G} \otimes \djo_{\G})((\ome' \otimes 1)\dcom_{\G}(\ome))) \\
& = \Pe^{2}(\jo_{\G}(g)\beta(m) \otimes \jo_{\G}(g')\alp(m'),(\djo_{\G}(\ome') \otimes 1)(\djo_{\G} \otimes \djo_{\G})(\dcom_{\G}(\ome)))
\end{align*}
for all $m,m' \in N$ and $g,g' \in \pG$. By non-degeneracy of the pairing $\Pe^{2}$, the last equality implies
\[
(\djo_{\G}(\ome') \otimes 1_{\du{A}})\dcom_{C}(\djo_{\G}(\ome)) = (\djo_{\G}(\ome') \otimes 1_{\du{A}})(\djo_{\G} \otimes \djo_{\G})(\dcom_{\G}(\ome)).
\]
By multiplying elements of the form $\iota_{\C}(\beta(m)) \odo 1_{\du{A}}$ by the left to the last equation, we have
\[
(\iota_{C}(\beta(m))\djo_{\G}(\ome') \otimes 1_{\du{A}})\dcom_{C}(\djo_{\G}(\ome)) = (\iota_{C}(\beta(m))\djo_{\G}(\ome') \otimes 1_{\du{A}})(\djo_{\G} \otimes \djo_{\G})(\dcom_{\G}(\ome))
\]
for all $m \in N$ and $\ome,\ome' \in \dpG$. Since the elements of the form $\iota_{C}(\beta(m))\djo_{\G}(\ome')$ span the \as-algebra $\du{A}$, thence
\[
\dcom_{C}(\djo_{\G}(\ome)) = (\djo_{\G} \otimes \djo_{\G})(\dcom_{\G}(\ome))
\]
for all $\ome \in \dpG$. Follows from the fact that $\djo_{\G}$ is an \as-homomorphism and because $\dcom_{B} = (* \overline{\times} *)\circ\dcom_{C}\circ *$, that we also have $\dcom_{B}(\djo_{\G}(\ome)) = (\djo_{\G} \otimes \djo_{\G})(\dcom_{\G}(\ome))$ for all $\ome \in \dpG$.
\fin


\medskip

The next proposition give another description of the \as-algebra $\du{A}$ using this time the smash product \as-algebra associated with the conjugate action $\theta^{\con}$. Recall that $N^{\op}_{\hat\gamma_{\theta}}$ is the vector space $N$ with \as-algebra structure given by $m^{\op}n^{\op} := (nm)^{\op}$ and $(m^{\op})^{*} := \hat\gamma_{\theta}(m^{*})^{\op}$ for all $m,n \in N$. Since we have the relation $\theta\circ\hat\gamma_{\theta} = (S^{-2}_{\Gc} \odo \hat\gamma_{\theta})\circ\theta$, we can construct the conjugate action  $\theta^{\con}: N^{\op}_{\hat\gamma_{\theta}} \to \pG \odo N^{\op}_{\hat\gamma_{\theta}}$, $m^{\op} \mapsto ({}^{\op} \odo {}^{\op})\theta(m)$, with dual multiplier Hopf \as-algebra action $\lhd_{\theta^{\con}}: N^{\op}_{\hat\gamma_{\theta}} \odo \dpG \to N^{\op}_{\hat\gamma_{\theta}}$ given by
\[
m^{\op} \lhd_{\theta^{\con}} \ome = (\p(\cdot,\ome) \odo \id)\theta^{\con}(m^{\op}) = \p(m_{[-1]},\ome)(m_{[0]})^{\op} = (m \lhd_{\theta} \ome)^{\op}
\]
for all $\ome \in \dpG$ and $m \in N$. 

\begin{rema}
If $(N,\theta,\dtheta)$ is a unital braided commutative Yetter--Drinfeld $\Gc$-\as-algebra, by Proposition~\ref{prop:dual_bc_yd} the triplet $(N^{\op}_{\hat\gamma_{\theta}},\dtheta^{\con},\theta^{\con})$ yields a unital braided commutative Yetter--Drinfeld $\dGco$-\as-algebra. Explicitly, we have the braided commutative relation
\begin{align*}\label{eq:calculo_yd_dual}
(n^{\op} \lhd_{\theta^{\con}} m^{[-1]})(m^{[0]})^{\op}
& = m^{\op}n^{\op}
\end{align*}
for all $m, n \in N$. This braided commutative relation will be useful for the next proposition.
\end{rema}

\begin{prop}\label{prop:dual_smash_MHAd}
Consider, the map $\sigma_{N}: N^{\op}_{\hat\gamma_{\theta}} \to C$, $m^{\op} \mapsto \beta(m)$, and the right actions $\lhd_{\theta^{\con}}$ and $\lhd$. Then, the linear map
\begin{equation}\label{eq:equiv_map}
\begin{array}{lccc}
\id \sm_{\theta^{\con}} \sigma_{N}: & \dpG \sm_{\theta^{\con}} N^{\op}_{\hat\gamma_{\theta}} & \to & \dpG \sm_{\lhd} C \\
& \ome \sm m^{\op} & \mapsto & \ome \sm \beta(m)
\end{array}
\end{equation}
yields an isomorphism of smash product \as-algebras, making the \as-algebra $\du{A}$ isomorphic to the smash product \as-algebra $\dpG \sm_{\theta^{\con}} N^{\op}_{\hat\gamma_{\theta}}$ via
\[
\mathcal{T} := (\id \sm_{\theta^{\con}} \sigma_{N})^{-1}\circ\Xi: (\alp(m)(h \sm 1_{N}))\cdot\phi \mapsto (h\cdot\varphi_{\G}) \sm_{\theta^{\con}} \hat\gamma_{\theta}(m)^{\op}.
\]
Moreover, if we consider the linear maps  $\tilde{\iota}_{C}: C \to \M(\dpG \sm_{\theta^{\con}} N^{\op}_{\hat\gamma_{\theta}})$ and $\tilde{\iota}_{B}: B \to \M(\dpG \sm_{\theta^{\con}} N^{\op}_{\hat\gamma_{\theta}})$ given respectively by
\[
\tilde{\iota}_{C}(\beta(m)) = \iota_{N^{\op}_{\hat\gamma_{\theta}}}(m^{\op}), \quad \text{ and } \quad \tilde{\iota}_{B}(\alp(m)) = m^{[-1]} \sm_{\theta^{\con}} (m^{[0]})^{\op}
\]
for all $m \in N$, i.e.
\[
(\ome \sm_{\theta^{\con}} n^{\op})[\tilde{\iota}_{C}(\beta(m))] = \ome \sm_{\theta^{\con}} (mn)^{\op} \quad \text{ and } \quad (\ome \sm_{\theta^{\con}} n^{\op})[\tilde{\iota}_{B}(\alp(m))] = \ome m^{[-1]} \sm_{\theta^{\con}} (nm^{[0]})^{\op},
\]
for all $m,n \in N$ and $\ome \in \dpG$. Then, we have $\mathcal{T}\circ\hat{\iota}_{C} = \tilde{\iota}_{C}$ and $\mathcal{T}\circ\hat{\iota}_{B} = \tilde{\iota}_{B}$.
\end{prop}
\pr
Observe that if the map $\sigma_{N}$ is an $(\dpG,\dcom_{\G})$-equivariant isomorphism of \as-algebras, then $\id \sm \sigma_{N}$ yields an isomorphism of \as-algebras. Since $\beta$ is an anti-isomorphism of algebras, it is enough to prove that $\sigma_{N}$ is compatible with the involution and $(\dpG,\dcom_{\G})$-equivariant. Indeed, on one hand we have $\sigma_{N}((m^{\op})^{*}) = \sigma_{N}(\hat\gamma_{\theta}(m^{*})^{\op}) = \beta(\hat\gamma_{\theta}(m^{*})) = \beta(m)^{*} = \sigma_{N}(m^{\op})^{*}$ for all $m \in N$. And, on the other hand we have $\sigma_{N}(m^{\op}) \lhd \ome = \beta(m) \lhd \ome = \beta(m \lhd_{\theta} \ome) = \sigma_{N}(m^{\op} \lhd_{\theta^{\con}} \ome)$ for all $m \in N$ and $\ome \in \dpG$. 

Now, consider the \as-isomorphism $\mathcal{T} := (\id \sm_{\theta^{\con}} \sigma_{N})^{-1}\circ\Xi: \du{A} \to \dpG \sm_{\theta^{\con}} N^{\op}_{\hat\gamma_{\theta}}$. By Proposition~\ref{prop:dual_smash}, we have
\begin{align*}
\mathcal{T}((\alp(n)(h \sm 1_{N}))\cdot\phi)\mathcal{T}([\hat\iota_{C}(\beta(m))]) & = (\id \sm_{\theta^{\con}} \sigma_{N})^{-1}(\Xi((\alp(n)(h \sm 1_{N}))\cdot\phi)\Xi([\hat\iota_{C}(\beta(m))])) \\
& = (\id \sm_{\theta^{\con}} \sigma_{N})^{-1}((h\cdot\varphi_{\G}) \sm_{\lhd} \beta(m\hat\gamma_{\theta}(n))) \\
& = (h\cdot\varphi_{\G}) \sm_{\theta^{\con}} (m\hat\gamma_{\theta}(n))^{\op} \\
& = \mathcal{T}((\alp(n)(h \sm 1_{N}))\cdot\phi)[\tilde\iota_{C}(\beta(m))]
\end{align*}
and
\begin{align*}
\mathcal{T}((\alp(n)(h \sm 1_{N}))\cdot\phi)\mathcal{T}([\hat\iota_{B}(\alp(m))]) & = (\id \sm_{\theta^{\con}} \sigma_{N})^{-1}(\Xi((\alp(n)(h \sm 1_{N}))\cdot\phi)\Xi([\hat\iota_{B}(\alp(m))])) \\
& = (\id \sm_{\theta^{\con}} \sigma_{N})^{-1}((h\cdot\varphi_{\G})m^{[-1]} \sm_{\lhd} \beta(\hat\gamma_{\theta}(n)m^{[0]})) \\
& = (h\cdot\varphi_{\G})m^{[-1]} \sm_{\theta^{\con}} (\hat\gamma_{\theta}(n)m^{[0]})^{\op} \\
& = \mathcal{T}((\alp(n)(h \sm 1_{N}))\cdot\phi)[\tilde\iota_{B}(\alp(m))]
\end{align*}
for all $m,n \in N$ and $h \in \pG$. The last two equalities implies that $\mathcal{T}\circ\hat{\iota}_{C} = \tilde\iota_{C}$ and $\mathcal{T}\circ\hat{\iota}_{B} = \tilde\iota_{B}$.
\fin

\medskip

The next technical lemma will be useful in the proof of Theorem~\ref{theo:DUAL_MHAd_YD}.

\begin{lemm}\label{lem:non-degenerancy_p_P}
For all $h,g \in \pG$ and $m,n \in N$, it hold
\[
((h\cdot\varphi_{\G}) \odo 1)\dcom_{\G}(g\cdot\varphi_{\G}) = (h_{(1)}\cdot\varphi_{\G}) \odo ((g\sigma_{\G}(S_{\G}(h_{(2)})))\cdot\varphi_{\G})
\]
and
\[
(\hat\iota_{C}(\beta(m))\hat\jo_{\G}(h\cdot\varphi_{\G}) \odo 1)\dcom_{C}((\alp(n)(g \sm 1))\cdot\phi) = \hat\iota_{C}(\beta(m))\hat\jo_{\G}(h_{(1)}\cdot\varphi_{\G}) \odo (\alp(n)(g\sigma_{\G}(S_{\G}(h_{(2)})) \sm 1))\cdot\phi.
\]
\end{lemm}
\pr
Fix $h,g \in \pG$ and $m,n\in N$. By using the canonical pairing $\p$ and the Remark~\ref{rem:left_right_invariance}, we have
\begin{align*}
\p^{2}(h' \odo g',((h\cdot\varphi_{\G}) \odo 1)\dcom_{\G}(g\cdot\varphi_{\G})) & = \p^{2}(\com_{\G}(h')(1 \odo g'),(h\cdot\varphi_{\G}) \odo (g\cdot\varphi_{\G})) \\
& = (\varphi_{\G} \odo \varphi_{\G})(\com_{\G}(h')(h \odo g'g)) \\
& = \varphi_{\G}((\varphi_{\G} \odo \id)(\com_{\G}(h')(h \odo 1))g'g) \\
& = \varphi_{\G}(S_{\G}((\varphi_{\G} \odo \id)((h' \odo 1)\com_{\G}(h)))g'g) \\
& = (\varphi_{\G} \odo \varphi_{\G})((h' \odo 1)(\id \odo S_{\G})(\com_{\G}(h))(1 \odo g'g)) \\
& = (\varphi_{\G} \odo \varphi_{\G})(h'h_{(1)} \odo S_{\G}(h_{(2)})g'g) \\
& = \p(h',h_{(1)}\cdot\varphi_{\G})\p(g',(g\sigma_{\G}(S_{\G}(h_{(2)})))\cdot\varphi_{\G}) \\
& = \p^{2}(h' \odo g',(h_{(1)}\cdot\varphi_{\G}) \odo ((g\sigma_{\G}(S_{\G}(h_{(2)})))\cdot\varphi_{\G}))
\end{align*}
for all $h',g' \in \pG$. By non-degeneracy of the pairing $\p$, it follows
\[
((h\cdot\varphi_{\G}) \odo 1)\dcom_{\G}(g\cdot\varphi_{\G}) = (h_{(1)}\cdot\varphi_{\G}) \odo ((g\sigma_{\G}(S_{\G}(h_{(2)})))\cdot\varphi_{\G}).
\]

Similarly, by using the canonical pairing $\Pe$, we have
\begin{align*}
\Pe^{2}(\jo_{\G}(h')\beta(m') &\odo \jo_{\G}(g')\alp(n'),(\hat\jo_{\G}(h\cdot\varphi_{\G}) \odo 1)\dcom_{C}((\alp(n)(g \sm 1))\cdot\phi)) \\
& = \Pe^{2}(\com_{C}((h' \sm 1)\beta(m'))(1_{A} \odo (g' \sm n')),\hat\jo_{\G}(h\cdot\varphi_{\G}) \odo (\alp(n)(g \sm 1))\cdot\phi) \\
& = \Pe^{2}((h'_{(1)} \sm 1)\beta(m') \odo (h'_{(2)}g' \sm n'),\hat\jo_{\G}(h\cdot\varphi_{\G}) \odo (\alp(n)(g \sm 1))\cdot\phi) \\
& = \Pe((h'_{(1)} \sm 1)\beta(m'),\hat\jo_{\G}(h\cdot\varphi_{\G}))\Pe(h'_{(2)}g' \sm n', \odo (\alp(n)(g \sm 1))\cdot\phi) \\
& = \phi((h'_{(1)} \sm 1)\beta(m')(h \sm 1))\phi((h'_{(2)}g' \sm n')\alp(n)(g \sm 1)) \\
& = \varphi_{\G}(h'_{(1)}h)\mu(m')\varphi_{\G}(h'_{(2)}g'g)\mu(n'n) \\
& = (\varphi_{\G} \odo \varphi_{\G})(\com_{\G}(h')(h \odo g'g))\mu(m')\mu(n'n) \\
& = (\varphi_{\G} \odo \varphi_{\G})((h' \odo 1)(\id \odo S_{\G})(\com_{\G}(h))(1 \odo g'g))\mu(m')\mu(n'n) \\
& = \varphi_{\G}(h'h_{(1)})\mu(m')\varphi_{\G}(\sigma^{-1}_{\G}(g)S_{\G}(h_{(2)})g')\mu(n'n) \\
& = \phi((h' \sm 1)\beta(m')(h_{(1)} \sm 1))\phi((g'\sm 1)\alp(n'n)(g\sigma_{\G}(S_{\G}(h_{(2)})) \sm 1)) \\
& = \Pe^{2}(\jo_{\G}(h')\beta(m') \odo \jo_{\G}(g')\alp(n'),\hat\jo_{\G}(h_{(1)}\cdot\varphi_{\G}) \odo (\alp(n)(g\sigma_{\G}(S_{\G}(h_{(2)})) \sm 1))\cdot\phi)
\end{align*}
for all $h',g' \in \pG$ and $m',n' \in N$. By non-degeneracy of the pairing $\Pe$, it follows
\[
(\hat\jo_{\G}(h\cdot\varphi_{\G}) \odo 1)\dcom_{C}((\alp(n)(g \sm 1))\cdot\phi) = \hat\jo_{\G}(h_{(1)}\cdot\varphi_{\G}) \odo (\alp(n)(g\sigma_{\G}(S_{\G}(h_{(2)})) \sm 1))\cdot\phi,
\]
which implies
\[
(\hat\iota_{C}(\beta(m))\hat\jo_{\G}(h\cdot\varphi_{\G}) \odo 1)\dcom_{C}((\alp(n)(g \sm 1))\cdot\phi) = \hat\iota_{C}(\beta(m))\hat\jo_{\G}(h_{(1)}\cdot\varphi_{\G}) \odo (\alp(n)(g\sigma_{\G}(S_{\G}(h_{(2)})) \sm 1))\cdot\phi
\]
after having multiplied $\hat\iota_{C}(\beta(m)) \odo 1$ by the left.
\fin

The next theorem is our main result of this section. It shows that the Pontrjagin dual $\du{\mathcal{A}}(N,\theta,\dtheta,\mu)$ have a similar algebraic structure with the original measured multiplier Hopf \as-algebroid $\mathcal{A}(N,\theta,\dtheta,\mu)$.

\begin{theo}\label{theo:DUAL_MHAd_YD}
Let $\G=(\pG,\com_{\G},\varphi_{\G})$ be an algebraic quantum group of compact type and $(N,\theta,\dtheta,\mu)$ be a unital braided commutative measured Yetter--Drinfeld $\Gc$-\as-algebra with canonical automorphisms denoted by $\gamma_{\theta}$ and $\hat\gamma_{\theta}$. Keep the notations given in Theorem~\ref{theo:MHAd_YD}.
\begin{enumerate}[label=\textup{(\roman*)}]
\item Consider the unital braided commutative measured Yetter--Drinfeld $\dGo$-\as-algebra $(N^{\op}_{\hat\gamma_{\theta}},\dtheta^{\con},\theta^{\con},\mu^{\ops})$, the non-degenerate \as-algebra $A' = \dpG \sm_{\theta^{\con}} N^{\op}_{\hat\gamma_{\theta}}$, the vector subspaces
\[
B'=\{\tilde\iota_{C}(\beta_{\theta}(m)) =: 1 \sm m^{\op} : m \in N\},
\]
\[
C' = \{\tilde\iota_{B}(\alp_{\dtheta}(m)) =: m^{[-1]} \sm (m^{[0]})^{\op} : m \in N\},
\]
of $\M(A')$, and the automorphisms
\[
\begin{array}{lccc}
\gamma_{\dtheta^{\con}}: & N^{\op}_{\hat\gamma_{\theta}} & \to & N^{\op}_{\hat\gamma_{\theta}} \\
& m^{\op} & \mapsto & \hat\gamma_{\theta}(m)^{\op}
\end{array}
,\quad\quad
\begin{array}{lccc}
\hat\gamma_{\dtheta^{\con}}: & N^{\op}_{\hat\gamma_{\theta}} & \to & N^{\op}_{\hat\gamma_{\theta}} \\
& m^{\op} & \mapsto & \gamma_{\theta}(m)^{\op}
\end{array}.
\]
Hence, the linear maps
\small
\[
\begin{array}{lccc}
\alp_{\theta^{\con}}: & N^{\op}_{\hat\gamma_{\theta}} & \to & \M(A') \\
& m^{\op} & \mapsto & \tilde\iota_{C}(\beta_{\theta}(m))
\end{array},
\qquad
\begin{array}{lccc}
\beta_{\dtheta^{\con}}: & N^{\op}_{\hat\gamma_{\theta}} & \to & \M(A') \\
& m^{\op} & \mapsto & \tilde\iota_{B}(\alp_{\dtheta}(m))
\end{array}
\]
are injective and satisfy the following conditions
\begin{enumerate}[label=\textup{(C'\arabic*)}]
\item $\alp_{\theta^{\con}}(m^{\op})\alp_{\theta^{\con}}(n^{\op}) = \alp_{\theta^{\con}}((nm)^{\op})$ and $\alp_{\theta^{\con}}(m^{\op})^{*} = \alp_{\theta^{\con}}((m^{\op})^{*})$,

\item $\beta_{\dtheta^{\con}}(m^{\op})\beta_{\dtheta^{\con}}(n^{\op}) = \beta_{\dtheta^{\con}}((mn)^{\op})$ and $\beta_{\dtheta^{\con}}(m^{\op})^{*} = \beta_{\dtheta^{\con}}(\hat\gamma_{\dtheta^{\con}}((m^{\op})^{*})) = \beta_{\dtheta^{\con}}(\gamma_{\dtheta^{\con}}(m^{\op})^{*})$,

\item $\alp_{\theta^{\con}}(m^{\op})\beta_{\dtheta^{\con}}(n^{\op}) = \beta_{\dtheta^{\con}}(n^{\op})\alp_{\theta^{\con}}(m^{\op})$
\end{enumerate}
for every $m,n \in N$. Moreover, considering the linear maps
\[
\begin{array}{lccc}
t_{B'}: & B' & \to & C' \\
& \alp_{\theta^{\con}}(m^{\op}) & \mapsto & \beta_{\dtheta^{\con}}(m^{\op})
\end{array},
\qquad
\begin{array}{lccc}
t_{C'}: & C' & \to & B' \\
& \beta_{\dtheta^{\con}}(m^{\op}) & \mapsto & \alp_{\theta^{\con}}(\gamma_{\dtheta^{\con}}(m^{\op}))
\end{array},
\]
\[
\begin{array}{lccc}
\com_{B'}: & A' & \to & A' \rtak{}{B'} A' \\
& \ome \sm m^{\op} & \mapsto & (\ome_{(1)} \sm 1^{\op}_{N}) \rtak{}{B'} (\ome_{(2)} \sm m^{\op})
\end{array},
\]
\[
\begin{array}{lccc}
\com_{C'}: & A' & \to & A' \ltak{}{C'} A' \\
& \ome \sm m^{\op} & \mapsto & (\ome_{(1)} \sm 1^{\op}_{N}) \ltak{}{C'} (\ome_{(2)} \sm m^{\op})
\end{array},
\]
\[
\begin{array}{lccc}
S': & A' & \to & A' \\
& \ome \odo m^{\op} & \mapsto & \beta_{\dtheta^{\con}}(\hat\gamma_{\dtheta^{\con}}(m^{\op}))(\dS_{\G}(\ome) \sm 1^{\op}_{N})
\end{array},
\]
\[
\begin{array}{lccc}
\cou_{B'}: & A' & \to & B' \\
& \alp_{\theta^{\con}}(m^{\op})(\ome \odo 1^{\op}_{N}) & \mapsto & \alp_{\theta^{\con}}(m^{\op} \lhd_{\theta^{\con}} \ome)
\end{array},
\]
\[
\begin{array}{lccc}
\cou_{C'}: & A' & \to & C' \\
& (\ome \odo 1^{\op}_{N})\beta_{\dtheta^{\con}}(m^{\op}) & \mapsto & \beta_{\dtheta^{\con}}(m^{\op} \lhd_{\theta^{\con}} \dS_{\G}(\ome)) 
\end{array},
\]
\[
\begin{array}{lccc}
\mu_{B'}: & B' & \to & \C \\
& \alp_{\theta^{\con}}(m^{\op}) & \mapsto & \mu(m)
\end{array},
\quad
\begin{array}{lccc}
\mu_{C'}: & C' & \to & \C \\
& \beta_{\dtheta^{\con}}(m^{\op}) & \mapsto & \mu(m)
\end{array},
\]
\[
\begin{array}{lccc}
{}_{B'}\psi_{B'}: & A' & \to & B' \\
& \alp_{\theta^{\con}}(m^{\op})(\ome \odo 1^{\op}_{N}) & \mapsto & \dvar_{\G}(\ome)\alp_{\theta^{\con}}(m^{\op})
\end{array},
\]
\[
\begin{array}{lccc}
{}_{C'}\phi_{C'}: & A' & \to & C' \\
& \beta_{\dtheta^{\con}}(m^{\op})(\ome \odo 1^{\op}_{N}) & \mapsto & \dvar_{\G}(\ome)\beta_{\dtheta^{\con}}(m^{\op})
\end{array},
\]
\normalsize
the collection
\[
\mathcal{A}(N^{\op}_{\hat\gamma_{\theta}},\dtheta^{\con},\theta^{\con},\mu^{\ops}) := (A',B',C',t_{B'},t_{C'},\com_{B'},\com_{C'},\mu_{B'},\mu_{C'},{}_{B'}\psi_{B'},{}_{C'}\phi_{C'});
\]
yields a measured multiplier Hopf \as-algebroid.

\item The linear map
\[
\begin{array}{lccc}
\mathcal{T} = (\id \sm_{\theta^{\con}} \sigma_{N})^{-1}\circ\Xi: & \du{A} & \to & \dpG \sm_{\theta^{\con}} N^{\op}_{\hat\gamma_{\theta}} \\
& (\alp(m)(h \sm 1_{N}))\cdot\phi & \mapsto & (h\cdot\varphi_{\G}) \sm \hat\gamma_{\theta}(m)^{\op}
\end{array}
\]
yields an isomorphism between the measured multiplier Hopf \as-algebroids
\[
\du{\mathcal{A}}(N,\theta,\dtheta,\mu) \quad \text{ and } \quad \mathcal{A}(N^{\op}_{\hat\gamma_{\theta}},\dtheta^{\con},\theta^{\con},\mu^{\ops})
\]
satisfying
\[
\mathcal{T} \circ \dS = S' \circ \mathcal{T}.
\]

\item The bilinear map
\[
\begin{array}{lccc}
\Pe_{\dtheta,\theta^{\con}}: & \pG \sm_{\dtheta} N \times \dpG \sm_{\theta^{\con}} N^{\op}_{\hat\gamma_{\theta}}& \to & \C \\
& (h \sm m) \times (\ome \sm n^{\op}) & \mapsto & \p(h,\ome)\mu(nm)
\end{array}
\]
yields a pairing in the sense of \cite[Definition~4.16]{TVDW22} between the measured multiplier Hopf \as-algebroids
\[
\mathcal{A}(N,\theta,\dtheta,\mu) \quad \text{ and } \quad \mathcal{A}(N^{\op}_{\hat\gamma_{\theta}},\dtheta^{\con},\theta^{\con},\mu^{\ops}).
\]
\end{enumerate}
\end{theo}
\pr
\begin{enumerate}[label=\textup{(\roman*)}]
\item Since $\beta_{\theta}$ and $\alp_{\dtheta}$ are injectives, by Proposition~\ref{prop:dual_smash_MHAd} the maps $\alp_{\theta^{\con}}$ and $\beta_{\dtheta^{\con}}$ are also injectives. Fix $m,n \in N$, we have
\begin{enumerate}[label=\textup{(C'\arabic*)}]
\item $\alp_{\theta^{\con}}(m^{\op})\alp_{\theta^{\con}}(n^{\op}) = \tilde\iota_{C}(\beta(m)\beta(n)) = \alp_{\theta^{\con}}((nm)^{\op})$ and
\[
\alp_{\theta^{\con}}(m^{\op})^{*} = \tilde\iota_{C}(\beta(m)^{*}) = \tilde\iota_{C}(\beta(\hat\gamma_{\theta}(m^{*}))) = \alp_{\theta^{\con}}(\hat\gamma_{\theta}(m^{*})^{\op}) = \alp_{\theta^{\con}}((m^{\op})^{*}),
\]

\item $\beta_{\dtheta^{\con}}(m^{\op})\beta_{\dtheta^{\con}}(n^{\op}) = \hat\iota_{B}(\alp(m)\alp(n)) = \beta_{\dtheta^{\con}}((mn)^{\op})$ and
\[
\beta_{\dtheta^{\con}}(m^{\op})^{*} = \tilde\iota_{B}(\alp(m^{*})) = \beta_{\dtheta^{\con}}(\hat\gamma_{\dtheta^{\con}}((m^{\op})^{*})) = \beta_{\dtheta^{\con}}(\gamma_{\dtheta^{\con}}(m^{\op})^{*}),
\]

\item $\alp_{\theta^{\con}}(m^{\op})\beta_{\dtheta^{\con}}(n^{\op}) = \mathcal{T}(\hat\iota_{C}(\beta(m))\hat\iota_{B}(\alp(n))) = \mathcal{T}(\hat\iota_{B}(\alp(n))\hat\iota_{C}(\beta(m))) = \beta_{\dtheta^{\con}}(n^{\op})\alp_{\theta^{\con}}(m^{\op})$.
\end{enumerate}

Now, observe that the linear maps defined in item~(i) for the collection $\mathcal{A}(N^{\op}_{\hat\gamma_{\theta}},\dtheta^{\con},\theta^{\con},\mu^{\ops})$ have precisely a similar construction as those one given in Theorem~\ref{theo:MHAd_YD} for the measured multiplier Hopf \as-algebroid $\mathcal{A}(N,\theta,\dtheta,\mu)$. The main difference is that in the current case the \as-algebra $A'$ is non unital and we need to deal with multipliers of $A'$ in the proof. That is the reason why, instead of repeating a similar proot to Theorem~\ref{theo:MHAd_YD} using multipliers, the proof is left to the reader. We will see in item~(ii) that the \as-isomorphism $\mathcal{T}$ of Proposition~\ref{prop:dual_smash_MHAd} is indeed an isomorphism of measured multiplier Hopf \as-algebroids between $\mathcal{A}(N,\theta,\dtheta,\mu)$ and $\mathcal{A}(N^{\op}_{\hat\gamma_{\theta}},\dtheta^{\con},\theta^{\con},\mu^{\ops})$ compatible with the antipode of each measured multiplier Hopf \as-algebroid.

\item Recall that the measured multiplier Hopf \as-algebroid $\du{\mathcal{A}}(N,\theta,\dtheta,\mu)$ is given by
\[
(\du{A},C,B,t^{-1}_{B},t^{-1}_{C},\dcom_{C},\dcom_{B},\mu_{C},\mu_{B},{}_{B}\psi_{B},{}_{C}\phi_{C}).
\]
Then, a linear map $\pi: \du{A} \to A'$ yields an isomorphism from $\du{\mathcal{A}}(N,\theta,\dtheta,\mu)$ to $\mathcal{A}(N^{\op}_{\hat\gamma_{\theta}},\dtheta^{\con},\theta^{\con},\mu^{\ops})$, if and only if $\pi$ is a \as-isomorphism satisfying the conditions
\begin{enumerate}[label=(\arabic*)]
\item $\pi(\hat\iota_{C}(C))B' = B'$ and $\pi(\hat\iota_{B}(B))C' = C'$;
\item $\pi\circ\hat\iota_{B}\circ t^{-1}_{B} = t_{B'}\circ\pi\circ\hat\iota_{C}$ and $\pi\circ\hat\iota_{C}\circ t^{-1}_{C} = t_{C'}\circ\pi\circ\hat\iota_{B}$;
\item $\com_{B'}\circ\pi = (\rtak{\pi}{C})\circ\dcom_{C}$ and $\com_{C'}\circ\pi = (\ltak{\pi}{B})\circ\dcom_{B}$.
\end{enumerate}

\medskip

In order to proof that $\mathcal{T}$ is a isomorphism from $\du{\mathcal{A}}(N,\theta,\dtheta,\mu)$ to $\mathcal{A}(N^{\op}_{\hat\gamma_{\theta}},\dtheta^{\con},\theta^{\con},\mu^{\ops})$, we need to check each condition mentioned above.

\begin{enumerate}[label=(\arabic*)]
\item By Proposition~\ref{prop:dual_smash_MHAd}, $\mathcal{T}$ is an \as-isomorphism from $\du{A}$ to $A'$ such that $\mathcal{T}(\hat\iota_{C}(C)) = \tilde\iota_{C}(C) = B'$ and $\mathcal{T}(\hat\iota_{B}(B)) = \tilde\iota_{B}(B) = C'$.

\item By definition, we have
\begin{align*}
\mathcal{T}(\hat\iota_{B}(t^{-1}_{B}(\beta_{\theta}(m)))) & = \tilde\iota_{B}(\alp_{\dtheta}(m)) = \beta_{\dtheta^{\con}}(m^{\op}) \\
& = t_{B'}(\alp_{\theta^{\con}}(m^{\op})) = t_{B'}(\tilde\iota_{C}(\beta_{\theta}(m))) \\
& = t_{B'}(\mathcal{T}(\hat\iota_{C}(\beta_{\theta}(m))))
\end{align*}
and
\begin{align*}
\mathcal{T}(\hat\iota_{C}(t^{-1}_{C}(\alp_{\dtheta}(m)))) & = \tilde\iota_{B}(\beta_{\theta}(\hat\gamma_{\theta}(m))) = \alp_{\theta^{\con}}(\hat\gamma_{\theta}(m)^{\op}) \\
& = t_{C'}(\beta_{\dtheta^{\con}}(m^{\op})) = t_{C'}(\tilde\iota_{B}(\alp_{\dtheta}(m))) \\
& = t_{C'}(\mathcal{T}(\hat\iota_{B}(\alp_{\dtheta}(m)))),
\end{align*}
for all $m \in N$.

\item Fix $n \in N$ and $g \in \pG$. By Lemma~\ref{lem:non-degenerancy_p_P}, we have
\begin{align*}
(\mathcal{T}(\hat\iota_{C}&(\beta(m))\hat\jo_{\G}(h\cdot\varphi_{\G})) \odo 1)((\rtak{\mathcal{T}}{C})\circ\dcom_{C})((\alp(n)(g \sm 1))\cdot\phi) \\
& = (\mathcal{T} \odo \mathcal{T})((\hat\iota_{C}(\beta(m))\hat\jo_{\G}(h\cdot\varphi_{\G}) \odo 1)\dcom_{C}((\alp(n)(g \sm 1))\cdot\phi)) \\
& = \mathcal{T}(\hat\iota_{C}(\beta(m))\hat\jo_{\G}(h_{(1)}\cdot\varphi_{\G})) \odo \mathcal{T}((\alp(n)(g\sigma_{\G}(S_{\G}(h_{(2)})) \sm 1))\cdot\phi) \\
& = \mathcal{T}(\hat\iota_{C}(\beta(m))\hat\jo_{\G}(h_{(1)}\cdot\varphi_{\G})) \odo ((g\sigma_{\G}(S_{\G}(h_{(2)})))\cdot\varphi_{\G}) \sm \hat\gamma_{\theta}(n)^{\op} \\
& = \mathcal{T}(\hat\iota_{C}(\beta(m))\hat\jo_{\G}(h\cdot\varphi_{\G})\hat\jo_{\G}((g\cdot\varphi_{\G})_{(1)})) \odo ((g\cdot\varphi_{\G})_{(2)} \sm \hat\gamma_{\theta}(n)^{\op}) \\
& = (\mathcal{T}(\hat\iota_{C}(\beta(m))\hat\jo_{\G}(h\cdot\varphi_{\G})) \odo 1)(\com_{B'}((g\cdot\varphi_{\G}) \sm \hat\gamma_{\theta}(n)^{\op}) \\
& = (\mathcal{T}(\hat\iota_{C}(\beta(m))\hat\jo_{\G}(h\cdot\varphi_{\G})) \odo 1)(\com_{B'}\circ\mathcal{T})((\alp(n)(g \sm 1))\cdot\phi)
\end{align*}
for all $m \in N$ and $h \in \pG$. The last equality implies that
\[
((\rtak{\mathcal{T}}{C})\circ\dcom_{C})((\alp(n)(g \sm 1))\cdot\phi) = (\com_{B'}\circ\mathcal{T})((\alp(n)(g \sm 1))\cdot\phi).
\]
Since the elements of the form $(\alp(n)(g \sm 1))\cdot\phi$ span the \as-algebra $\du{A}$, thus we have
\[
(\rtak{\mathcal{T}}{C})\circ\dcom_{C} = \com_{B'}\circ\mathcal{T}.
\]
The equality $(\ltak{\mathcal{T}}{B})\circ\dcom_{B} = \com_{C'}\circ\mathcal{T}$ can be proven in a similar way.
\end{enumerate}

\smallskip

Finally, since $\du{A} = \mathrm{span}\{\djo_{\G}(\ome)\hat{\iota}_{C}(\beta(m)) : \ome \in \pG, m \in N\}$ and
\begin{align*}
(S' \circ \mathcal{T})(\hat\jo_{\G}(\ome)\hat\iota_{C}(\beta(m))) & = S'(\mathcal{T}(\hat\jo_{\G}(\ome))\mathcal{T}(\hat\iota_{C}(\beta(m)))) \\
& = S'((\ome \odo 1^{\op}_{N})\tilde\iota_{C}(\beta(m))) \\
& = \beta_{\dtheta^{\con}}(\hat\gamma_{\dtheta^{\con}}(m^{\op}))(\dS_{\G}(\ome) \sm 1^{\op}_{N}) \\
& = \tilde\iota_{B}(\alp(\gamma_{\theta}(m)))\mathcal{T}(\hat\jo_{\G}(\dS_{\G}(\ome))) \\
& = \mathcal{T}(\hat\iota_{B}(\alp(\gamma_{\theta}(m)))\hat\jo_{\G}(\dS_{\G}(\ome))) \\
& = \mathcal{T}(\dS(\hat\iota_{C}(\beta(m)))\dS(\hat\jo_{\G}(\ome))) \\
& = \mathcal{T}(\dS(\hat\iota_{C}(\beta(m)))\dS(\hat\jo_{\G}(\ome))) \\
& = (\mathcal{T}\circ\dS)(\hat\jo_{\G}(\ome)\hat\iota_{C}(\beta(m)))
\end{align*}
for all $\ome \in \dpG$ and $m \in N$, we have $\mathcal{T} \circ \dS = S' \circ \mathcal{T}$.

\item Consider the canonical pairing between $\mathcal{A}$ and its Pontrjagin dual $\du{\mathcal{A}}$, $\Pe: A \times \du{A} \to \C$, given by $\Pe(h \sm m,(g \sm n)\cdot\phi) = \phi((h \sm m)(g \sm n))$. By direct computation, we have
\begin{align*}
\Pe(h \sm m,\mathcal{T}^{-1}((g\cdot\varphi_{\G}) \sm n^{\op})) & = \Pe(h \sm m,(\alp(\gamma_{\theta}(n))(g \sm 1_{N}))\cdot\phi) \\
& = \phi((h \sm m)\alp(\gamma_{\theta}(n))(g \sm 1_{N})) \\
& = \phi(\sigma^{-1}_{\G}(g)h \sm m\gamma_{\theta}(n)) \\
& = \p(h,g\cdot\varphi_{\G})\mu(nm) \\
& = \Pe_{\dtheta,\theta^{\con}}(h \sm m,(g\cdot\varphi_{\G}) \sm n^{\op})
\end{align*}
for all $h,g \in \pG$ and $m,n \in N$. Since $\mathcal{T}$ is an isomorphism of measured multiplier Hopf \as-algebroids satisfying $\mathcal{T} \circ \dS = S' \circ \mathcal{T}$, then $\Pe_{\dtheta,\theta^{\con}}$ yields a pairing of multiplier Hopf \as-algebroids in the sense of \cite[Definition 4.16]{TVDW22}.
\end{enumerate}
\fin

\begin{exam}[Pontrjagin duals of finite transformation groupoids]\label{ex:MHAd_finite_classic_dual}
Let $G$ be a finite group acting by the left on a finite space $X$ and $\nu:X \to \R^{+}_{0}$ be a non-zero $G$-invariant function. Consider the unital braided commutative measured Yetter--Drinfeld $\Gc$-\as-algebra $(K(X),\theta,\hat\theta,\mu_{\nu})$ arising from the action of $G$ on $X$ (see Example~\ref{ex:bc_yd_finite_transformation_groupoid}). The Pontrjagin dual of the measured multiplier Hopf \as-algebroid $\mathcal{A}(K(X),\theta,\dtheta,\mu_{\nu})$ of Example~\ref{ex:MHAd_finite_classic} is the measured multiplier Hopf \as-algebroid $\mathcal{A}(K(X),\dtheta^{\con},\theta,\mu_{\nu})$ with structure
\[
\begin{array}{cccc}
A' = & \C[G]\,\#_{\theta^{\con}}\,K(X) & \cong  & \C[G \ltimes X] \\
& \lambda_{g}\,\#\,f & \mapsto & \displaystyle\sum_{x \in X} f(x) \lambda_{(g,x)}
\end{array},
\]
\[
\begin{array}{rccc}
\alpha_{\theta^{\con}}: & K(X) & \to & \C[G]\,\#_{\theta^{\con}}\,K(X) \\
& f & \mapsto & \displaystyle\sum_{x \in X} f(x)\lambda_{(e,x)} 
\end{array},
\quad\quad
\begin{array}{rccc}
\beta_{\dtheta^{\con}}: & K(X) & \to & \C[G]\,\#_{\theta^{\con}}\,K(X) \\
& f & \mapsto & \displaystyle\sum_{x \in X} f(x)\lambda_{(e,x)}
\end{array},
\]
\[
B := \{\alpha_{\dtheta}(f) = d^{\bullet}(f) : f \in K(X)\} \cong K((G \ltimes X)^{(0)}) \cong \{\beta_{\theta}(f) = r^{\bullet}(f) : f \in K(X)\} =: C,
\]
\[
t_{B} = t_{C} = \id_{K(X)}.
\]
Given $\lambda_{g} \sm f \in \C[G \ltimes X]$,
\begin{enumerate}
\item Since $A \rtak{}{B} A = A \ltak{}{C} A = K((G \ltimes X)^{(2)})$, we have
\begin{align*}
\com_{B}(\lambda_{g} \sm f)((g,x),(g',x')) = \com_{C}(\lambda_{g} \sm f)((g,x),(g',x')) & = (\lambda_{g} \sm 1_{X})(g,x)(\lambda_{g} \sm f)(g',x') \\
& = p_{(1)}(g)p_{(2)}(g')f(x') \\
& = p(gg')f(x') \\
& = (p \odo f)((g,x)(g',x'))
\end{align*}
for all $((g,x),(g',x')) \in (G \ltimes X)^{(2)}$.

\item It hold
\begin{align*}
S(p \odo f)(g,x) & = \beta_{\theta}(f)(S_{\G}(p) \odo 1_{X})(g,x) \\
& = r^{\bullet}(f)(g,x)p(g^{-1}) \\
& = (p \odo f)(g^{-1},g\cdot x) \\
& = (p \odo f)((g,x)^{-1}),
\end{align*}
for all $(g,x) \in G \ltimes X$, and
\begin{align*}
\cou_{B}(\alp(f)(p \odo 1_{X}))(e,x) & = \cou_{\G}(p)f(x) = p(e)f(x) = \cou_{\G}(p)f(e\cdot x) = \cou_{C}((p \odo 1_{X})\beta(f))(e,x), \\\\
\mu_{B}(\alp(f)) & = \mu_{\nu}(f) = \sum_{x \in X}\nu(x)f(x) = \mu_{\nu}(f) = \mu_{C}(\beta(f)), \\\\
{}_{B}\psi_{B}(\alp(f)(p \odo 1_{X}))(e,x) & = \sum_{g \in G}p(g)f(x) = \sum_{g \in G}p(g)f(e \cdot x) = {}_{C}\phi_{C}((p \odo 1_{X})\beta(f))(e,x)
\end{align*}
for all $(g,e) \in (G \ltimes X)^{(0)}$.
\end{enumerate}
\end{exam}

\begin{exam}[Pontrjagin duals associated with Fell bundles]\label{ex:MHAd_fell_dual}
Let $\fell{A}$ be a separable Fell bundle over a discrete group $G$ and $\nu: \fell{A}_{e} \to \C$ be a faithful state on the unital C*-algebra $\fell{A}_{e}$. Assume that $\fell{A}$ is equipped with a family of Fell bundle isomorphisms $\rho = \{\rho_{g}:\fell{A} \to \fell{A}^{\ad(g^{-1})}\}_{g\in G}$ verifying the conditions
\begin{enumerate}
\item $\rho_{e} = \id_{\fell{A}}$;
\item $\rho_{g}\circ\rho_{g'} = \rho_{g'g}$ for all $g,g' \in G$;
\item for any $g,h \in G$, we have $ab = \rho_{g^{-1}}(b)a$ for all $a \in \fell{A}_{g}$ and $b \in \fell{A}_{h}$.
\end{enumerate}

Applying an adapted version of \cite[Theorem~5.1]{KMQW10}, the total algebra of the Pontrjagin dual of the measured multiplier Hopf \as-algebroid $\mathcal{A}(\Gamma_{c}(G,\fell{A}),\theta_{\fell{A}},\hat\theta_{\fell{A},\rho},\mu_{\nu})$ (Example~\ref{ex:bc_yd_fell_bundle}) is the opposite \as-algebra of continuous compactly supported cross sections of the transformation Fell bundle $\fell{A} \times G$. Indeed, we have
\[
\begin{array}{cccccc}
A' = & K(G) \#_{\theta^{\con}_{\fell{A}}} \Gamma_{c}(G,\fell{A})^{\op}_{\hat\gamma_{\theta_{\fell{A}}}} & \cong & \Gamma_{c}(G \ltimes G, \fell{A} \times G)^{\op}_{\tau_{\rho}} \\
& p \# f^{\op} & \mapsto & (\gamma_{\theta_{\fell{A}}}(f) \boxtimes S(p))^{\op}: (g,g') \mapsto (\gamma_{\theta_{\fell{A}}}(f)(g)p(g'^{-1}),g')
\end{array}.
\]
Recall that $\fell{A} \times G$ is the Fell bundle over the transformation groupoid $G \ltimes G$ constructed as pullback of $\fell{A}$ through the groupoid homomorphism $\pi_{2}: G \ltimes G \to G$, $(g_{1},g_{2}) \mapsto g_{2}$ (cf. \cite{KMQW10}).
\end{exam}

\begin{exam}[Algebraic quantum groups of discrete type]\label{ex:MHAd_dqg}
Let $\G$ is an algebraic quantum group of compact type. Consider the trivial unital braided commutative measured Yetter--Drinfeld $\Gc$-\as-algebra $(\C,\tr_{\Gc},\tr_{\dGo},\id_{\C})$. Thus, the Pontrjagin dual of the measured multiplier Hopf \as-algebroid $\mathcal{A}(\C,\tr_{\Gc},\tr_{\dGo},\id_{\C})$ is the measured multiplier Hopf \as-algebroid $\mathcal{A}(\C,\tr_{\dGco},\tr_{\G},\id_{\C})$ given by
\[
A' = \dpG \sm_{\tr} \C = \dpG, \quad \alp = \beta: 1_{\C} \mapsto 1_{\M(\dpG)}, \\
\]
\[
\quad B' = \alp(\C) = \C\cdot1_{\dpG} = \beta(\C) = C',
\]
\[
t_{B'} = t_{C'} = \id_{\C}, \quad \com_{B'} = \com_{C'} = \dcom_{\G},
\]
\[
S' = \dS_{\G}, \quad \cou_{B'} = \cou_{C'} = \dcou_{\G},
\]
\[
\mu_{B'} = \mu_{C'} = \id_{\C}, \quad {}_{B'}\psi_{B'} = \psi = \du{\psi}_{\G}, \quad {}_{C'}\phi_{C'} = \phi = \dvar_{\G}.
\]
In other words, we have
\[
\mathcal{A}(\C,\tr_{\dGco},\tr_{\G},\id_{\C}) = \mathcal{A}(\dGo).
\]
\end{exam}

\begin{exam}[Heisenberg algebras as Pontrjagin duals of quantum transformation groupoids]
Let $\G$ be an algebraic quantum group of compact type and $\HH$ be a quantum subgroup defined by the surjective \as-morphism $\pi:\pG \to \pol(\HH)$.
Consider the unital braided commutative measured Yetter--Drinfeld $\Gc$-\as-algebra
\[
(\pol(\HH\backslash\G),\theta=(S^{-1}_{\G} \odo \id)\circ\Sigma\circ\com_{\G}|_{\pol(\HH\backslash\G)},\dtheta=\ad_{\Sigma(U^{*})}|_{\pol(\HH\backslash\G)},\cou_{\G}|_{\pol(\HH\backslash\G)})
\]
given by Example~\ref{ex:MHAd_quotient}. It was shown that the canonical automorphism for the mentioned Yetter--Drinfeld \as-algebra are given by $\gamma_{\theta} = S^{-2}_{\G}$, $\hat\gamma_{\theta} = S^{2}_{\G}$.

\begin{prop}
Let $\G$ be an algebraic quantum group of compact type and $\HH$ be a quantum subgroup with respect to the surjective \as-morphism $\pi:\pG \to \pol(\HH)$. The linear map
\[
\begin{array}{lccc}
\mathcal{F}_{\HH,\G}: & \dpG \sm_{\theta^{\con}} \pol(\HH\backslash\G)^{\op}_{\hat\gamma_{\theta}} & \to & \He(\G)^{\op}_{S^{2}_{\G} \sm \dS^{-2}_{\G}} \\
& \ome \sm h^{\op} & \mapsto & (h \sm_{\brhd} \dS^{-1}_{\G}(\ome))^{\op}
\end{array}
\]
is a faithful embedding of \as-algebras, and using this embedding, we have the identification
\[
\dpG \sm_{\theta^{\con}} \pol(\HH\backslash\G)^{\op}_{\hat\gamma_{\theta}} \cong (\pol(\HH\backslash\G) \sm_{\brhd} \dpG)^{\op}_{S^{2}_{\G} \sm \dS^{-2}_{\G}}.
\]
\end{prop}
\pr
Denote $\mathcal{F}_{\HH,\G}$ simply by $\mathcal{F}$. Recall that $\He(\G)$ is the smash product \as-algebra associated with the left action $\com_{\G}$ of $\G$ on $\pG$. Explicitly $\He(\G) = \pG \sm_{\brhd} \dpG$, where the left module map $\brhd:\dpG \odo \pG \to \pG$ is given by $\ome \brhd h = (\id \odo \p(\cdot,\ome))\com_{\G}(h) = h_{(1)}\p(h_{(2)},\ome)$ for all $\ome \in \dpG$ and $h \in \pG$. First observe that
\begin{align*}
h^{\op} \lhd_{\theta^{\con}} \ome & = (\p(\cdot,\ome) \odo \id)(\theta^{\con}(h^{\op})) \\
& = (\p(\cdot,\ome) \odo \id)({}^{\op} \odo {}^{\op})(S^{-1}_{\G} \odo \id)\Sigma\com_{\G}(h) \\
& = ((\p(S^{-1}_{\G}(h_{(2)}),\ome)h_{(1)})^{\op} \\
& = (\dS^{-1}_{\G}(\ome) \brhd h)^{\op}
\end{align*}
for all $h \in \pol(\HH\backslash\G)$ and $\ome \in \dpG$. Then, we have
\begin{align*}
\mathcal{F}((\ome \sm h^{\op})(\ome' \odo h'^{\op})) & = \mathcal{F}(\ome\ome'_{(1)} \sm (h^{\op} \lhd_{\theta^{\con}} \ome'_{(2)})h'^{\op}) \\
& = (h'(\dS^{-1}_{\G}(\ome'_{(2)}) \brhd h) \sm_{\brhd} \dS^{-1}_{\G}(\ome'_{(1)})\dS^{-1}_{\G}(\ome))^{\op} \\
& = \mathcal{F}(\ome \sm h^{\op})\mathcal{F}(\ome' \sm h'^{\op})
\end{align*}
and
\begin{align*}
\mathcal{F}((\ome \sm h^{\op})^{*}) & = \mathcal{F}(\ome^{*}_{(1)} \sm ((h^{\op})^{*} \lhd_{\theta^{\con}} \ome^{*}_{(2)})) \\
& = \mathcal{F}(\ome^{*}_{(1)} \sm (S^{2}_{\G}(h^{*})^{\op} \lhd_{\theta^{\con}} \ome^{*}_{(2)})) \\
& = \mathcal{F}(\ome^{*}_{(1)} \sm (\dS^{-1}_{\G}(\ome^{*}_{(2)}) \brhd S^{2}_{\G}(h^{*}))^{\op}) \\
& = ((\dS^{-1}_{\G}(\ome^{*}_{(2)}) \brhd S^{2}_{\G}(h^{*})) \sm_{\brhd} \dS^{-1}_{\G}(\ome^{*}_{(1)}))^{\op} \\
& = ((S^{-2}_{\G}(h) \sm_{\brhd} \dS_{\G}(\ome))^{*})^{\op} \\
& = (((S^{-2}_{\G} \sm_{\brhd} \dS^{2}_{\G})(h \sm_{\brhd} \dS^{-1}_{\G}(\ome)))^{*})^{\op} \\
& = ((h \sm_{\brhd} \dS^{-1}_{\G}(\ome))^{\op})^{*} \\
& = \mathcal{F}(\ome \sm h^{\op})^{*}
\end{align*}
for all $h,h' \in \pol(\HH\backslash\G)$ and $\ome,\ome' \in \dpG$. Finally, $\mathcal{F}$ is a faithful embedding since the linear map $\mathcal{F}': (\pol(\HH\backslash\G) \sm_{\brhd} \dpG)^{\op}_{S^{2}_{\G}\sm \dS^{-2}_{\G}} \to \dpG \sm_{\theta^{\con}} \pol(\HH\backslash\G)^{\op}_{S^{2}_{\G}}$ given by $h \sm_{\brhd} \ome \mapsto \dS_{\G}(\ome) \sm h^{\op}$ satisfies
\[
\mathcal{F}'\circ\mathcal{F} = \id_{\dpG \sm_{\theta^{\con}} \pol(\HH\backslash\G)^{\op}_{\hat\gamma_{\theta}}} \quad \text{ and } \quad \mathcal{F}\circ\mathcal{F}' = \id_{(\pol(\HH\backslash\G) \sm_{\brhd} \dpG)^{\op}_{S^{2}_{\G}\sm \dS^{-2}_{\G}}}.
\]
\fin

Now, using the proposition above, the total algebra of the Pontrjagin dual of the measured multiplier Hopf \as-algebroid
\[
\mathcal{A}(\pol(\HH\backslash\G),(S^{-1}_{\G} \odo \id)\circ\Sigma\circ\com_{\G}|_{\pol(\HH\backslash\G)},\ad_{\Sigma(U^{*})}|_{\pol(\HH\backslash\G)},\cou_{\G}|_{\pol(\HH\backslash\G)}),
\]
is given by
\[
A' = \dpG \sm_{\theta^{\con}} \pol(\HH\backslash\G)^{\op}_{\hat\gamma_{\theta}} \cong (\pol(\HH\backslash\G) \sm_{\brhd} \dpG)^{\op}_{S^{2}_{\G} \sm \dS^{-2}_{\G}}.
\]

In particular, if we consider the trivial algebraic quantum group $\HH_{\bullet} = (\C,\id_{\C},\id_{\C})$ of $\G$, then the total algebra of the Pontrjagin dual of the measured multiplier Hopf \as-algebroid
\[
\mathcal{A}(\pol(\G),(S^{-1}_{\G} \odo \id)\circ\Sigma\circ\com_{\G},\ad_{\Sigma(U^{*})},\cou_{\G}),
\]
is given by the opposite Heinseberg algebra
\[
A' \cong (\pol(\G) \sm_{\brhd} \dpG)^{\op}_{S^{2}_{\G} \sm \dS^{-2}_{\G}}.
\]
\end{exam}

\bigskip
\appendix

\section{The opposite, the  co-opposite  and the bi-opposite measured multiplier Hopf \as-algebroids of a measured multiplier Hopf \as-algebroid}

Let $\mathcal{A} = (A,B,C,t_{B},t_{C},\com_{B},\com_{C},\mu_{B},\mu_{C},{}_{B}\psi_{B},{}_{C}\phi_{C})$ be a measured multiplier Hopf \as-algebroid. Consider the linear maps
\small
\[
\begin{array}{lccc}
\Sigma_{A_{B},{}^{B}A}: & \rtak{A}{B} & \to & \rtak{(A^{\co})}{C} \\
& a \rtak{}{B} a' & \mapsto & a' \rtak{}{C} a
\end{array},
\;\;
\begin{array}{lccc}
\Sigma_{A^{C},{}_{C}A}: & \ltak{A}{C}& \to & \ltak{(A^{\co})}{B} \\
& a \ltak{}{C} a' & \mapsto & a' \ltak{}{B} a
\end{array}
\]
\normalsize
and
\small
\[
\begin{array}{lccc}
{\op} \rtak{}{B} {\op}: & \rtak{A}{B} & \to & \ltak{(A^{\op})}{C^{\op}} \\
& a \rtak{}{B} a' & \mapsto & a^{\op} \ltak{}{C^{\op}} a'^{\op}
\end{array},
\;\;
\begin{array}{lccc}
{\op} \ltak{}{C} {\op}: & \ltak{A}{C}& \to & \rtak{(A^{\op})}{B^{\op}} \\
& a \ltak{}{C} a' & \mapsto & a^{\op} \rtak{}{B^{\op}} a'^{\op}
\end{array}.
\]
\normalsize
Following \cite{TVD18}, the collections
\[
\mathcal{A}^{\co} := (A,C,B,t^{-1}_{B},t^{-1}_{C},\Sigma_{A_{B},{}^{B}A}\circ\com_{B},\Sigma_{A^{C},{}_{C}A}\circ\com_{C}),
\]
\[
\mathcal{A}^{\op} := (A^{\op},B^{\op},C^{\op},(t^{\ops}_{C})^{-1},(t^{\ops}_{B})^{-1},(\op \ltak{}{C} \op)\circ\com_{C},(\op \rtak{}{B} \op)\circ\com_{B})
\]
and
\[
\mathcal{A}^{\op,\co} := (A^{\op},C^{\op},B^{\op},t^{\ops}_{C},t^{\ops}_{B},(\op \ltak{}{B} \op)\circ\Sigma_{A^{C},{}_{C}A}\circ\com_{C},(\op \rtak{}{C} \op)\circ\Sigma_{A_{B},{}^{B}A}\circ\com_{B})
\]
yield measured multiplier Hopf \as-algebroids, called the {\em co-opposite}, the {\em opposite} and the {\em bi-opposite} measured multiplier Hopf \as-algebroid of $\mathcal{A}$, respectively.

\section{Measured multiplier Hopf \as-algebroids arising from unital braided commutative (left) Yetter--Drinfeld \as-algebras}

Let $\G=(\pG,\com_{\G},\varphi_{\G})$  be an algebraic quantum group of compact type. Let $(N,\theta,\dtheta)$ be a unital braided commutative (left) Yetter--Drinfeld $\Gc$-\as-algebra. Consider the left action $\rhd_{\dtheta}: \pG \odo N \to N$, defined by $h \rhd_{\dtheta} m = (\id \odo \p(h^{},\cdot))\dtheta(m)$, of $(\pG,\com_{\G})$ on $N$ induced by the action $\dtheta$ and the right coaction $\theta: N \to N \odo \pG^{\op}$ of $(\pG^{\op},\com^{\op}_{\G})$ on $N$. Keep in mind also the Sweedler type leg notation $m_{[0]} \odo m^{\op}_{[1]}:= \theta(m)$ for all $m \in N$. By an equivalent version of Proposition~\ref{prop:bc_yd_equivalence}, the triplet $(N,\rhd_{\dtheta},\theta)$ yields a unital braided commutative (left-right) Yetter--Drinfeld $(\pG,\com_{\G})$-\as-algebra, i.e. we have
\begin{equation}\label{eq:yd_condition_l}
\theta(h \rhd_{\dtheta} m) = (h_{(2)} \rhd_{\dtheta} m_{[0]}) \;\odo\; (h_{(3)}m_{[1]}S^{-1}_{\G}(h_{(1)}))^{\op},
\end{equation}
\begin{equation}\label{eq:bc_condition_l}
mn = (S_{\G}(m_{[1]}) \rhd_{\dtheta} n)m_{[0]} = n_{[0]}(n_{[1]}\rhd_{\dtheta} m)
\end{equation}
for any $h \in \pG$, $m,n \in N$.

\begin{theo}\label{theo:gamma_hatgamma_l}
Let $(N,\theta,\dtheta)$ be a unital braided commutative (left) Yetter--Drinfeld $\Gc$-\as-algebra. Consider the linear maps
\[
\begin{array}{lccc}
\gamma_{\theta}: & N & \to & N \\
& m & \mapsto & S^{-1}_{\G}(m_{[1]}) \rhd_{\dtheta} m_{[0]}
\end{array},
\quad
\begin{array}{lccc}
\hat\gamma_{\theta}: & N & \to & N \\
& m & \mapsto & S^{2}_{\G}(m_{[1]}) \rhd_{\dtheta} m_{[0]}
\end{array}
.
\]
It holds
\begin{enumerate}[label=\textup{(CA'\arabic*)}]
\item $* \circ \gamma_{\theta} = \hat\gamma_{\theta} \circ *$, $\theta\circ\gamma_{\theta} = (\gamma_{\theta} \odo S^{2}_{\Gc})\circ\theta$, and $\theta\circ\hat\gamma_{\theta} = (\hat\gamma_{\theta} \odo S^{-2}_{\Gc})\circ\theta$.
\end{enumerate}
Moreover, we have
\begin{enumerate}[label=\textup{(CA'\arabic*)},resume]
\item $\gamma_{\theta}(h \rhd_{\dtheta} m) =  S^{-2}_{\G}(h) \rhd_{\dtheta} \gamma_{\theta}(m)$ and $\hat\gamma_{\theta}(h \rhd_{\dtheta} m) = S^{2}_{\G}(h) \rhd_{\dtheta} \hat\gamma_{\theta}(m)$, for all $m \in N$, $h \in \pG$,
\item $S^{-1}_{\G}(m_{[1]}) \rhd_{\dtheta} nm_{[0]}  = \gamma_{\theta}(m)n$ and $S^{2}_{\G}(m_{[1]}) \rhd_{\dtheta} m_{[0]}n = n\hat\gamma_{\theta}(m)$ for all $m, n \in N$,
\item $\gamma_{\theta}(\hat{\gamma_{\theta}}(m)) = m$ and $\hat\gamma_{\theta}(\gamma_{\theta}(m)) = m,$ for all $m \in N$,
\item $\gamma_{\theta}(mn) = \gamma_{\theta}(m)\gamma_{\theta}(n)$, $\hat\gamma_{\theta}(mn) = \hat\gamma_{\theta}(m)\hat\gamma_{\theta}(n)$ for all $m, n \in N$,
\item $\gamma_{\theta}(\gamma_{\theta}(m^{*})^{*}) = m$ and $\hat\gamma_{\theta}(\hat\gamma_{\theta}(m^{*})^{*}) = m$ for all $m \in N$,
\item\label{eq:main} $\dtheta\circ\gamma_{\theta} = (\gamma_{\theta} \odo \dS^{2}_{\G})\circ\dtheta$ and $\dtheta\circ\hat\gamma_{\theta} = (\hat\gamma_{\theta} \odo \dS^{-2}_{\G})\circ\dtheta$,
\item $\gamma_{\theta}(\ome \rhd_{\theta} m) =  \dS^{2}_{\G}(\ome) \rhd_{\theta} \gamma_{\theta}(m)$ and $\hat\gamma_{\theta}(\ome \rhd_{\theta} m) = \dS^{-2}_{\G}(\ome) \rhd_{\theta} \hat\gamma_{\theta}(m)$, for all $m \in N$, $\ome \in \dpG$,
\item If $\mu$ is a $\dtheta$-invariant non-zero faithful functional on $N$, it holds $\mu = \mu\circ\gamma_{\theta} = \mu\circ\hat\gamma_{\theta}$ and $\mu(mn)=\mu(\gamma_{\theta}(n)m) = \mu(n\hat\gamma_{\theta}(m))$ for all $m, n \in N$, i.e. $\mu$ satisfies the weak KMS condition with modular automorphism $\hat\gamma_{\theta}$.
\end{enumerate}
\end{theo}
\pr
The proof is similar to the proof of Theorem~\ref{theo:gamma_hatgamma}, then it is left to the reader.
\fin
\begin{theo}\label{theo:MHAd_YD_l}
Let $\G=(\pG,\com_{\G},\varphi_{\G})$ be an algebraic quantum group of compact type and $(N,\theta,\dtheta,\mu)$ be a unital braided commutative measured (left) Yetter--Drinfeld $\Gc$-\as-algebra. Denote by $A_{\lambda}$ the smash product \as-algebra $N \sm_{\dtheta} \pG$. Then, the linear maps
\[
\begin{array}{lccc}
\alp: & N & \to & A_{\lambda} \\
& m & \mapsto & m \sm 1
\end{array},
\qquad
\begin{array}{lccc}
\beta: & N & \to & A_{\lambda} \\
& m & \mapsto & m_{[0]} \sm m_{[1]}
\end{array}
\]
satisfy the following conditions
\begin{enumerate}[label=\textup{(\arabic*)}]
\item $\alp(m)\alp(n) = \alp(mn)$ and $\alp(m)^{*} = \alp(m^{*})$,
\item $\beta(m)\beta(n) = \beta(nm)$ and $\beta(m)^{*} = \beta(\gamma_{\theta}(m)^{*}) = \beta(\hat\gamma_{\theta}(m^{*}))$,
\item $\alp(m)\beta(n) = \beta(n)\alp(m)$
\end{enumerate}
for every $m, n \in N$. Moreover, if we take $B_{\lambda} = \beta(N)$, $C_{\lambda} = \alp(N)$ and the following linear maps
\small
\[
\begin{array}{lccc}
t_{B_{\lambda}}: & B_{\lambda} & \to & C_{\lambda} \\
& \beta(m) & \mapsto & \alp(\gamma_{\theta}(m))
\end{array},
\qquad
\begin{array}{lccc}
t_{C_{\lambda}}: & C_{\lambda} & \to & B_{\lambda} \\
& \alp(m) & \mapsto & \beta(m)
\end{array},
\]
\[
\begin{array}{lccc}
\com_{B_{\lambda}}: & A_{\lambda} & \to & A_{\lambda} \rtak{}{B_{\lambda}} A_{\lambda} \\
& m \sm h & \mapsto & (m \sm h_{(1)}) \rtak{}{B_{\lambda}} (1_{N} \sm h_{(2)})
\end{array},
\,
\begin{array}{lccc}
\com_{C_{\lambda}}: & A_{\lambda} & \to & A_{\lambda} \ltak{}{C_{\lambda}} A_{\lambda} \\
& m \sm h & \mapsto & (m \sm h_{(1)}) \ltak{}{C_{\lambda}} (1_{N} \sm h_{(2)})
\end{array},
\]
\[
\begin{array}{lccc}
S: & A_{\lambda} & \to & A_{\lambda} \\
& m \sm h & \mapsto & (1_{N} \sm S_{\G}(h))\beta(\hat\gamma_{\theta}(m))
\end{array},
\]
\[
\begin{array}{lccc}
\cou_{B_{\lambda}}: & A_{\lambda} & \to & B_{\lambda} \\
& (1_{N} \sm h)\beta(m) & \mapsto & \cou_{\G}(h)\beta(m)
\end{array},
\qquad
\begin{array}{lccc}
\cou_{C_{\lambda}}: & A_{\lambda} & \to & C_{\lambda} \\
& \alp(m)(1_{N} \sm h) & \mapsto & \cou_{\G}(h)\alp(m)
\end{array},
\]
\[
\begin{array}{lccc}
\mu_{B_{\lambda}}: & B_{\lambda} & \to & \C \\
& \beta(m) & \mapsto & \mu(m)
\end{array},
\qquad
\begin{array}{lccc}
\mu_{C_{\lambda}}: & C_{\lambda} & \to & \C \\
& \alp(m) & \mapsto & \mu(m)
\end{array},
\]
\[
\begin{array}{lccc}
{}_{B_{\lambda}}\psi_{B_{\lambda}}: & A_{\lambda} & \to & B_{\lambda} \\
& \beta(m)(1_{N} \sm h) & \mapsto & \varphi_{\G}(h)\beta(m)
\end{array},
\qquad
\begin{array}{lccc}
{}_{C_{\lambda}}\phi_{C_{\lambda}}: & A_{\lambda} & \to & C_{\lambda} \\
& \alp(m)(1_{N} \sm h)  & \mapsto & \varphi_{\G}(h)\alp(m)
\end{array},
\]
\normalsize
then the collection
\[
\mathcal{A}_{\lambda}(N,\theta,\dtheta,\mu) := (A_{\lambda},B_{\lambda},C_{\lambda},t_{B_{\lambda}},t_{C_{\lambda}},\com_{B_{\lambda}},\com_{C_{\lambda}},\mu_{B_{\lambda}},\mu_{C_{\lambda}},{}_{B_{\lambda}}\psi_{B_{\lambda}},{}_{C_{\lambda}}\phi_{C_{\lambda}})
\]
yields a unital measured multiplier Hopf \as-algebroid.
\end{theo}
\pr
The proof is similar to the proof of Theorem~\ref{theo:MHAd_YD}, then it is left to the reader.
\fin

In what follows, we will show that the constructions given in Theorem~\ref{theo:MHAd_YD} and Theorem~\ref{theo:MHAd_YD_l} are indeed equivalent. Let $\G = (\pG,\com_{\pG},\varphi_{\G})$ be an algebraic quantum group of compact type, $N$ be a unital \as-algebras, $\theta: N \to N \odo \pG^{\op}$ be a left action of $\Gc$ on $N$, and $\dtheta:N \to \M(N \odo \dpG)$ be a left action of $\dGo$ on $N$.

\begin{prop}
The following statements are equivalent
\begin{enumerate}[label=\textup{(\roman*)}]
\item $(N,\theta,\dtheta,\mu)$ is a unital braided commutative measured (left) Yetter--Drinfeld $\Gc$-\as-algebra with canonical automorphisms $\gamma_{\theta}$ and $\hat\gamma_{\theta}$.
\item $(N^{\op}_{\gamma_{\theta}},\theta':=(\Sigma\circ\theta)^{\ops},\dtheta':=(\Sigma\circ\dtheta)^{\ops},\mu^{\ops})$ is a unital braided commutative measured (right) Yetter--Drinfeld $\Gc$-\as-algebra with canonical automorphisms $\gamma_{\theta'} :=  (\gamma_{\theta})^{\ops}$ and $\hat\gamma_{\theta'} := (\hat\gamma_{\theta})^{\ops}$.
\end{enumerate}
\end{prop}
\pr
The equivalence between the braided commutativity Yetter--Drinfeld conditions follows directly from \cite[Proposition~B.5]{Ta22_1}. We show the (i) implies (ii), the converse can be shown in a similar way. Recall that
\[
\theta' = (\Sigma\theta)^{\ops} =  ({}^{\op}\circ S_{\G} \odo {\;}^{\op})\circ\Sigma\circ\theta\circ{}^{\op}: N^{\op}_{\gamma_{\theta}} \to \pG^{\op} \odo N^{\op}_{\gamma_{\theta}},
\]
\[
\dtheta' = (\Sigma\dtheta)^{\ops} = (\dS^{-1}_{\G} \odo {\;}^{\op})\circ\Sigma\circ\dtheta\circ{}^{\op}: N^{\op}_{\gamma_{\theta}} \to \dpG \odo N^{\op}_{\gamma_{\theta}}.
\]
On one hand we have 
\begin{align*}
m^{\op} \lhd_{\dtheta'} h & = (\p(h,\cdot) \odo \id)\circ\dtheta'(m^{\op}) \\
& = (\p(h,\cdot) \odo \id)(\dS^{-1}_{\G} \odo {\;}^{\op})\circ\Sigma\circ\dtheta(m) \\
& = ((\id \odo \p(S^{-1}_{\G}(h),\cdot))(\dtheta(m))^{\op} \\
& = (S^{-1}_{\G}(h) \rhd_{\dtheta} m)^{\op},
\end{align*}
for all $m \in N$ and $h \in \pG$, and on the other,
\small
\begin{align*}
(m^{\op})^{\op}_{[-1]'} \odo (m^{\op})_{[0]'} & = \theta'(m^{\op}) = (\Sigma\circ\theta)^{\ops}(m^{\op}) = ({}^{\op}\circ S_{\G} \odo {\;}^{\op})\circ\Sigma\circ\theta(m) = (S_{\G}(m_{[1]}))^{\op} \odo (m_{[0]})^{\op},
\end{align*}
\normalsize
for all $m \in N$. Then
\begin{align*}
\gamma_{\theta'}(m^{\op}) & = (m^{\op})_{[0]'} \lhd_{\dtheta'} S^{-1}_{\G}((m^{\op})_{[-1]'}) \\
& = (m_{[0]})^{\op} \lhd_{\dtheta'} m_{[1]} = (S^{-1}_{\G}(m_{[1]}) \rhd_{\dtheta} m_{[0]})^{\op} \\
& = \gamma_{\theta}(m)^{\op} 
\end{align*}
and
\begin{align*}
\hat\gamma_{\theta'}(m^{\op}) & = (m^{\op})_{[0]'} \lhd_{\dtheta'} S^{2}_{\G}((m^{\op})_{[-1]'}) \\
&  = (m_{[0]})^{\op} \lhd_{\dtheta'} S^{3}_{\G}(m_{[1]}) = (S^{2}_{\G}(m_{[1]}) \rhd_{\dtheta} m_{[0]})^{\op} \\
& = \hat\gamma_{\theta}(m)^{\op} 
\end{align*}
for all $m \in N$. Finally, it is evident that $\mu$ is a Yetter--Drinfeld integral for $(N,\theta,\dtheta)$ if and only if $\mu^{\ops}$ is a Yetter--Drinfeld integral for $(N^{\op}_{\gamma_{\theta}},\theta':=(\Sigma\circ\theta)^{\ops},\dtheta':=(\Sigma\circ\dtheta)^{\ops})$.
\fin

The next proposition shows the relation between the smash products \as-algebras associated with $\dtheta$ and $\dtheta'$.

\begin{prop}
The linear map
\[
\begin{array}{lccc}
\mathcal{L}:& \pG \sm_{\dtheta'} N^{\op}_{\gamma_{\theta}} & \to & (N \sm_{\dtheta} \pG)^{\op}_{(\hat\gamma_{\theta} \odo S^{2}_{\G})} \\
& h \sm m^{\op} & \mapsto & (\hat\gamma_{\theta}(m) \sm S_{\G}(h))^{\op}
\end{array}
\]
yields an isomorphism of \as-algebras.
\end{prop}
\pr
We have
\begin{align*}
\mathcal{L}((h \sm m^{\op})(h' \sm n^{\op})) & = \mathcal{L}(hh'_{(1)} \sm (m^{\op} \lhd_{\dtheta'} h'_{(2)})^{\op}n^{\op}) \\
& = \mathcal{L}(hh'_{(1)} \sm (n(S^{-1}_{\G}(h'_{(2)}) \rhd_{\dtheta} m))^{\op}) \\
& = (\hat\gamma_{\theta}(n)\hat\gamma_{\theta}(S^{-1}_{\G}(h'_{(2)}) \rhd_{\dtheta} m) \sm S_{\G}(hh'_{(1)}))^{\op} \\
& = (\hat\gamma_{\theta}(n)(S_{\G}(h'_{(2)}) \rhd_{\dtheta} \hat\gamma_{\theta}(m)) \sm S_{\G}(h'_{(1)})S_{\G}(h))^{\op} \\
& = (\hat\gamma_{\theta}(m) \sm S_{\G}(h))^{\op}(\hat\gamma_{\theta}(n) \sm S_{\G}(h'))^{\op} \\
& = \mathcal{L}(h \sm m^{\op})\mathcal{L}(h' \sm n^{\op})
\end{align*}
and
\begin{align*}
\mathcal{L}((h \sm m^{\op})^{*}) & = \mathcal{L}(h^{*}_{(1)} \sm (m^{\op})^{*} \lhd_{\dtheta'} h^{*}_{(2)}) \\
& = \mathcal{L}(h^{*}_{(1)} \sm \gamma_{\theta}(m^{*})^{\op} \lhd_{\dtheta'} h^{*}_{(2)}) \\
& = \mathcal{L}(h^{*}_{(1)} \sm (S^{-1}_{\G}(h^{*}_{(2)}) \rhd_{\dtheta} \gamma_{\theta}(m^{*}))^{\op}) \\
& = \mathcal{L}(h^{*}_{(1)} \sm \gamma_{\theta}(S_{\G}(h^{*}_{(2)}) \rhd_{\dtheta} m^{*})^{\op}) \\
& = ((S_{\G}(h^{*}_{(2)}) \rhd_{\dtheta} m^{*}) \sm S_{\G}(h^{*}_{(1)}))^{\op} \\
& = ((S^{-1}_{\G}(h)^{*}_{(1)} \rhd_{\dtheta} m^{*}) \sm S^{-1}_{\G}(h)^{*}_{(2)})^{\op} \\
& = ((m \sm S^{-1}_{\G}(h))^{*})^{\op} \\
& = (((\gamma_{\theta} \# S^{-2}_{\G})(\gamma_{\theta}(m) \sm S_{\G}(h)))^{*})^{\op} \\
& = ((\hat\gamma_{\theta} \# S^{2}_{\G})((\hat\gamma_{\theta}(m) \sm S_{\G}(h))^{*}))^{\op} \\
& = (\mathcal{L}(h \sm m^{\op}))^{*}
\end{align*}
for all $m, n \in N$ and $h, h' \in \pG$. For the last equality, first observe that  $(\hat\gamma_{\theta} \# S^{2}_{\G})\circ * = * \circ (\gamma_{\theta} \# S^{-2}_{\G})$ on $N \sm_{\dtheta} \pG$, thus the anti-linear map
\[
(m \sm h)^{\op} \mapsto ((\hat\gamma_{\theta} \sm S^{2}_{\G})((m \sm h)^{*}))^{\op}
\]
gives the involution of the $(\hat\gamma_{\theta} \odo S^{2}_{\G})$-opposite \as-algebra $(N \sm_{\dtheta} \pG)^{\op}_{(\hat\gamma_{\theta} \odo S^{2}_{\G})}$.
\fin

It is not hard to see that the linear map $(\hat\gamma_{\theta})^{\ops}: N^{\op}_{\gamma_{\theta}} \to N^{\op}_{\hat\gamma_{\theta}}$, $m^{\op} \mapsto \hat\gamma_{\theta}(m)^{\op}$ is an isomorphism of \as-algebras, indeed we have
\begin{align*}
(\hat\gamma_{\theta})^{\ops}((m^{\op})^{*}) = (\hat\gamma_{\theta})^{\ops}(\gamma_{\theta}(m^{*})^{\op}) = (m^{*})^{\op}, \qquad (\hat\gamma_{\theta})^{\ops}(m^{\op})^{*} = (\hat\gamma_{\theta}(m)^{\op})^{*} = \hat\gamma_{\theta}(\hat\gamma_{\theta}(m)^{*})^{\op} = (m^{*})^{\op}
\end{align*}
for all $m \in N$.

\begin{prop}\label{prop:iso1}
Consider the linear maps associated with the measured multiplier Hopf \as-algebroid structure of $\mathcal{A}_{\lambda}(N,\theta,\dtheta,\mu)$, $\alp^{\ops}_{\dtheta}: N^{\op}_{\hat\gamma_{\theta}} \to (N \sm_{\dtheta} \pG)^{\op}_{(\hat\gamma_{\theta} \odo S^{2}_{\G})}$ and $\beta^{\ops}_{\theta}: N^{\op}_{\hat\gamma_{\theta}} \to (N \sm_{\dtheta} \pG)^{\op}_{(\hat\gamma_{\theta} \odo S^{2}_{\G})}$, then we have
\begin{enumerate}[label=\textup{(\roman*)}]
\item $\mathcal{L} \circ \alp_{\dtheta'} = \alp^{\ops}_{\dtheta} \circ (\hat\gamma_{\theta})^{\ops}$;
\item $\mathcal{L} \circ \beta_{\theta'} = \beta^{\ops}_{\theta} \circ (\hat\gamma_{\theta})^{\ops}$.
\end{enumerate}
\end{prop}
\pr
Fix $m \in N$. By definition $\mathcal{L}(\alp_{\dtheta'}(m^{\op})) = (\hat\gamma_{\theta}(m) \sm 1_{\pG})^{\op} = \alp^{\ops}_{\dtheta}((\hat\gamma_{\theta})^{\ops}(m^{\op}))$ for all $m \in N$. On the other hand,
\begin{align*}
\beta^{\ops}_{\theta}((\hat\gamma_{\theta})^{\ops}(m^{\op})) & = \beta_{\theta}(\hat\gamma_{\theta}(m))^{\op} = (\hat\gamma_{\theta}(m_{[0]}) \sm S^{2}_{\G}(m_{[1]}))^{\op} \\
& = \mathcal{L}(S_{\G}(m_{[1]}) \sm (m_{[0]})^{\op}) \\
& = \mathcal{L}((m^{\op})_{[-1]'} \sm (m^{\op})_{[0]'}) \\
& = \mathcal{L}(\beta_{\theta'}(m^{\op}))
\end{align*}
for all $m \in N$.
\fin

Fix the following notations $A_{\lambda}=N \sm_{\dtheta} \pG$, $C_{\lambda}=\alp_{\dtheta}(N)$, $B_{\lambda}=\beta_{\theta}(N)$, $A =  \pG \sm_{\dtheta'} N^{\op}_{\gamma_{\theta}}$, $B=\alp_{\dtheta'}(N^{\op}_{\gamma_{\theta}})$ and $C=\beta_{\theta'}(N^{\op}_{\gamma_{\theta}})$. Consider the maps
\[
\begin{array}{lccc}
(\mathcal{L} \ltak{}{\bullet} \mathcal{L}): & A \ltak{}{C} A & \to & A^{\op}_{\lambda,(\hat\gamma_{\theta} \odo S^{2}_{\G})} \ltak{}{C^{\op}_{\lambda}} A^{\op}_{\lambda,(\hat\gamma_{\theta} \odo S^{2}_{\G})} \\
& a' \ltak{}{C} b' & \mapsto & \mathcal{L}(a) \ltak{}{C^{\op}_{\lambda}} \mathcal{L}(b')
\end{array}
\]
\[
\begin{array}{lccc}
(\mathcal{L} \rtak{}{\bullet} \mathcal{L}): & A \rtak{}{B} A & \to & A^{\op}_{\lambda,(\hat\gamma_{\theta} \odo S^{2}_{\G})} \rtak{}{B^{\op}_{\lambda}} A^{\op}_{\lambda,(\hat\gamma_{\theta} \odo S^{2}_{\G})} \\
& a' \rtak{}{B} b' & \mapsto & \mathcal{L}(a) \rtak{}{B^{\op}_{\lambda}} \mathcal{L}(b')
\end{array}
\]

\begin{prop}\label{prop:iso2}
It hold
\begin{enumerate}[label=\textup{(\roman*)}]
\item $(\mathcal{L} \ltak{}{\bullet} \mathcal{L})\circ\com_{C} = \com^{\co}_{C^{\op}_{\lambda}}\circ\mathcal{L}$;
\item $(\mathcal{L} \rtak{}{\bullet} \mathcal{L})\circ\com_{B} = \com^{\co}_{B^{\op}_{\lambda}}\circ\mathcal{L}$;
\item $S_{\mathcal{A}^{\op,\co}_{\lambda}}\circ\mathcal{L} = \mathcal{L}\circ S_{\mathcal{A}}.$
\end{enumerate}
\end{prop}
\pr
Indeed, we have
\begin{align*}
(\mathcal{L} \ltak{}{\bullet} \mathcal{L})(\com_{C}(h \sm m^{\op})) & = (\mathcal{L} \ltak{}{\bullet} \mathcal{L})((h_{(1)} \sm 1^{\op}_{N}) \ltak{}{C} (h_{(2)} \sm m^{\op})) \\
& = (1_{N} \sm S_{\G}(h)_{(2)}) \ltak{}{C^{\op}_{\lambda}}(\hat\gamma_{\theta}(m) \sm S_{\G}(h)_{(1)}) \\
& = \com^{\co}_{C^{\op}_{\lambda}}(\hat\gamma_{\theta}(m) \sm S_{\G}(h)) \\
& = \com^{\co}_{C^{\op}_{\lambda}}(\mathcal{L}(h \sm m^{\op})),
\end{align*}
\begin{align*}
(\mathcal{L} \rtak{}{\bullet} \mathcal{L})(\com_{B}(h \sm m^{\op})) & = (\mathcal{L} \rtak{}{\bullet} \mathcal{L})((h_{(1)} \sm 1^{\op}_{N}) \ltak{}{B} (h_{(2)} \sm m^{\op})) \\
& = (1_{N} \sm S_{\G}(h)_{(2)}) \rtak{}{B^{\op}_{\lambda}}(\hat\gamma_{\theta}(m) \sm S_{\G}(h)_{(1)}) \\
& = \com^{\co}_{B^{\op}_{\lambda}}(\hat\gamma_{\theta}(m) \sm S_{\G}(h)) \\
& = \com^{\co}_{B^{\op}_{\lambda}}(\mathcal{L}(h \sm m^{\op}))
\end{align*}
and
\begin{align*}
S_{\mathcal{A}^{\op,\co}_{\lambda}}(\mathcal{L}(h \sm m^{\op})) & = S_{\mathcal{A}^{\op,\co}}((\hat\gamma_{\theta}(m) \sm S_{\G}(h))^{\op}) \\
& = (((1_{N} \sm S^{2}_{\G}(h))\beta_{\theta}(\hat\gamma^{2}_{\theta}(m)))^{\op} \\
& = \beta_{\theta}(\hat\gamma^{2}_{\theta}(m))^{\op}(1_{N} \sm S^{2}_{\G}(h))^{\op} \\
& = \mathcal{L}(\beta_{\theta'}(\hat\gamma_{\theta}(m)^{\op})(S^{2}_{\G}(h) \sm 1^{\op}_{N})) \\
& = \mathcal{L}(S_{\mathcal{A}}(h \sm m^{\op}))
\end{align*}
for all $h \in \pG$ and $m \in N$.
\fin

\begin{theo}
The \as-isomorphism $\mathcal{L}$ yields an isomorphism between the measured multiplier Hopf \as-algebroids $\mathcal{A}(N^{\op}_{\gamma_{\theta}},\theta',\dtheta',\mu^{\ops})$ and $\mathcal{A}^{\op,\co}_{\lambda}(N,\theta,\dtheta,\mu)$.
\end{theo}
\pr
Follows directly from Proposition~\ref{prop:iso1} and Proposition~\ref{prop:iso2}.
\fin

\begin{exam}
Consider the unital braided commutative measured (left) Yetter--Drinfeld $\Gc$-\as-algebra $(\C,\tr,\du{\tr},\id_{\C})$. We have
\[
\mathcal{A}_{\lambda}(\C,\tr,\du{\tr},\id_{\C}) = \mathcal{A}(\G), \quad\quad \du{\mathcal{A}}_{\lambda}(\C,\du{\tr},\id_{\C}) = \mathcal{A}(\dGo)
\]
and
\[
\mathcal{A}(\C,\tr,\du{\tr},\id_{\C}) = \mathcal{A}(\Gco), \quad\quad \du{\mathcal{A}}(\C,\du{\tr},\tr,\id_{\C}) = \mathcal{A}(\dGc).
\]
\end{exam}


\bibliographystyle{abbrv}
\bibliography{biblio}

\addresseshere


\end{document}